\documentclass[12pt]{article}
\usepackage{amsmath, amssymb, amsthm}
\numberwithin{equation}{section}
\usepackage{geometry}
\usepackage{xcolor}
\usepackage{graphicx}
\usepackage{hyperref}
\usepackage{amssymb, amsthm,amsfonts, amsmath,color}
\usepackage{amstext}
\newtheorem{theorem}{Theorem}[section]
\newtheorem{lemma}[theorem]{Lemma}
\newtheorem{proposition}[theorem]{Proposition}
\newtheorem{corollary}[theorem]{Corollary}
\newtheorem{remark}[theorem]{Remark}

\newcommand{\RR}{\mathbb{R}}

\newcommand{\eps}{\varepsilon}

\definecolor{darkblue}{rgb}{0,0,0.75}
\definecolor{darkred}{rgb}{0.75,0,0}

\begin{document}

\title{Improved Singular Adams-type inequalities and the existence of extremals in dimension four}

\author{ Xumin Wang,  Maochun Zhu* }
\date{}\maketitle

\makeatletter
\renewcommand{\@makefntext}[1]{%
  \parindent 1em \noindent  #1
}
\makeatother
\footnotetext{
The first  author  was supported by  the National Natural Science Foundation of China (No.12101050). The second author was supported by Natural Science Foundation of China (No.12571122) and the Fundamental research funds for the Central Universities (No.30926WZCT202601009).}
\makeatother

\begin{abstract}
Let $H_0^2(\Omega)$ denote the Sobolev space on $\Omega$, where $\Omega\subset\mathbb{R}^4$ is a smooth bounded domain containing the origin, and let $\lambda(\Omega)$ denote the first eigenvalue of the bi-Laplacian $\Delta^2$ on $\Omega$ subject to homogeneous Dirichlet boundary conditions.  We establish the following singular Adams inequality of Adimurthi-Druet type: for $0<\beta<4$ and $0\leq\alpha<\lambda(\Omega)$, 
\[
\sup\limits_{{u\in H_0^2(\Omega),\|\Delta u\|_2^2=1}}\int_\Omega \frac{e^{32\pi^2(1-\frac{\beta}{4}) u^2(1+\alpha\|u\|_2^2)}}{|x|^\beta}dx < +\infty.
\]
Moreover, this supremum is infinite for $\alpha>\lambda_1(\Omega)$.  Furthermore, we show the existence of extremal functions, that is, for sufficiently small $\alpha>0$, there exists $u^*\in H_0^2(\Omega)\cap C_{\mathrm{loc}}^4(\Omega\setminus{0})$ with $\|\Delta u^*\|_2^2=1$ such that
\[
\int_\Omega \frac{e^{32\pi^2(1-\frac{\beta}{4}) {|u^*|}^2(1+\alpha\|u^*\|_2^2)}}{|x|^\beta}dx=\sup\limits_{{u\in H_0^2(\Omega),\|\Delta u\|_2^2=1}}\int_\Omega \frac{e^{32\pi^2(1-\frac{\beta}{4}) u^2(1+\alpha\|u\|_2^2)}}{|x|^\beta}dx.
\]

These results can be regarded as a further development, in the singular setting, of the theory on the existence of extremals established by Lu and Yang [Adv. Math. (2009)]. Our proof relies on the blow-up analysis of the subcritical maximizers, combined with biharmonic truncation adapted to the singular weight, capacity-type estimates on annular regions  and a test function argument.

 A key novelty of our approach lies in carrying out the blow-up analysis without relying on an explicit expression for the limiting profile near the blow-up point, which is unavailable in the present setting. To the best of our knowledge, this is the first attempt to perform such a  blow-up analysis  without the explicit form of  the limiting profile.



\end{abstract}

\section{Introduction and main results}

Let $1\leq m<n$, and let $\Omega\subset\mathbb{R}^n$, $n\geq2$,
be a bounded domain. For $1<p<\infty$, let $W_0^{m,p}(\Omega)$
denote the closure of $C_c^\infty(\Omega)$ in $W^{m,p}(\Omega)$
with respect to the standard Sobolev norm
\[
\|u\|_{W^{m,p}(\Omega)}
=
\left(\sum_{|\gamma|\leq m}
\|D^\gamma u\|_{L^p(\Omega)}^p\right)^{1/p}.
\]
For $1<p<n/m$, the classical Sobolev embedding theorem gives
\[
W_0^{m,p}(\Omega)\hookrightarrow L^{p^*}(\Omega), \text{ with }
 p^*=\frac{np}{n-mp}.
\]
As $p\rightarrow n/m$, the exponent $p^*$ tends to infinity. At the
critical exponent $p=n/m$, one has
\[
W_0^{m,n/m}(\Omega)\hookrightarrow L^q(\Omega)
\quad\text{for every }1\leq q<\infty,
\]
but the embedding into $L^\infty(\Omega)$ fails.

\subsection{Trudinger-Moser and Adams inequalities and the existence of extremals}

It is well known that the endpoint counterpart of the Sobolev embedding
for $W_0^{m,\frac{n}{m}}(\Omega)$ is the Trudinger--Moser inequality
when $m=1$ (\cite{Tr1967,Moser1971}) and the Adams inequality when
$1<m<n$ (\cite{Adams1988}). These inequalities can
be stated in the following unified form:
\begin{equation}\label{eqad}
\sup_{\substack{u\in W_0^{m,\frac{n}{m}}(\Omega)\\
\|\nabla^m u\|_{L^{\frac{n}{m}}(\Omega)}\leq1}}
\frac{1}{|\Omega|}
\int_\Omega
\exp\left(\alpha |u(x)|^{\frac{n}{n-m}}\right)\,dx
<\infty
\text{ iff } 
\alpha\leq\alpha(n,m),
\end{equation}
where
\[
\nabla^m u=
\begin{cases}
\Delta^{\frac{m}{2}}u, & \text{if $m$ is even},\\
\nabla\Delta^{\frac{m-1}{2}}u, & \text{if $m$ is odd},
\end{cases}
\]
and
\[
\alpha(n,m)=
\begin{cases}
\displaystyle
\frac{n}{\omega_{n-1}}
\left[
\frac{\pi^{n/2}2^m\Gamma\left(\frac{m+1}{2}\right)}
{\Gamma\left(\frac{n-m+1}{2}\right)}
\right]^{\frac{n}{n-m}},
& \text{if $m$ is odd},\\[1.2em]
\displaystyle
\frac{n}{\omega_{n-1}}
\left[
\frac{\pi^{n/2}2^m\Gamma\left(\frac{m}{2}\right)}
{\Gamma\left(\frac{n-m}{2}\right)}
\right]^{\frac{n}{n-m}},
& \text{if $m$ is even}.
\end{cases}
\]

Due to a wide range of applications in geometric analysis and partial differential equations 
(see \cite{Cp2003,Figueiredo,LamLu20121,LamLu20131}  and references therein), so far,
the Trudinger-Moser and Adams inequalities have been extended
to a wide range of settings, including unbounded domains, compact and
noncompact Riemannian manifolds, hyperbolic spaces, and Heisenberg
groups, etc., to just name a few from along list of extensive works we refer the interested readers to \cite{AT2020,Ruf2005, RS2013, LR2008, LamLu2013, LamLu2012JDE,LT2016, YLi-JPDE01,YLi-SciChina05,  FT1993, ManciniSandeep-CCM2010,Lutang, LT2016,  LuYang-AIM2017,LuYang-AIM2022b, CohnLu2001,LL2012Adv,FT2018},  and the reference therein.
 
\bigskip

A natural and interesting question  associated with Trudinger-Moser-Adams inequalities is whether there exist extremal functions for these kind of  inequalities. The first result was due to the  work of Carleson and Chang \cite{CC1986}, who proved the existence of extremal functions  when $m=1$, $\Omega$ is a ball in $\mathbb{R}^n$,  by the symmetrization argument. Later, Carleson and Chang's results were extended by Flucher \cite{Fl1992} to arbitrary bounded domain in $\mathbb{R}^2$ and by Lin \cite{Lin1996} to general bounded domains in $\mathbb{R}^n$ using $n$-harmonic transplantation. 

In 2004, for a bounded domain $\Omega\subset\mathbb{R}^2$, Adimurthi and Druet \cite{AD2004} established an improved Trudinger--Moser inequality involving the $L^2$ norm: \begin{equation*} 
\underset{u\in H_{0}%
^{1}\left(  \Omega\right),\left\Vert \nabla u\right\Vert _{2}\leq1}{\sup}\int_{\Omega}\exp\big({  4\pi\left\vert u\right\vert ^{2}(
1+\alpha\left\Vert u\right\Vert _{2}^{2})}\big)
dx<\infty, \text{ iff }0\leq \alpha<\inf\limits_{{u\in H_0^{1}(\Omega), \|u\|_2^2=1}} \|\nabla u\|_2^2.
\end{equation*}
They also proved that the supremum is attained for sufficiently
small $\alpha$.  
Their proof is based on  blow-up analysis with Lions' concentration--compactness
principle for the Trudinger-Moser inequality (see \cite{Lions1985}). Extensions
involving $L^p$ norms were subsequently obtained in two dimensions
(see \cite{LY2009,Zhu2014}) and in higher dimensions (see \cite{Yang2006}).

In the higher order case $m>1$, the first result of the existence of extremals for the inequality \eqref{eqad} on bounded domain  was obtained by Lu and Yang  in \cite{LY-2009}. Indeed, in \cite{LY-2009}, they also obtained an improved Adams type inequalities  the spirit of Adimurthi and Druet: 
\begin{equation*} 
\sup\limits_{{u\in H_0^2(\Omega),\|\Delta u\|_2^2=1}}\int_\Omega \exp\bigl({32\pi^2 u^2(1+\alpha\|u\|_2^2)}\bigl)dx < +\infty \text{ iff } 0\leq\alpha< \lambda(\Omega),
\end{equation*}
where \begin{equation}\label{evalu1025}
    \lambda(\Omega):=\inf\limits_{{u\in H_0^2(\Omega), \|u\|_2^2=1}} \|\Delta u\|_2^2.
\end{equation} Furthermore, for sufficiently small $\alpha$, the extremal functions exist. Recently, DelaTorre and Mancini 
\cite{DG2021} extended the results of \cite{LY-2009} to arbitrary even dimension. We note that the work of Carleson and Chang was based on the rearrangement argument to reduce the problem to the one-dimensional case. However, the symmetrization technique cannot be used for the Adams inequality, since there is no corresponding P\'olya--Szeg\H{o} type inequality in the higher order case $m>1$. Instead, the authors
of  \cite{LY-2009,DG2021} applied  blow-up analysis, combined with the capacity-type estimates and the Pohozaev type identity to  rule out concentration and establish the existence of extremals. 
For further developments concerning the existence of extremals for Trudinger-Moser-Adams inequalities on unbounded domains and compact Riemannian manifolds, we refer the reader to \cite{Ruf2005,DSM2014,LZ2019JDE,YLi-SciChina05,LR2008,ChLuZhu2020} and the references therein.

 \subsection{Singular type of Trudinger-Moser-Adams inequalities and  the existence of extremals}

 Motivated by the Hardy type inequality in the critical dimension (see \cite{ACR2002}):  for any 
 $u \in W_0^{1,n}(\Omega)$ with $0 \in \Omega\subset\mathbb{R}^n$,
\[
\int_{\Omega} |\nabla u|^n  dx \geq \left(\frac{n-1}{n}\right)^n\int_{\Omega}\frac{|u|^n}{|x|^n \left(\log \frac{R}{|x|}\right)^n}dx,
\]
where $R\geq e \sup_{x\in\Omega} |x|$. Adimurthi and Sandeep proved in \cite{AS2007} a singular version Trudinger-Moser inequality: for any $\alpha<n$,
\begin{equation}\label{singu1025}
    \sup_{u\in W^{1,n}_0(\Omega),\int_{\Omega}|\nabla u|^n\,dx \leq 1}\int_{\Omega}
\frac{\exp\left(\beta |u|^{\frac{n}{n-1}}\right)}
{|x|^{\alpha}}dx
\leq +\infty \text{ iff } 0 \leq \beta \leq \left(1-\frac{\alpha}{n}\right)\alpha({n,1}).
\end{equation}
The existence of extremals for the singular Trudinger--Moser inequality \eqref{singu1025} was established Csat\'{o} and Roy in~\cite{Roy2015,CNR2021} by using the harmonic transplantation. Similar results were also established through the blowing up analysis method (see \cite{Nguyen-PA20,YangZhu2017,Lula}). We also mention a recent existence result for anisotropic singular Trudinger--Moser inequalities, established in \cite{LSXZ2024} through blow-up analysis and a weighted anisotropic isoperimetric inequality.

  The results discussed above concern the first-order setting. A higher-order
counterpart of \eqref{singu1025} was established   by Lam and Lu \cite{LamLu2012} in the critical
Sobolev space $W_0^{m,\frac{n}{m}}(\Omega)$ by the Riesz potential estimates combined with O'Neil's convolution rearrangement inequality \cite{ONeil1963}.
Their result can be stated as follows: for $0<\beta<n$, 
\begin{equation}\label{eqsadm}
\sup_{u \in W_0^{m,\frac{n}{m}}(\Omega),\; \|\nabla^m u\|_{\frac{n}{m}} \leq 1}
\int_\Omega
\frac{\exp\left(\alpha(n,m)\left(1-\frac{\beta}{n}\right)
|u|^{\frac{n}{n-m}}\right)}{|x|^\beta}dx
< +\infty.
\end{equation}

Now, an interesting and intriguing question is whether or not the  singular  Adams inequality \eqref{eqsadm} has an extremal. In this paper, we  devote ourselves  to studying this problem in the case $n=4$ and $m=2$, and give an affirmative answer to this question.

\subsection{Main results in this paper}

To state our main results, let $\Omega\subset\mathbb{R}^4$ be a smooth  oriented  bounded domain with $0\in\Omega$. Let $H_0^2(\Omega)$ be the completion of $C_c^\infty(\Omega)$ with respect to the Dirichlet norm
$\|u\|_{H_0^2(\Omega)}=\|\Delta u\|_2,$
where $\|\cdot\|_2$ denotes the standard $L^2(\Omega)$-norm, and $\lambda(\Omega)$ be defined by \eqref{evalu1025}.  
 
We first establish a refinement of \eqref{eqsadm} involving the $L^2$-norm.
\begin{theorem}\label{thm1}
Let
$$q(r)=1+a_1r+a_2r^2+\cdots+a_k r^k\ \ \ (k\ge 1)$$
be a polynomial of order $k$ in $\RR$. Suppose $0<\beta<4$ and
\[
0\le a_1<\lambda(\Omega),\quad 0\le a_2\le \lambda(\Omega)a_1,\quad \ldots,\quad 0\le a_k\le \lambda(\Omega)a_{k-1},
\]
then there holds
\begin{equation}\label{eqsa1}
S_{\Omega,q,\beta}:=\sup\limits_{u\in H_0^2(\Omega),\|\Delta u\|_2^2=1}\int_\Omega \frac{\exp\Big({32\pi^2(1-\frac{\beta}{4}) u^2 q(\|u\|_2^2)}\Big)}{|x|^{\beta}}dx < +\infty.
\end{equation} 
\end{theorem}
The condition $a_1<\lambda(\Omega)$ is almost sharp for \eqref{eqsa1} in the following sense.
\begin{theorem}[Divergence criteria]\label{Diver crit}
$S_{\Omega,q,\beta}=+\infty$
if either
\begin{enumerate}
    \item[(i)] $a_{1}>\lambda(\Omega)$; or
    \item[(ii)] $a_{1}=\lambda(\Omega)$ and one of the following holds:
    \begin{enumerate}
        \item[(a)] some first eigenfunction $u_0$ of the Dirichlet bi-Laplacian
        satisfies $u_0(0)\neq 0$;
        
        \item[(b)] every first eigenfunction satisfies $u_0(0)=0$, and there exists
        $ s_0\in\left(0,\frac{1}{\lambda(\Omega)}\right)$
        such that
        \begin{equation}
        \bigl(1-\lambda(\Omega)s_0\bigr)q(s_0)>1.
        \label{eq:condition-div}
        \end{equation}
    \end{enumerate}
\end{enumerate}
Moreover, in case $\mathrm{(ii)(b)}$, condition
\eqref{eq:condition-div} is satisfied whenever
$a_2>\lambda^2(\Omega).$
\end{theorem}

\begin{remark}
 It is worth emphasizing that the singular case exhibits a phenomenon that is fundamentally different from its unweighted case of Lu and Yang \cite{LY-2009}. In the unweighted case, the concentrating bubble can be centered at an arbitrary interior point. Hence, if a first eigenfunction $u_0$ vanishes at the origin, one may simply choose another point $x_0\in\Omega$ with $u_0(x_0)\neq0$  and construct the test function near $x_0$. This relocation preserves the best constant and yields, in the borderline case $a_1=\lambda(\Omega)$, a positive interaction term of order $t_\varepsilon\sqrt{\log(1/\varepsilon)}$, which forces the corresponding functional to diverge. In contrast, in the singular case, the best constant $32\pi^2(1-\beta/4)$ is tied to concentration at the singular point $x=0$. Moving the bubble to a fixed point $x_0\neq0$ makes the weight locally bounded, while the coefficient $32\pi^2(1-\beta/4)$ becomes strictly subcritical relative to the usual best constant $32\pi^2$ of Adams inequality. Consequently, the concentration point cannot be freely relocated. 

In the singular case, the analysis at the endpoint $a_1=\lambda(\Omega)$ is more delicate. Whether the integral diverges depends not only on the value of $\lambda(\Omega)$, but also on the behavior at the singular point $x=0$ of the eigenfunctions associated with $\lambda(\Omega)$. In particular, if there exists $$\phi\in E_1: = \lbrace \phi\in H_0^2(\Omega):\Delta^2\phi=\lambda(\Omega)\phi \rbrace$$ with $\phi(0)\neq0$, the above construction still yields divergence. In contrast, if every  $\phi\in E_1$ vanishes at the origin, the contribution from the first eigenspace no longer produces the same divergence, and the above construction fails. In this case, whether the integral diverges depends on the higher-order terms of $q(r)$ and the local behavior of the first eigenspace near the singular point.
\end{remark}

Our third result concerns the existence of extremals for the inequality \eqref{eqsa1}.

\begin{theorem}\label{thm2}
Under the assumptions and notation of Theorem~\ref{thm1}, there exists
a constant $\eps_0=\eps_0(\Omega,\beta,k)\in(0,\lambda(\Omega))$ such that if
\[
0\le a_1\le \eps_0,\quad 0\le a_2\le \lambda(\Omega)a_1,\quad \ldots,\quad 0\le a_k\le \lambda(\Omega)a_{k-1},
\]
then the supremum $S_{\Omega,\beta,q}$ is attained. More precisely,
there exists a function $u^*\in H_0^2(\Omega)$ satisfying $\|\Delta u^*\|_2=1$ such that
\[
S_{\Omega,q,\beta}
=\int_\Omega\frac{\exp\Big(32\pi^2\left(1-\frac{\beta}{4}\right)(u^*)^2q(\|u^*\|_2^2)\Big)}{|x|^\beta}dx.
\]
Moreover, $u^*\in C_{\mathrm{loc}}^4(\Omega\setminus\{0\})$.
\end{theorem}
As an immediate consequence of Theorem~\ref{thm1}, by taking the linear
polynomial $q(r)=1+\alpha r$, we obtain the singular Adams inequality of Adimurthi-Druet type.

\begin{corollary} 
Let $0<\beta<4$. Then, for every $0\leq\alpha<\lambda(\Omega)$, there
holds
\begin{equation}\label{eqsaad}
\sup\limits_{{u\in H_0^2(\Omega), \|\Delta u\|_2\leq1}}
\int_\Omega\frac{\exp\Big(32\pi^2\left(1-\frac{\beta}{4}\right)u^2\left(1+\alpha\|u\|_2^2\right)\Big)}{|x|^\beta}dx
<+\infty.
\end{equation}
 Moreover, the restriction $\alpha<\lambda(\Omega)$ is sharp in the sense of Theorem \ref{Diver crit}.
\end{corollary}
As a further consequence of Theorem~\ref{thm2}, by taking $q(r)=1+\alpha r$, we obtain the existence of extremal functions for inequality \eqref{eqsaad}.
\begin{corollary} 
Let $0<\beta<4$. Then, for sufficiently small $\alpha\geq 0$, there exists a function $u^*\in H_0^2(\Omega)$ satisfying $\|\Delta u^*\|_2=1$ such that
\begin{equation*} 
S_{\Omega,\alpha,\beta}=
\int_\Omega\frac{\exp\Big(32\pi^2\left(1-\frac{\beta}{4}\right)(u^*)^2\left(1+\alpha\|u^*\|_2^2\right)\Big)}{|x|^\beta}dx.
\end{equation*}
Moreover, $u^*\in C_{\mathrm{loc}}^4(\Omega\setminus\{0\})$.
\end{corollary}

 The general strategy we use for Theorems \ref{thm1} and \ref{thm2} is exploiting the blow-up analysis. We first prove the subcritical singular Adams inequalities and the existence of extremals through a weighted version of the Lions-type concentration-compactness principle. Then, we find a maximizing sequence $\{u_\eps\}\subseteq H^2_0(\Omega)$ for the supremum $S_{\Omega,q,\beta}$. If $u_\eps$  is bounded in $L^{\infty}(\Omega)$,  i.e., $c_\eps:=u_\eps(x_\eps)
=\max_{\overline{\Omega}}u_\eps
<C
 \text{ as }\eps\to0,$ standard elliptic estimates and compactness arguments yield, up to a subsequence, the strong convergence of $u_\varepsilon$ to some $u_0$ in $H^2_0(\Omega)$. While, if $c_\eps\rightarrow \infty$, i.e., the
 blow-up arises. In this case, we employ blow-up analysis to investigate the asymptotic behavior of $u_\varepsilon$ both near and away from the blow-up point $0\in\Omega$. This analysis leads to the following sharp upper bound under the blow-up assumption:

\begin{equation}\label{sharpupp1028}
\begin{aligned}
S_{\Omega,q,\beta}\leq\int_\Omega\frac{1}{|x|^{\beta}}dx+\exp\Big(32\pi^2\big(1-\frac \beta4\big)\left(A_0-C_0-J_0\right)\Big),
\end{aligned}
\end{equation}
where $A_0$ denotes the value at the origin of the trace of the regular part of the Green function associated with the biharmonic operator $\Delta^2$, and
\begin{equation*}\begin{cases}
\displaystyle
C_0=\frac{1}{8\pi^2}\int_{\RR^4}
\frac{\log|y|e^{64\pi^2(1-\frac{\beta}{4})\phi(y)}}{|y|^\beta}dy,\\
\displaystyle
J_0=\int_{\mathbb R^4}\frac{\phi(x)e^{64\pi^2(1-\frac \beta4)\phi(x)}}{|x|^{\beta}}dx,
\end{cases}\end{equation*}
here $\phi$ is the solution to \begin{equation}\label{liuw1028}
    \Delta^2\phi=\frac{\exp\bigl(64\pi^2(1-\frac{\beta}{4})\phi\bigr)}{|x|^\beta}\quad\text{in }\RR^4,
\end{equation}
under the  finite growth condition \begin{equation}\label{eqlimitDe1028}
\int_{B_R}|\Delta\phi|dx\leq CR^2.
\end{equation} Finally, we construct a function sequence in $H^2_0(\Omega)$ to show that the upper bound above can actually be surpassed.  This contradicts the upper bound derived under the blow-up assumption, thereby ruling out the possibility of concentration and establishing the existence of an extremal.

Neither the blow-up strategy in the study of second-order Adams inequality without singular weight (see  \cite{LY-2009}) nor the blow-up method for the first-order singular Trudinger-Moser inequality and existence of their extremals (see \cite{LY-2018JDE}, \cite{Nguyen-PA20,YangZhu2017}) can be easily generalized to second-order singular Adams inequality case.  In what follows, we highlight the key differences between the blow-up analysis for the second-order singular Adams inequality and the corresponding analyses for the second-order Adams and singular Trudinger-Moser inequalities. We will also outline the key novel aspects of our approach in the present setting.

To derive the upper bound in \eqref{sharpupp1028} for the singualr Adams functional  when $c_\eps\rightarrow+\infty$, one need to analyze the asymptotic behavior of $u_\eps$ near and far away from the blow up point $0$. For the asymptotic behavior of $u_\eps$ near $0$, a crucial step is to classify the solutions to the limit equation \eqref{liuw1028}. Fortunately, the classification result for \eqref{liuw1028} has been  obtained by Ahmedou, Wu and Zhang in \cite{AWZ}  under the assumption of subquadratic growth at infinity. In our setting, to verify the desired subquadratic growth, a key step is to  establish the finite-growth condition \eqref{eqlimitDe1028}. This relies on the following local estimate for $\Delta u_\eps$:
\[
c_\eps \int_{B_{Rr_\eps^{\frac{4}{4-\beta}}}} |\Delta u_\eps| dx\leq C \left(Rr_\eps^{\frac{4}{4-\beta}}\right)^2,
\]
which is obtained by combining an $L(\log L)^{1/2}$-type estimate for $\Delta^2 u_\eps$ with Lorentz-space estimates for higher-order derivatives of $u_\eps$.  Due to the presence of singular weight, the singular behavior near the origin becomes more pronounced, which making the local $L(\log L)^{1/2}$ estimate more delicate. 

To investigate the asymptotic behavior of $u_\eps$ away from the blow-up point, we aim to adapt the truncation approach developed for the singular Trudinger-Moser inequality in \cite{Nguyen-PA20} and \cite{YangZhu2017}. However, in the second-order Sobolev space $H_0^2(\Omega)$, standard truncation techniques cannot be directly applied, as they generally fail to preserve the required higher-order regularity. To overcome this difficulty, we adapt the biharmonic truncation method introduced by DelaTorre and Mancini in \cite{DG2021} to our singular setting. A key step is to establish the pointwise gradient estimates for $u_\eps$:
 for each $l\in\{1,2,3\}$, 
\begin{equation*} 
\min\{|x|,|x-x_\eps|\}^{l}|u_\eps(x)|\,|\nabla^l u_\eps(x)|\le C \text{ for all } x\in\Omega\setminus\{0,x_\eps\}.
\end{equation*}
 which allow us to construct biharmonic truncations which have many nice 
properties: on the one hand, they preserve the high order regularity; on the other hand, 
they behave approximately like constants within suitable balls.

 When we try to obtain the sharp upper bound  for the singular Adams inequality if the concentration phenomenon occur, an important step consists in finding the sharp lower bounds of the integral of $|\Delta u_\eps|^2$
 on some annular regions when we carry out the capacity estimates. In the earlier work of \cite{Nguyen-PA20,YangZhu2017}, this could be achieved by comparing the energy of $u_\eps$ with the quantity \[
\min_{u \in \{ u : u|_{\partial B_\delta(x_\eps)}=a, u|_{\partial B_{r_\eps^{\frac{n}{n-\beta}}}(x_\eps)}=b
 \}} \int_{B_{\delta}(x_\eps)\setminus B_{r_\eps^{\frac{n}{n-\beta}}}(x_\eps)} |\nabla u|^n  dx,
\]
whose extremal function is some $n$-harmonic function which can be explicitly obtained by solving some  equation. However, in the second-order setting, the explicit expression of the corresponding extremal appears to be rather complicated. In this work, we will compute the lower bound by directly comparing the Dirichlet energy of $u_\eps$ with some bi-harmonic function in the annular region, this idea was first used by  DelaTorre and Mancini~ in \cite{DG2021}. In previous works \cite{LY-2009,DG2021,Nguyen-PA20,YangZhu2017}, such computations rely heavily on the explicit form of the limiting profile near the blow-up point. However, unlike in the unweighted Adams inequality setting (see \cite{LY-2009,DG2021}) or the first-order singular Trudinger-Moser setting (see \cite{Nguyen-PA20,YangZhu2017}), the limiting profile in our second-order singular setting is not explicitly known (see \cite{AWZ}). This lack of an explicit representation introduces substantial difficulties in carrying out the corresponding comparison estimates. For instance, using comparing idea and the asymptotic behavior of $u_\eps$ near and far away from the blow up point $0$, one can only obtain \begin{equation*} 
\begin{aligned}
 S_{\Omega,q,\beta}
 \leq&\lim_{R\rightarrow +\infty}32\pi^2\left(1-\frac\beta4\right)
 \Bigg(A_0-\int_{B_R}|\Delta\phi|^2dx-2\phi(R)-R\phi'(R)\\
 &\ \
 -\frac1{8\pi^2}\log R-\frac1{16\pi^2}\Bigg).
\end{aligned}
\end{equation*}
To derive the sharp upper bound for \(S_{\Omega,q,\beta}\), we need to evaluate the limit
\begin{equation*}
\lim_{R\to +\infty}\int_{B_R}|\Delta\phi|^2dx+2\phi(R)+R\phi'(R)+\frac{1}{8\pi^2}\log R.
\end{equation*}
By combining the representation formula and the mass identity with a delicate analysis of the asymptotic behavior of the limiting profile $\phi(x)$ at the infinity, we obtain an integral representation of this limit, which in turn yields an explicit formula for the sharp upper bound \eqref{sharpupp1028}.

 Finally, unlike previous studies on the existence of extremals for the Trudinger-Moser-Adams inequality (see \cite{Roy2015}, \cite{YangZhu2017}, \cite{LY-2009} and  \cite{ChLuZhu2020} for example), where test functions were constructed based on the explicit limiting profile and the Green's function, the singular Adams inequality considered here lacks such an explicit expression. This absence poses significant difficulties in constructing test functions and performing the associated estimates. Nevertheless, by applying a suitable scaling to the limit function, we are able to construct the inner part of the test function. Together with the radial symmetry, the unit weighted mass, the asymptotic behavior of the limiting profile at infinity (\cite{AWZ}), energy estimates of the Green's function, and a careful  error analysis, we derive a precise energy estimate for the test functions. This enables us to prove that the sharp upper bound obtained when blow-up occurs is surpassed, which in turn implies the existence of extremal functions.

The remainder of this paper is organized as follows. In Section 2, we establish the sharpness of the singular Adams inequality by constructing suitable test functions. In Section 3, we prove the subcritical Adams inequalities and establish the existence of maximizers. In Section 4, we investigate the asymptotic behavior of these maximizers both near and away from the origin in the blow-up case, and introduce the notion of biharmonic truncation. In Section 5, we derive a sharp upper bound under the blow-up assumption via capacity estimates and complete the proof of Theorem~\ref{thm1}. In Section 6, we construct test functions whose functional values strictly exceed this upper bound, thereby ruling out the possibility of blow-up and establishing the existence of extremals asserted in Theorem~\ref{thm2}. Finally, in the Appendix, we establish pointwise gradient estimates for the subcritical maximizers and prove several derived  properties of biharmonic truncations.

\section{The sharpness of the singular Adams inequality}
Throughout the rest of this paper, $C$ will stand for a positive constant, may vary from line to line.

Suppose $0<\beta<4$, in this section, we will prove that  for $a_1> \lambda(\Omega)$ and arbitrary $a_2,\dots,a_k$, it holds
\begin{equation*} 
\sup\limits_{u\in H_0^2(\Omega),\|\Delta u\|_2^2=1}\int_\Omega \frac{\exp\Big({32\pi^2(1-\frac{\beta}{4}) u^2 q(\|u\|_2^2)}\Big)}{|x|^{\beta}}dx = +\infty.
\end{equation*}

Let $u_0\in H_0^2(\Omega)\cap C^4(\overline{\Omega})$ be an eigenfunction of bi-Laplacian operator $\Delta^2$ satisfying
\[
\begin{cases}
\Delta^2 u_0 = \lambda(\Omega) u_0 & \text{in } \Omega,\\
\|u_0\|_2=1,\quad u_0=\frac{\partial u_0}{\partial \nu}=0 & \text{on } \partial\Omega.
\end{cases}
\]
This equation is solvable by applying the direct method of variation. Without loss of generality we assume the unit ball $\mathbb{B}\subseteq \Omega$.

Define
\[
u_\eps(x)=
\begin{cases}
\sqrt{\frac{1}{32\pi^2}\log\frac{1}{\eps}} - \frac{|x|^2}{\sqrt{8\pi^2\eps\log\frac{1}{\eps}}} + \frac{1}{\sqrt{8\pi^2\log\frac{1}{\eps}}}, & |x|\le \eps^{1/4},\\[1.2ex]
\frac{1}{\sqrt{2\pi^2\log\frac{1}{\eps}}}\log\frac{1}{|x|}, & \eps^{1/4}<|x|\le 1,\\[1.2ex]
\xi_\eps, & |x|>1,
\end{cases}
\]
where $\xi_\eps$ is a smooth function satisfying $\xi_\eps|_{\partial \mathbb{B}}=\xi_\eps|_{\partial \Omega}=0$, $\frac{\partial \xi_\eps}{\partial \nu}|_{\partial \mathbb{B}} = \frac{1}{\sqrt{2\pi^2\log\frac{1}{\eps}}}$, $\frac{\partial \xi_\eps}{\partial \nu}|_{\partial \Omega}=0$ and $\xi_\eps, \nabla \xi_\eps, \Delta \xi_\eps$ are all $O(1/\sqrt{\log(1/\eps)})$. By direct computating,  one can easily check that $u_\eps\in H_0^2(\Omega)$, and
\[
\|u_\eps\|_2^2 = O\left(\frac{1}{\log(1/\eps)}\right),\qquad \|\Delta u_\eps\|_2^2 = 1+O\left(\frac{1}{\log(1/\eps)}\right).
\]

Let $t_\varepsilon \geq 0$ be such that $t_\varepsilon \to 0, 
t_\varepsilon^2 \log(1/\varepsilon) \to +\infty,
 \text{ and }  
t_\varepsilon^2 \bigl(\log(1/\varepsilon)\bigr)^{1/2} \to 0$, and define
$v_\varepsilon := u_\varepsilon + t_\varepsilon u_0.$ Then we have
\begin{equation}\label{add1-103}
    \begin{aligned}
\|v_\eps\|_2^2 &= \|u_\eps\|_2^2 + t_\eps^2\|u_0\|_2^2 + 2t_\eps\int_\Omega u_\eps u_0dx\\
&= t_\eps^2 + 2t_\eps\int_\Omega u_\eps u_0 ~dx + O\left(\frac{1}{\log(1/\eps)}\right)
\end{aligned}
\end{equation}
and
\begin{equation}\label{add2-103}
\begin{aligned}
\|\Delta v_\eps\|_2^2 &= \|\Delta u_\eps\|_2^2 + t_\eps^2\|\Delta u_0\|_2^2 + 2t_\eps\int_\Omega \Delta u_\eps \Delta u_0 ~dx\\
&= 1 + 2\lambda(\Omega)t_\eps\int_\Omega u_\eps u_0 ~dx + \lambda(\Omega)t_\eps^2 + O\left(\frac{1}{\log(1/\eps)}\right).
\end{aligned}
\end{equation}
A straightforward calculation shows
\begin{equation}\label{eqqu}
\frac{1}{\|\Delta v_\eps\|_2^2} q\left(\frac{\|v_\eps\|_2^2}{\|\Delta v_\eps\|_2^2}\right)
= 1 + (a_1-\lambda(\Omega))\left(t_\eps^2 + 2t_\eps\int_\Omega u_\eps u_0dx\right)
+ o\left(\frac{t_\eps}{(\log(1/\eps))^{1/2}}\right).
\end{equation}

Noting that $\int_\Omega |u_\eps|dx = O(1/\sqrt{\log(1/\eps)})$, $\int_\Omega u_\eps u_0dx=O\left(\frac{1}{\sqrt{\log(1/\eps)}}\right)$ and $t_\eps^2\log\frac{1}{\eps}\rightarrow+\infty,$ we have $t_\eps\sqrt{\log\frac{1}{\eps}}\rightarrow+\infty.$
Consequently, we get
\[
t_\eps\int_\Omega u_\eps u_0dx=O\left(\frac{t_\eps}{\sqrt{\log(1/\eps)}}\right)=o(t_\eps^2)
\]
and
\[
\frac{t_\eps}{\sqrt{\log(1/\eps)}}=o(t_\eps^2).
\]

Since on $\mathbb{B}_{\eps^{1/4}}$, 
\[
v_\epsilon=\sqrt{\frac{1}{32\pi^2}\log\frac{1}{\epsilon}}+t_\epsilon u_0(x)+o(1),
\]
we have
\begin{equation}\label{eqve}
v_\epsilon^2=\frac{1}{32\pi^2}\log\frac{1}{\epsilon}+2t_\epsilon u_0(x)\sqrt{\frac{1}{32\pi^2}\log\frac{1}{\epsilon}}
+t_\epsilon^2u_0^2(x)+o\left(t_\epsilon\sqrt{\log\frac{1}{\epsilon}}\right).
\end{equation}
Combining \eqref{add2-103}, \eqref{add1-103}, \eqref{eqqu} and \eqref{eqve}, for any  $x\in\mathbb{B}_{\eps^{1/4}}$, we have 
\begin{equation}\begin{split}\label{eqveq}
&32\pi^2\left(1-\frac{\beta}{4}\right)
\frac{v_\eps^2}{\|\Delta v_\eps\|_2^2}q\left(
\frac{\|v_\eps\|_2^2}{\|\Delta v_\eps\|_2^2}\right)\\
&=\left(1-\frac{\beta}{4}\right)\log\frac{1}{\eps}+\left(1-\frac{\beta}{4}\right)
\bigl(a_1-\lambda(\Omega)\bigr)t_\eps^2\log\frac{1}{\eps}\\
&\quad+8\sqrt{2}\pi\left(1-\frac{\beta}{4}\right)
u_0(x)t_\eps\sqrt{\log\frac{1}{\eps}}
+o\left(\bigl(a_1-\lambda(\Omega)\bigr)t_\eps^2\log\frac{1}{\eps}
+t_\eps\sqrt{\log\frac{1}{\eps}}\right).
\end{split}\end{equation}
Employing  \eqref{eqveq} and the fact that 
\begin{equation}\label{add310-3}
    \int_{B_{\eps^{1/4}}(0)}\frac{dx}{|x|^\beta}=\frac{2\pi^2}{4-\beta}\eps^{1-\frac{\beta}{4}},
\end{equation}
we obtain

\begin{equation}\label{eqia1}
\begin{aligned}
&\int_{\Omega}\frac{\exp\!\left\{32\pi^2\left(1-\frac{\beta}{4}\right)
\frac{v_\eps^2}{\|\Delta v_\eps\|_2^2}q\!\left(
\frac{\|v_\eps\|_2^2}{\|\Delta v_\eps\|_2^2}\right)\right\}}{|x|^\beta}dx\\
&\ge\int_{\mathbb B_{\eps^{1/4}}}
\frac{\eps^{-\left(1-\frac{\beta}{4}\right)}}{|x|^\beta}
\exp\!\Bigg\{\left(1-\frac{\beta}{4}\right)
\bigl(a_1-\lambda(\Omega)\bigr)t_\eps^2\log\frac{1}{\eps}\\
&\quad+8\sqrt{2}\pi\left(1-\frac{\beta}{4}\right)
u_0(x)t_\eps\sqrt{\log\frac{1}{\eps}}
+o\!\left(
\bigl(a_1-\lambda(\Omega)\bigr)t_\eps^2\log\frac{1}{\eps}+t_\eps\sqrt{\log\frac{1}{\eps}}
\right)\Bigg\}dx\\
&\ge\frac{2\pi^2}{4-\beta}
\exp\!\Bigg\{\left(1-\frac{\beta}{4}\right)\bigl(a_1-\lambda(\Omega)\bigr)t_\eps^2\log\frac{1}{\eps}\\
& \quad+8\sqrt{2}\pi\left(1-\frac{\beta}{4}\right)
\eps^{-\left(1-\frac{\beta}{4}\right)}
\left(\int_{\mathbb B_{\eps^{1/4}}}\frac{u_0(x)}{|x|^\beta}dx\right)
t_\eps\sqrt{\log\frac{1}{\eps}}\\
&\quad+o\!\left(
\bigl(a_1-\lambda(\Omega)\bigr)t_\eps^2\log\frac{1}{\eps}
+t_\eps\sqrt{\log\frac{1}{\eps}}
\right)
\Bigg\}.
\end{aligned}
\end{equation}

We now distinguish the following three cases:

Case 1. If $a_1>\lambda(\Omega)$: without loss of generality, we can assume $u_0(0)\geq 0$. Then
\[\bigl(a_1-\lambda(\Omega)\bigr)t_\eps^2\log\frac{1}{\eps}\rightarrow+\infty.\]
Therefore, the right-hand side of \eqref{eqia1} tends to $+\infty$, regardless of whether $u_0(0)$ vanishes.

Case 2. If $a_1=\lambda(\Omega)$ and there exists a first eigenfunction of the Dirichlet bi-Laplacian such that $u_0(0)>0$, then
\[u_0(0)t_\eps\sqrt{\log\frac{1}{\eps}}\rightarrow+\infty,\]
hence the right-hand side of \eqref{eqia1} again tends to $+\infty$.

Case 3. If $a_1=\lambda(\Omega)$ and every first eigenfunction satisfies $u_0(0)=0$,
assume that there exists
$s_{0}\in\left(0,1/\lambda(\Omega)\right)$ such that
$F(s_{0})=\bigl(1-\lambda(\Omega)s_{0}\bigr)q(s_{0})
>1.$
Choose a fixed constant $t>0$ satisfying
\begin{equation*} 
s_{0}=\frac{t^{2}}{1+\lambda(\Omega)t^{2}},
\end{equation*}
equivalently,
\[
t^{2}=\frac{s_{0}}{1-\lambda(\Omega)s_{0}}.
\]
Notice that
\begin{equation}\label{eq:identity-t}
\frac{1}{1+\lambda(\Omega)t^{2}}=1-\lambda(\Omega)s_{0}.
\end{equation}

Define
\begin{equation*} 
v_{\varepsilon}:=u_{\varepsilon}+t u_{0}.
\end{equation*}
Arguing as in \eqref{add2-103} and \eqref{add1-103}, we have
\begin{align}
\|\Delta v_{\varepsilon}\|_{2}^{2}
&=1 + 2\lambda(\Omega)t\int_\Omega u_\eps u_0dx + \lambda(\Omega)t^2 + O\left(\frac{1}{\log(1/\eps)}\right)
\nonumber\\
&=1+\lambda(\Omega)t^{2}+o(1)\label{eq:v-energy}
\end{align}
and
\begin{align}\label{eq:v-l2}
\|v_{\varepsilon}\|_{2}^{2}=\|u_{\varepsilon}\|_{2}^{2}+2t\int_{\Omega}u_{\varepsilon}u_{0}dx
+t^{2}\|u_{0}\|_{2}^{2}=t^{2}+o(1).
\end{align}
It follows from \eqref{eq:v-energy} and \eqref{eq:v-l2} that
\begin{equation}\label{eq:s-eps-limit}
\frac{\|v_{\varepsilon}\|_{2}^{2}}
{\|\Delta v_{\varepsilon}\|_{2}^{2}}
\rightarrow
\frac{t^{2}}{1+\lambda(\Omega)t^{2}}
=
s_{0},
\end{equation}
as $\varepsilon\rightarrow0$.
Since $q$ is continuous, as $\varepsilon\rightarrow0$, we get
\begin{align}\label{eq:effective-factor}
\frac{1}{\|\Delta v_{\varepsilon}\|_{2}^{2}}
q\!\left(\frac{\|v_{\varepsilon}\|_{2}^{2}}
{\|\Delta v_{\varepsilon}\|_{2}^{2}}\right)\rightarrow
\frac{q(s_{0})}{1+\lambda(\Omega)t^{2}}
=\bigl(1-\lambda(\Omega)s_{0}\bigr)q(s_{0})
=F(s_{0})>1, 
\end{align}
where in the second equality we have used \eqref{eq:identity-t}.

We next analyze $v_{\varepsilon}$ near $0$.
Since $u_{0}\in C^{1}(\overline{\Omega})$ and $u_{0}(0)=0$,
there exist $C>0$ and $r_{0}>0$ such that
\[
|u_{0}(x)|\leq C|x|,
\text{ for }x\in B_{r_{0}}(0).
\]
Hence, for $\varepsilon>0$ sufficiently small and
$x\in B_{\varepsilon^{1/4}}(0)$,
\begin{equation*} 
u_{0}(x)
=
O(\varepsilon^{1/4})
\end{equation*}
uniformly in $x$. Therefore, we obtain
\[
u_{\varepsilon}(x)u_{0}(x)
=
O\!\left(
\varepsilon^{1/4}
\sqrt{\log\frac{1}{\varepsilon}}
\right)
=
o(1)
\]
uniformly on $B_{\varepsilon^{1/4}}(0)$. Since $t$ is fixed number, 
$t^{2}u_{0}^{2}(x)=o(1)$ uniformly on the same ball. Thus,
\begin{align}\label{eq:v-core}
v^2_{\varepsilon}(x)
&=u^2_{\varepsilon}(x)+2t u_{\varepsilon}(x)u_{0}(x)+t^{2}u_{0}^{2}(x)\nonumber\\
&=u^2_{\varepsilon}(x)+o(1)
=\frac{1}{32\pi^{2}}\log\frac{1}{\varepsilon}+O(1)
\end{align}
uniformly for $x\in B_{\varepsilon^{1/4}}(0)$.

Now normalize $v_{\varepsilon}$ by setting
\begin{equation*} 
w_{\varepsilon}
:=
\frac{v_{\varepsilon}}
{\|\Delta v_{\varepsilon}\|_{2}},
\end{equation*}
and use \eqref{eq:s-eps-limit}, we have
\[
\|w_{\varepsilon}\|_{2}^{2}
=
\frac{\|v_{\varepsilon}\|_{2}^{2}}
{\|\Delta v_{\varepsilon}\|_{2}^{2}}
=
s_{0}+o(1).
\]
Combining \eqref{eq:effective-factor} and \eqref{eq:v-core}, 
  for $x\in B_{\varepsilon^{1/4}}(0)$, we obtain 
\begin{align*}
32\pi^{2}\left(1-\frac{\beta}{4}\right)w^2_{\varepsilon}(x)q(\|w_{\varepsilon}\|_{2}^{2})
&=32\pi^{2}\left(1-\frac{\beta}{4}\right)
v^2_{\varepsilon}(x)\frac{1}{\|\Delta v_{\varepsilon}\|_{2}^{2}}
q\!\left(\frac{\|v_{\varepsilon}\|_{2}^{2}}
{\|\Delta v_{\varepsilon}\|_{2}^{2}}\right)
\nonumber\\
&\quad=\left(1-\frac{\beta}{4}\right)F(s_{0})\log\frac{1}{\varepsilon}
+o\!\left(\log\frac{1}{\varepsilon}\right),
\end{align*}
as $\varepsilon\rightarrow0$.
This together with \eqref{add310-3} implies
\begin{align}
S_{\Omega,q,\beta}
&\geq
\int_{B_{\varepsilon^{1/4}}(0)}
\frac{
\exp\!\left(
32\pi^{2}\left(1-\frac{\beta}{4}\right)
w^2_{\varepsilon}
q(\|w_{\varepsilon}\|_{2}^{2})
\right)
}{|x|^{\beta}}
dx
\nonumber\\
&\geq
\exp\!\left(
\left(1-\frac{\beta}{4}\right)
F(s_{0})
\log\frac{1}{\varepsilon}
+
o\!\left(
\log\frac{1}{\varepsilon}
\right)
\right)
\int_{B_{\varepsilon^{1/4}}(0)}
|x|^{-\beta}dx\nonumber\\
&\geq
\frac{2\pi^{2}}{4-\beta}
\exp\!\Bigg(
\left(1-\frac{\beta}{4}\right)
F(s_{0})
\log\frac{1}{\varepsilon}
-
\left(1-\frac{\beta}{4}\right)
\log\frac{1}{\varepsilon}+
o\!\left(
\log\frac{1}{\varepsilon}
\right)
\Bigg)
\nonumber\\
&=
\frac{2\pi^{2}}{4-\beta}
\exp\!\Bigg(
\left(1-\frac{\beta}{4}\right)
\bigl(F(s_{0})-1\bigr)
\log\frac{1}{\varepsilon}
+
o\!\left(
\log\frac{1}{\varepsilon}
\right)
\Bigg).
\label{eq:divergence-final}
\end{align}

Since $1-\frac{\beta}{4}>0 \text{ and } 
F(s_{0})-1>0,$ the right-hand side of \eqref{eq:divergence-final} tends to
$+\infty$ as $\varepsilon\to0$. Hence
$S_{\Omega,q,\beta}=+\infty.$

Finally, it remains to verify the explicit sufficient condition $a_{2}>\lambda^{2}(\Omega)$ for \eqref{eq:condition-div}. Since
$q(s)=1+\lambda(\Omega)s+a_{2}s^{2}+\cdots+a_{k}s^{k},$ we have
\begin{align*}
F(s)&=\bigl(1-\lambda(\Omega)s\bigr)q(s)
\nonumber\\
&=1+\bigl(a_{2}-\lambda^{2}(\Omega)\bigr)s^{2}
+\bigl(a_{3}-\lambda(\Omega)a_{2}\bigr)s^{3}
+\cdots+\bigl(a_{k}-\lambda(\Omega)a_{k-1}\bigr)s^{k}
-\lambda(\Omega)a_{k}s^{k+1}.
\end{align*}
If $a_{2}>\lambda^{2}(\Omega),$ then, as $s\to0^{+}$, we have 
\[
F(s)=1+\bigl(a_{2}-\lambda^{2}(\Omega)\bigr)s^{2}+O(s^{3}).
\]
Because $a_{2}-\lambda^{2}(\Omega)>0,$ there exists a sufficiently small $s_{0}\in\left(0, 1/\lambda(\Omega)\right)$ such that $F(s_{0})>1,$ and the proof is finished.

\section{Extremals for the subcritical singular Adams inequality}\label{boundedcase1020}
In this section, we establish the existence of maximizers for the
subcritical  supremum
\begin{equation}\label{eqsubf}
S_{\Omega,q,\beta,\varepsilon}:=\sup\limits_{{u\in H_0^2(\Omega), \|\Delta u\|_2\leq1}}
\int_\Omega\frac{\exp\left(\left(32\pi^2\left(1-\frac{\beta}{4}\right)-\eps\right)u^2q(\|u\|_2^2)\right)}{|x|^\beta}dx,
\end{equation}
where  $0<\eps<32\pi^2\left(1-\frac{\beta}{4}\right).$

We begin by establishing a weighted version of the concentration--compactness principle originally due to Lions
\cite{Lions1985}.

\begin{proposition}\label{prop3.1}
Let $\{u_\eps\}_{\eps>0}\subset H_0^2(\Omega)$ be a sequence of functions such that $\|\Delta u_\eps\|_2=1$ and $u_\eps\rightharpoonup u_0$ weakly in $H_0^2(\Omega)$. Then for any $0<p<\frac{1}{1-\|\Delta u_0\|_2^2}$, we have
\[
\limsup_{\eps\to 0}\int_\Omega \frac{\exp\Big({32\pi^2(1-\frac{\beta}{4}) p u_\eps^2}\Big)}{|x|^{\beta}}dx < +\infty.
\]
\end{proposition}

\begin{proof}
If $u_0=0$, then for any $p<1$, and the conclusion follows directly from the singular Adams inequality \eqref{eqsadm} from Lam and Lu \cite{LamLu2012}.

Suppose that $u_0\neq0$. Since
$u_\eps\rightharpoonup u_0$ weakly in $H_0^2(\Omega)$ and
$\|\Delta u_\eps\|_2=1$, the following holds as $\eps\rightarrow 0$: 
\[
\begin{aligned}
\|\Delta(u_\eps-u_0)\|_2^2
&=
\|\Delta u_\eps\|_2^2+\|\Delta u_0\|_2^2
-2\int_\Omega\Delta u_\eps\Delta u_0dx\\
&\rightarrow
1-\|\Delta u_0\|_2^2<1.
\end{aligned}
\]
Since $p\left(1-\|\Delta u_0\|_2^2\right)<1,$ we can choose $\delta>0$ and $r>1$, with both $\delta$ and $r-1$
sufficiently small, such that
$rp(1+\delta)\left(1-\|\Delta u_0\|_2^2\right)<1.$ Hence, for these choices of $\delta$ and $r$, when $\varepsilon$ is sufficiently small, we have
\[
rp(1+\delta)\|\Delta(u_\eps-u_0)\|_2^2<1.
\]
Let $s>1$ satisfy $\frac1r+\frac1s=1.$
Using the elementary inequality
\[
u_\eps^2
\leq
(1+\delta)(u_\eps-u_0)^2
+\left(1+\frac1\delta\right)u_0^2
\]
and H\"{o}lder's inequality, we obtain
\[
\begin{aligned}
&\int_\Omega \frac{\exp\Big({32\pi^2(1-\frac{\beta}{4}) p u_\eps^2}\Big)}{|x|^{\beta}}dx\\
&\le \int_\Omega \frac{\exp\Big({32\pi^2(1-\frac{\beta}{4}) p(1+\delta)(u_\eps-u_0)^2 + 32\pi^2(1-\frac{\beta}{4}) p(1+1/\delta)u_0^2}\Big)}{|x|^{\beta}}dx \\
&\leq\left(\int_\Omega\frac{\exp\left(32\pi^2\left(1-\frac{\beta}{4}\right)rp(1+\delta)(u_\eps-u_0)^2\right)}{|x|^\beta}dx
\right)^{1/r}\\
&\ \ \times\left(\int_\Omega\frac{\exp\left(32\pi^2\left(1-\frac{\beta}{4}\right)sp\left(1+\frac1\delta\right)u_0^2\right)}{|x|^\beta}dx
\right)^{1/s}\\
&\le \left(\int_\Omega \frac{\exp\Big({32\pi^2(1-\frac{\beta}{4}) \frac{(u_\eps-u_0)^2}{\|\Delta(u_\eps-u_0)\|_2^2}}\Big)}{|x|^{\beta}}dx\right)^{1/r}
\left(\int_\Omega \frac{\exp\Big({32\pi^2(1-\frac{\beta}{4}) p' u_0^2}\Big)}{|x|^{\beta}}dx\right)^{1/s},
\end{aligned}
\]
where $p'=sp\left(1+\frac1\delta\right)>p.$

The singular Adams inequality \eqref{eqsadm} implies that the first term on the right-hand side of the above inequality is bounded independent of $\varepsilon$. For the second term,
since $0<\beta<4$, we may choose $\ell>1$ such that
$\beta\ell<4$. Let $\ell'=\ell/(\ell-1)$. Since
$u_0\in H_0^2(\Omega)$ is fixed, the Orlicz embedding implies that
$e^{C u_0^2}\in L^{\ell'}(\Omega)$ for every fixed $C>0$. Therefore,
by H\"older's inequality, we have
\[
\begin{aligned}
&\int_\Omega
\frac{
\exp\left(
32\pi^2\left(1-\frac{\beta}{4}\right)p'u_0^2
\right)}
{|x|^\beta}dx\\
&\leq
\left(\int_\Omega |x|^{-\beta\ell}dx\right)^{1/\ell}
\left(
\int_\Omega
\exp\left(
32\pi^2\left(1-\frac{\beta}{4}\right)
p'\ell' u_0^2
\right)dx
\right)^{1/\ell'}
<+\infty.
\end{aligned}
\]
Therefore, we derive
\[
\limsup_{\eps\to 0}\int_\Omega \frac{\exp\Big({32\pi^2(1-\frac{\beta}{4}) p u_\eps^2}\Big)}{|x|^{\beta}}dx < +\infty.
\]
\end{proof}

Next we prove the existence of maximizers for subcritical functionals defined in \eqref{eqsubf}.

\begin{proposition}\label{prop3.2}
Assume that $q(t)$ satisfies the assumptions of Theorem~\ref{thm1}. Then for any $0<\eps<32\pi^2\left(1-\frac{\beta}{4}\right)$, $S_{\Omega,q,\beta,\varepsilon}<\infty$. Furthermore, there exists  $u_\eps\in H_0^2(\Omega)$
such that $\int_\Omega |\Delta u_\eps|^2dx=1$, and
\[
\int_\Omega \frac{\exp\Big({(32\pi^2(1-\frac{\beta}{4})-\eps)u_\eps^2 q(\|u_\eps\|_2^2)}\Big)}{|x|^{\beta}}dx =S_{\Omega,q,\beta,\varepsilon}.
\]
\end{proposition}

\begin{proof}

For any fixed $\eps>0$,  we can choose a maximizing sequence $\{u_j\}\subset H_0^2(\Omega)$ such that $\int_\Omega |\Delta u_j|^2dx=1$, and
\[
\int_\Omega \frac{\exp\Big({(32\pi^2(1-\frac{\beta}{4})-\eps)u_j^2 q(\|u_j\|_2^2)}\Big)}{|x|^{\beta}}dx \rightarrow
S_{\Omega,q,\beta,\varepsilon}
\]
as $j\to +\infty$. Since $\{u_j\}$ is bounded in $H_0^2(\Omega)$,  after passing to a subsequence, we have
\[
u_j\rightharpoonup u_\eps \  \text{weakly in } H_0^2(\Omega),
u_j\to u_\eps \  \text{strongly in } L^2(\Omega),
u_j\to u_\eps \  \text{a.e. in } \Omega.
\]
Hence
\[\begin{split}
f_j :&= \frac{\exp\Big({(32\pi^2(1-\frac{\beta}{4})-\eps)u_j^2 q(\|u_j\|_2^2)}\Big)}{|x|^{\beta}} \\
&\to f_\eps := \frac{\exp\Big({(32\pi^2(1-\frac{\beta}{4})-\eps)u_\eps^2 q(\|u_\eps\|_2^2)}\Big)}{|x|^{\beta}} \quad \text{a.e. in } \Omega.
\end{split}\]

We fisrt claim that $u_\eps\neq0$. Indeed, if $u_\eps=0$, the strong convergence in $L^2(\Omega)$ gives
$q(\|u_j\|_2^2)\rightarrow q(0)=1$, as $j\rightarrow \infty$.
Hence, for all sufficiently large $j$,
\begin{equation}\label{eqqe1}
\left(
32\pi^2\left(1-\frac{\beta}{4}\right)-\eps
\right)
q(\|u_j\|_2^2)
\leq
32\pi^2\left(1-\frac{\beta}{4}\right)-\frac{\eps}{2}.
\end{equation}
Choose $s>1$, sufficiently close to $1$, such that $\beta s<4$ and
\begin{equation}\label{eqqe2}
s\left(
32\pi^2\left(1-\frac{\beta}{4}\right)-\frac{\eps}{2}
\right)
<
32\pi^2\left(1-\frac{\beta s}{4}\right).
\end{equation}
Combining equations \eqref{eqqe1} and \eqref{eqqe2}, from the singular Adams inequality \eqref{eqsadm}, we get
\[\begin{aligned}
\int_\Omega f_j^sdx&=\int_\Omega \frac{\exp\Big({(32\pi^2(1-\frac{\beta}{4})-\eps)u_j^2 q(\|u_j\|_2^2)s}\Big)}{|x|^{\beta s}} dx\\
&\le\int_\Omega \frac{\exp\Big({(32\pi^2(1-\frac{\beta}{4})-\frac{\eps}{2})u_j^2 s}\Big)}{|x|^{\beta s}} dx\\
&\leq\int_\Omega\frac{\exp\left(32\pi^2\left(1-\frac{\beta s}{4}\right)u_j^2\right)}
{|x|^{\beta s}}dx\leq C.
\end{aligned}
\]

Thus $f_j$ is bounded in $L^s(\Omega)$ for some $s>1$. Since
$f_j\to |x|^{-\beta}$ almost everywhere, Vitali's convergence
theorem yields $f_j\to |x|^{-\beta}$ {strongly in } $L^1(\Omega).$
One has
\[ \int_\Omega \frac{1}{|x|^{\beta}}dx= S_{\Omega,q,\beta,\varepsilon},\]
which is impossible. Therefore,  the claim is proved.

From Proposition \ref{prop3.1}, we see that for any $p<1/(1-\|\Delta u_\eps\|_2^2)$, it holds
\begin{equation}\label{eqsampj}
\limsup_{j\to\infty}\int_\Omega \frac{\exp\Big({32\pi^2(1-\frac{\beta}{4}) p u_j^2}\Big)}{|x|^{\beta}}dx < +\infty.
\end{equation}
By assumption on $q(t)$, $0\le a_1<\lambda(\Omega)$, $0\le a_2\le \lambda(\Omega)a_1$, $\dots$, $0\le a_k\le \lambda(\Omega)a_{k-1}$, for any $u\in H_0^2(\Omega)$ with $\|u\|_2\ne 0$, there holds 
\begin{equation}\label{eqsampj1}
\begin{aligned}
q(\|u\|_2^2)\bigl(1-\|\Delta u\|_2^2\bigr)
&= 1+\sum_{j=1}^{k}a_j\|u\|_2^{2j}
   -\|\Delta u\|_2^2
   -\sum_{j=1}^{k}a_j\|u\|_2^{2j}\|\Delta u\|_2^2 \\
&< 1-a_k\|u\|_2^{2k}\|\Delta u\|_2^2<1,
\end{aligned}
\end{equation}
which implies that as $j\rightarrow +\infty$,
\begin{equation}\label{eqsampj2}
q(\|u_j\|_2^2)\to q(\|u_\eps\|_2^2) < \frac{1}{1-\|\Delta u_\eps\|_2^2}.
\end{equation}
Combining \eqref{eqsampj} and \eqref{eqsampj2}, we obtain $
S_{\Omega,q,\beta,\varepsilon}<\infty$ and ${f_j}$ is bounded in $L^r(\Omega)$ for some $r>1$. Therefore, by Vitali's convergence theorem, we get
$$
f_j\to f_\varepsilon   \text{ strongly in } L^1(\Omega).
$$
Consequently, $u_\varepsilon$ is an extremal function for $S_{\Omega,q,\beta,\varepsilon}$. This completes the proof.
\end{proof}

Let $\{u_\eps\}$ be the maximizers for $S_{\Omega,q,\beta,\varepsilon}$ obtained in Proposition~\ref{prop3.2}. A direct computation shows that $u_\eps$ satisfies the following Euler--Lagrange equation in the weak sense: 
\begin{equation}\label{eqjzfc}
\begin{cases}
\displaystyle
\Delta^2 u_\eps=\frac{d_\eps}{\lambda_\eps}\frac{u_\eps\exp(\alpha_\eps u_\eps^2)}{|x|^\beta}+\gamma_\eps u_\eps&\text{in }\Omega,\\
\displaystyle
u_\eps=\frac{\partial u_\eps}{\partial\nu}=0&\text{on }\partial\Omega,\\\|\Delta u_\eps\|_2=1,\\
\displaystyle
\alpha_\eps=\left(32\pi^2\left(1-\frac{\beta}{4}\right)-\eps\right)q(\|u_\eps\|_2^2),\\
\displaystyle
d_\eps=\frac{q(\|u_\eps\|_2^2)}{q(\|u_\eps\|_2^2)+q'(\|u_\eps\|_2^2)\int_\Omega u_\eps^2dx},\\
\displaystyle
\gamma_\eps=\frac{q'(\|u_\eps\|_2^2)}{q(\|u_\eps\|_2^2)+q'(\|u_\eps\|_2^2)\int_\Omega u_\eps^2dx},\\
\displaystyle
\lambda_\eps=\int_\Omega\frac{u_\eps^2\exp(\alpha_\eps u_\eps^2)}{|x|^\beta}dx,
\end{cases}
\end{equation}
where $q'$ denotes the derivative of $q$.

Since the constants appearing in \eqref{eqjzfc} will play an important role in the subsequent analysis, we first establish the following estimates.
\begin{lemma}\label{lem3.3}
The sequences $\{\alpha_\eps\}$, $\{d_\eps\}$, and
$\{\gamma_\eps\}$ are bounded. Moreover, if all the coefficients
$a_1,\dots,a_k$ of $q$ are nonnegative, then
$\liminf\limits_{\eps\to0}\alpha_\eps>0.$
\end{lemma}

\begin{proof}
Since $\|\Delta u_\eps\|_2=1$, the definition of
$\lambda(\Omega)$ gives $\|u_\eps\|_2^2\leq  1/ \lambda(\Omega).$
Thus both $q(\|u_\eps\|_2^2)$ and
$q'(\|u_\eps\|_2^2)$ are uniformly bounded. Moreover, since
$a_1,\dots,a_k$ are nonnegative, we get
\[q(\|u_\eps\|_2^2)+q'(\|u_\eps\|_2^2)\|u_\eps\|_2^2
=1+\sum_{j=1}^k(j+1)a_j\|u_\eps\|_2^{2j}\geq1.
\]
It follows from their definitions that $d_\eps$ and $\gamma_\eps$
are uniformly bounded. The boundedness of $\alpha_\eps$ follows from
the boundedness of $q(\|u_\eps\|_2^2)$.

Finally, since $q(\|u_\eps\|_2^2)\geq1,$
we have $\alpha_\eps\geq32\pi^2\left(1-\frac{\beta}{4}\right)-\eps.$
Consequently,
\[
\liminf_{\eps\to0}\alpha_\eps
\geq
32\pi^2\left(1-\frac{\beta}{4}\right)>0.
\]
\end{proof}

Next, we establish a uniform positive lower bound for $\lambda_\eps$.
\begin{lemma}\label{lem3.4}
There holds $\liminf\limits_{\eps\to 0}\lambda_\eps > 0$.
\end{lemma}

\begin{proof}
We first claim that 
\begin{equation}\begin{split}\label{eqlims}
\lim_{\eps\to 0}\int_\Omega \frac{\exp({\alpha_\eps u_\eps^2})}{|x|^{\beta}}dx =
\sup\limits_{{u\in H_0^2(\Omega),\|\Delta u\|_2^2=1}}
\int_\Omega \frac{\exp\Big({32\pi^2(1-\frac{\beta}{4}) u^2 q(\|u\|_2^2)}\Big)}{|x|^{\beta}}dx.
\end{split}\end{equation}
Indeed, by Proposition \ref{prop3.2}, for any $\eps>0$ it holds
\begin{equation}\label{eq082403}
\int_\Omega \frac{\exp({\alpha_\eps u_\eps^2})}{|x|^{\beta}}dx \le
\sup\limits_{{u\in H_0^2(\Omega),\|\Delta u\|_2=1}}
\int_\Omega \frac{\exp\Big({32\pi^2(1-\frac{\beta}{4}) u^2 q(\|u\|_2^2)}\Big)}{|x|^{\beta}}dx.
\end{equation}
On the other hand, for any fixed $u\in H_0^2(\Omega)$ satisfying
$\|\Delta u\|_2=1$, Proposition~\ref{prop3.2} and Fatou's lemma yield
\begin{equation}\begin{split}\label{eq082404}
&\int_\Omega \frac{\exp\Big({32\pi^2(1-\frac{\beta}{4}) u^2 q(\|u\|_2^2)}\Big)}{|x|^{\beta}}dx\\
&\le \lim_{\eps\to 0}\int_\Omega \frac{\exp\Big((32\pi^2(1-\frac{\beta}{4})-\eps)u^2 q(\|u\|_2^2)\Big)}{|x|^{\beta}}dx\\
&\le \lim_{\eps\to 0}\int_\Omega \frac{\exp({\alpha_\eps u_\eps^2})}{|x|^{\beta}}dx.
\end{split}\end{equation}
Combining inequalities \eqref{eq082403} and \eqref{eq082404}, the claim \eqref{eqlims} is proved. 

Now, applying  the inequality $e^t\le 1+te^t$ for $t\ge 0$, we have
\begin{equation*} 
\int_\Omega \frac{\exp({\alpha_\eps u_\eps^2})}{|x|^{\beta}}dx \le  \int_\Omega \frac{1}{|x|^{\beta}}dx + \alpha_\eps \lambda_\eps.
\end{equation*}
Since the right-hand side is strictly greater than
$\int_\Omega |x|^{-\beta}dx$, the boundedness of $\alpha_\eps$ established in Lemma~\ref{lem3.3} yields
\[\liminf_{\eps\to0}\lambda_\eps>0.\]
\end{proof}
 
Denote
$|c_\eps|=|u_\eps(x_\eps)|=\max_{\overline{\Omega}}|u_\eps|.$
If $\{c_\eps\}$ is bounded, then the right-hand side of
\eqref{eqjzfc} is bounded in $L^s(\Omega)$ for some $s>1$. Standard
elliptic estimates imply, after passing to a subsequence, that
$u_\eps\to u^*\ \text{strongly in }H_0^2(\Omega)$
and
$u_\eps\to u^*\ \text{in }C_{\mathrm{loc}}^4(\overline{\Omega}\setminus\{0\}),$
where $u^*\in H_0^2(\Omega)$ and $\|\Delta u^*\|_2=1$. This together with \eqref{eqlims} leads to the conclusions of both Theorem \ref{thm1} and Theorem \ref{thm2}.

In the following, we consider the blow-up case, that is, $|c_\eps|\rightarrow +\infty$ as $\eps\to0.$  Without loss of generality, we may assume that
\begin{equation}\label{eq4.3}
c_\eps=u_\eps(x_\eps)=\max_{\overline{\Omega}}u_\eps\to+\infty \text{ as }\eps\to0.
\end{equation}
Otherwise, we replace $u_\eps$ by $-u_\eps$. In the following, we  will analyze the asymptotic behavior of $u_\eps$ as $\eps\rightarrow 0$. 

\section{Asymptotic behavior of the Subcritical Maximizers}

In this section,  we will show the uniqueness of the blow-up point and study the asymptotic behavior
of $u_\eps$ both near and away from the blow-up point.
 
 We first show that the blow-up point coincides with the singular point, namely, $x_\varepsilon\rightarrow0$ (the origin).  Indeed, if $x_\varepsilon \rightarrow p\neq0$, then the weight $|x|^{-\beta}$ is uniformly bounded in a neighborhood of $p$. Since
$\alpha_\eps\to32\pi^2\left(1-\beta/4\right)<32\pi^2,$ the classical Adams inequality, together with local interior or boundary
elliptic estimates, implies that $u_\eps$ is uniformly bounded in a
neighborhood of $p$. This contradicts the assumption
$c_\eps\to+\infty$, hence,  $x_\varepsilon\rightarrow0$. On the other hand,  since $0\in\Omega$,
boundary blow-up is automatically excluded.

Since $u_\eps$ is bounded in $H_0^2(\Omega)$, there exists some $u_0\in H_0^2(\Omega)$ such that $u_\eps\rightharpoonup u_0$ weakly in $H_0^2(\Omega)$, and $u_\eps\to u_0$ strongly in $L^s(\Omega)$ for any $s>1$. In the following result, we show $u_0=0$.

\begin{lemma} 
Under the assumption that $c_\eps\to+\infty$, we have $u_\eps\rightharpoonup 0$ weakly in $H_0^2(\Omega)$, $u_\eps\rightarrow 0$ strongly in $L^s(\Omega)$ for any $s>1$, and
\begin{equation*} 
|\Delta u_\eps|^2dx\rightharpoonup\delta_0 \text{ in the sense of measures},
\end{equation*}
 where $\delta_0$ is the usual Dirac measure supported at $0$.
\end{lemma}
\begin{proof}
 Suppose $u_0\ne 0$, then $\|\Delta u_0\|_2\ne 0$. We have by \eqref{eqsampj1} that 
\[
q(\|u_\eps\|^2_2)\to q(\|u_0\|^2_2) < \frac{1}{1-\|\Delta u_0\|_2^2},
\]
as $\eps\rightarrow 0$. Choosing $r>1$ sufficiently close to $1$ such that $\beta r<4$ and applying Proposition~\ref{prop3.1}, we obtain
\[\sup_{\eps>0}\int_\Omega\left(\frac{e^{\alpha_\eps u_\eps^2}}{|x|^\beta}\right)^r dx<+\infty.\]
By the Sobolev embedding, H\"older's inequality, and
Lemmas~\ref{lem3.3}--\ref{lem3.4}, the right-hand side of the bi-harmonic equation in \eqref{eqjzfc} is bounded in $L^s(\Omega)$ for some $s>1$. Standard
elliptic estimates then imply that $\{u_\eps\}$ is uniformly bounded
in $L^\infty(\Omega)$, contradicting \eqref{eq4.3}. Hence, we must have $u_0=0$. Therefore
\begin{equation}\label{eq4.4}
\begin{cases}
u_\eps\rightharpoonup 0 & \text{weakly in } H_0^2(\Omega),\\
\nabla u_\eps\to0 &\text{strongly in } L^2(\Omega;\mathbb{R}^4),\\
u_\eps\to 0 & \text{strongly in } L^s(\Omega), \forall s>1,\\
\alpha_\eps\to 32\pi^2(1-\frac{\beta}{4}), & d_\eps\to 1,\quad \gamma_\eps\to a_1,
\end{cases}
\end{equation}
as $\eps\to0$. 

Next, we show $|\Delta u_\eps|^2dx\rightharpoonup\delta_0$. Suppose it is not true. Since
$\|\Delta u_\eps\|_2=1$, there exist $r>0$ and $\eta>0$ such that
$B_r(0)\subset\Omega$ and
\[
\limsup_{\eps\to0}
\int_{B_r(0)}|\Delta u_\eps|^2dx
\leq1-\eta.
\]
Let $\phi\in C_0^2(B_r(0))$ be a cut-off function satisfying $0\le \phi\le 1$ on $B_r(0)$ and $\phi\equiv 1$ on $B_{r/2}(0)$.
Since
\[
\Delta(\phi u_\eps)=\phi\Delta u_\eps+2\nabla\phi\cdot\nabla u_\eps+u_\eps\Delta\phi,
\]
it follows from \eqref{eq4.4} that
\[
\limsup_{\eps\to0}\int_{B_r(0)}|\Delta(\phi u_\eps)|^2dx\leq1-\eta<1.
\]

Choosing $s>1$ sufficiently close to $1$ such that  $\beta s<4$ and
\[
s\Big(1-\frac{\eta}{2}\Big)32\pi^2\Big(1-\frac{\beta}{4}\Big)<32\pi^2\Big(1-\frac{\beta s}{4}\Big).
\]
Since $\alpha_\eps\to32\pi^2\left(1-\frac{\beta}{4}\right),$ as $\eps\to0$, hence for all sufficiently small $\eps$,
\[
s\alpha_\eps\|\Delta(\phi u_\eps)\|_2^2<s\alpha_\eps\Big(1-\frac{\eta}{2}\Big)\leq32\pi^2\Big(1-\frac{\beta s}{4}\Big).
\]
Applying the singular Adams inequality \eqref{eqsadm} with $\beta$ replaced by $\beta s$ to the normalized function $\phi u_\eps/\|\Delta(\phi u_\eps)\|_2$, we obtain
\[\int_\Omega\frac{e^{s\alpha_\eps\phi^2u_\eps^2}}{|x|^{\beta s}}dx\leq C.\]
Since $\phi\equiv1$ on $B_{r/2}(0)$, it follows that
$\frac{e^{\alpha_\eps u_\eps^2}}{|x|^\beta}$ is uniformly bounded in $L^s(B_{r/2}(0))$ for some $s>1$. Using the Sobolev embedding and H\"older's inequality, we obtain that
$\frac{d_\eps}{\lambda_\eps}\frac{u_\eps e^{\alpha_\eps u_\eps^2}}{|x|^\beta}+\gamma_\eps u_\eps$
is uniformly bounded in $L^q(B_{r/2}(0))$ for some $q>1$. The interior $L^q$-estimate for the biharmonic operator yields
$\|u_\varepsilon\|_{W^{4,q}(B_{r/4}(0))}\leq C.$
Since $q>1$, the Sobolev embedding theorem implies that
$W^{4,q}(B_{r/4}(0))\hookrightarrow C^\theta(B_{r/4}(0))$
for some $\theta>0$. It follows that
$\sup_{B_{r/4}(0)}|u_\varepsilon|\leq C,$
which contradicts the assumption that $c_\varepsilon\to+\infty$. Therefore,
$|\Delta u_\varepsilon|^2\,dx\rightharpoonup\delta_0$
in the sense of measures. This completes the proof.

\end{proof}

To study the asymptotic behavior of $u_\eps$ near the blow up point $0$, we introduce the quantity $r_\eps$ as  
\begin{equation*} 
r_\eps^4 = \frac{\lambda_\eps}{d_\eps c_\eps^2} e^{-\alpha_\eps c_\eps^2}.
\end{equation*}
We claim that $r_\eps$ converges to zero rapidly. Indeed, for any $\gamma$ satisfying $0<\gamma<32\pi^2(1-\frac{\beta}{4})$,  we have 
\begin{equation}\label{eqr0speed}
r_\eps^4 c_\eps^2 e^{\gamma c_\eps^2}
= \frac{1}{d_\eps} e^{(\gamma-\alpha_\eps)c_\eps^2} \int_\Omega \frac{u_\eps^2 e^{\alpha_\eps u_\eps^2}}{|x|^{\beta}}dx
\le \frac{1}{d_\eps} \int_\Omega \frac{u_\eps^2 e^{\gamma u_\eps^2}}{|x|^{\beta}}dx \to 0,
\end{equation}
here we have used the singular Adams inequality, H\"older's inequality and \eqref{eq4.4}. Consequently, $r_\eps\to0$  as $\eps\to0$.

We now define the following rescaled functions on the corresponding rescaled domains: $\Omega_\eps = \{x\in \RR^4: x_\eps+r_\eps^{\frac{4}{4-\beta}} x\in \Omega\}$:
\begin{equation*}\begin{cases} 
\psi_\eps(x) = \frac{u_\eps\Big(x_\eps+r_\eps^{\frac{4}{4-\beta}} x\Big)}{c_\eps},\\
\phi_\eps(x) = c_\eps\bigl(u_\eps(x_\eps+r_\eps^{\frac{4}{4-\beta}} x)-c_\eps\bigr).
\end{cases}\end{equation*}

We first establish the following convergence result for $\psi_\varepsilon$.

\begin{lemma}\label{lem4.1}
For every $1<q<\frac{4}{\beta}$, as $\eps\to 0$, there holds
\[
\psi_\varepsilon\rightarrow 1
\quad\text{in }W_{\mathrm{loc}}^{4,q}(\mathbb R^4).
\]
Consequently,
\[
\psi_\varepsilon\rightarrow 1
\quad\text{in }C_{\mathrm{loc}}^{m,\alpha}(\mathbb R^4)
\]
for every integer $m\geq0$ and every $\alpha\in(0,1)$ satisfying
$m+\alpha<4-\beta$.
\end{lemma}

\begin{proof}
Using the Euler--Lagrange equation \eqref{eqjzfc}, we have

\begin{equation*}
\begin{aligned} 
\Delta^2\psi_\varepsilon(x)
&=\frac{1}{c_\varepsilon}\frac{d_\varepsilon}{\lambda_\varepsilon}r_\varepsilon^{\frac{16}{4-\beta}}
\frac{u_\varepsilon \exp({\alpha_\varepsilon u_\varepsilon^2})}{|x_\varepsilon+r_\varepsilon^{\frac{4}{4-\beta}}x|^\beta}
+\frac{r_\varepsilon^{\frac{16}{4-\beta}}}{c_\varepsilon}\gamma_\varepsilon u_\varepsilon\\
&=\frac{1}{c_\varepsilon^3}r_\varepsilon^{\frac{4\beta}{4-\beta}}
\frac{u_\varepsilon \exp({\alpha_\varepsilon u_\varepsilon^2-\alpha_\varepsilon c_\varepsilon^2})}{|x_\varepsilon+r_\varepsilon^{\frac{4}{4-\beta}}x|^\beta}
+\frac{r_\varepsilon^{\frac{16}{4-\beta}}}{c_\varepsilon}\gamma_\varepsilon u_\varepsilon\\
&=\frac{1}{c_\varepsilon^2} \psi_\varepsilon(x)
e^{\alpha_\varepsilon c_\varepsilon^2
(\psi_\varepsilon^2-1)}
 \frac{1 }{\left |x_\varepsilon r_\varepsilon^{\frac{4}{\beta-4}} +  x\right |^{\beta}} + r_\varepsilon^{\frac{16}{4-\beta}} \gamma_\varepsilon \psi_\varepsilon(x).
\end{aligned}
\end{equation*}
Since $|\psi_\varepsilon|\leq1$, it follows that
$\exp({\alpha_\varepsilon c_\varepsilon^2(\psi_\varepsilon^2-1)})\leq1.$
Moreover, for every $R>0$ and every $1<q<4/\beta$, we have
\begin{equation}\label{add912}
\sup_{a\in\mathbb{R}^4}\int_{B_R(a)}|x|^{-\beta q}dx\leq C(R,q,\beta).
\end{equation}
Hence, using the boundedness of
$\gamma_\varepsilon$, we obtain
\[
\|\Delta^2\psi_\varepsilon\|_{L^q(B_R)}\leq
C \left(\frac{1}{c_\varepsilon^2}+ r_\varepsilon^{\frac{16}{4-\beta}}\right)\rightarrow0
\]
as $\varepsilon\to 0$. The interior $L^q$-estimate for the bi-harmonic operator and
$|\psi_\varepsilon|\leq1$ imply that
$\{\psi_\varepsilon\}$ is bounded in
$W^{4,q}(B_R)$ for every $R>0$ and every
$1<q<4/\beta$. Since $q>1$, the embedding
$W^{4,q}(B_R)\hookrightarrow C^{\theta}(B_R)$
is compact for every
$0<\theta<\min\{1,4-\frac4q\}$.
Now, for every interger $m\geq0$  and $\alpha\in(0,1)$ such that
$m+\alpha<4-\beta,$ we may choose $q<4/\beta$ sufficiently close to $4/\beta$ so that
$m+\alpha<4-\frac4q.$
The Sobolev embedding
$W^{4,q}_{\mathrm{loc}}(\mathbb R^4)
\hookrightarrow
C_{\mathrm{loc}}^{m,\alpha}(\mathbb R^4)$
then gives
$\psi_\varepsilon\to \psi
\text{ in }
C_{\mathrm{loc}}^{m,\alpha}(\mathbb R^4).$

Moreover, for any $R>0$, using a change of variables, we obtain
\begin{equation*} 
\begin{aligned}
\int_{B_R}
|\Delta\psi_\varepsilon|^2dx
&=\frac1{c_\varepsilon^2}
\int_{B_{Rr_\varepsilon^{\frac{4}{4-\beta}}}(x_\varepsilon)}
|\Delta u_\varepsilon(y)|^2dy\\
&\leq\frac{\|\Delta u_\varepsilon\|_2^2}{c_\varepsilon^2}=\frac1{c_\varepsilon^2}\rightarrow0
\end{aligned}
\end{equation*}
as $\varepsilon\to 0$. Hence, we have
$\Delta\psi=0
 \text{ in }\mathbb R^4.$ Furthermore, since $|\psi|\leq1$ and $\psi(0)=1,$ the Liouville theorem for bounded harmonic functions yields
$\psi\equiv1 \text{ in }\mathbb R^4.$ This completes the proof.

\end{proof}


\begin{lemma}\label{lemvk}
For every integer $m\geq0$ and every $\alpha\in(0,1)$ satisfying $m+\alpha<4-\beta$, we have
\[
v_\eps (x) = u_\eps(x_\eps + r_\eps^{\frac{4}{4-\beta}} x) - u_\eps(x_\eps) \to 0
\quad \text{in } C_{\rm loc}^{m,\alpha}(\mathbb R^4)
\]
as $\eps \to 0$. Consequently, for every positive integer \(i<4-\beta\) and every \(R>0\),
\[
\sup_{x\in B_{Rr_\eps^{4/(4-\beta)}}(x_\eps)}|\nabla^iu_\eps(x)|
=o\left(r_\eps^{-\frac{4i}{4-\beta}}\right).
\]
\end{lemma}

\begin{proof}
By \eqref{eqjzfc} and the definition of \(r_\eps\) and  \(v_\eps\),   we have
\[
\begin{aligned}
\Delta^2v_\eps(x)
&=r_\eps^{\frac{16}{4-\beta}}\Delta^2u_\eps\left(x_\eps+r_\eps^{\frac4{4-\beta}}x\right)\\
&=\frac{\psi_\varepsilon}{c_\eps}\frac{\exp\left(\alpha_\eps c^2_\eps \left(\psi^2_\varepsilon
-1\right)\right)}{
\left|x+\dfrac{x_\eps}{r_\eps^{4/(4-\beta)}}
\right|^\beta}+\gamma_\eps r_\eps^{\frac{16}{4-\beta}}c_\eps\psi_\varepsilon.
\end{aligned}
\]
Since
$\psi_\varepsilon\leq 1,$ we have $\exp\left(\alpha_\eps c^2_\eps \left(\psi^2_\varepsilon
-1\right)\right)\leq1.$
Therefore, for every fixed \(R>0\) and every $1<p<\frac4\beta,$
by using \eqref{add912} and \eqref{eqr0speed}, we obtain
\begin{equation}\label{adddelta10-3}
    \begin{aligned}
\|\Delta^2v_\eps\|_{L^p(B_R)}
&\leq\frac1{c_\eps}\left(
\int_{B_R}\frac{dx}{
\left|x+\dfrac{x_\eps}{r_\eps^{4/(4-\beta)}}\right|^{\beta p}}\right)^{1/p}+C r_\eps^{\frac{16}{4-\beta}}c_\eps  \leq\frac{C}{c_\eps}+C r_\eps^{\frac{16}{4-\beta}}c_\eps=o(1),
\end{aligned}
\end{equation}
 as $\eps\to 0$. Moreover, since
\[
\Delta v_\eps(x)
=
r_\eps^{\frac8{4-\beta}}
\Delta u_\eps\left(
x_\eps+r_\eps^{\frac4{4-\beta}}x
\right),
\]
the change of variables $y=x_\eps+r_\eps^{\frac4{4-\beta}}x$ gives
\begin{equation}\label{add9122}
    \begin{aligned}
\int_{\Omega_\eps}|\Delta v_\eps|^2dx
&=\int_{\Omega_\eps}r_\eps^{\frac{16}{4-\beta}}
\left|\Delta u_\eps\left(x_\eps+r_\eps^{\frac4{4-\beta}}x\right)\right|^2dx
=\int_\Omega|\Delta u_\eps(y)|^2dy=1.
\end{aligned}
\end{equation}

We now claim that $v_\eps$ is bounded in $L_{\rm loc}^1(\mathbb R^4)$. Indeed, fix $R>0$,  since $x_\eps\to0$ and $r_\eps\to0$, we have $B_{4R}\subset\Omega_\eps$ for all sufficiently small $\eps$. 
 Let $w_\eps$ be a weak solution of
\[
\begin{cases}
\Delta^2w_\eps=\Delta^2v_\eps&\text{in }B_{4R},\\
\displaystyle
w_\eps=\frac{\partial w_\eps}{\partial\nu}=0
&\text{on }\partial B_{4R}.
\end{cases}
\]
Using the local $L^p$ estimate for the biharmonic operator with homogeneous Dirichlet boundary conditions and \eqref{adddelta10-3},  for any $1<p<4/\beta$, we have
\begin{equation*} 
    \|w_\eps\|_{W^{4,p}(B_{4R})}
\leq C(R)\|\Delta^2v_\eps\|_{L^p(B_{4R})}
=o(1),
\end{equation*}
as $\eps\to 0$. Since $p>1$, by the Sobolev compact embedding theorem, we have
\begin{equation}\label{addcovergence}
w_\eps\rightarrow0
 \text{ in }C^0(\overline{B_{4R}})\ \text{ and }\  \|\Delta w_\eps\|_{L^2(B_{4R})}\to 0,
\end{equation}
as $\eps\to 0$. Set
\[
h_\eps=v_\eps-w_\eps.
\]
Then
\[
\Delta^2h_\eps=0\text{ in }B_{4R}. 
\]
Hence, $\Delta h_\eps$ is harmonic in $B_{4R}$.
By the mean value property of harmonic functions and \eqref{add9122}, we have
\begin{equation}\label{addloc1}
\begin{aligned}
|\Delta h_\eps(0)|
&\leq C(R)\|\Delta h_\eps\|_{L^2(B_{2R})}\\
&\leq C(R)\left(
\|\Delta v_\eps\|_{L^2(B_{2R})}
+\|\Delta w_\eps\|_{L^2(B_{2R})}
\right)\leq C(R).
\end{aligned}
\end{equation}
Moreover, from \eqref{addcovergence}, we have
\begin{equation}\label{addloc2}
h_\eps(0)
=v_\eps(0)-w_\eps(0)
=-w_\eps(0)
=o(1).
\end{equation}
Since $h_\eps$ is biharmonic,  using Pizzetti's formula (see \cite{Pizzetti}), \eqref{addloc1} and \eqref{addloc2},  we get
\[
\int_{B_{2R}(0)}h_\eps dx=C_1 R^4h_\eps(0)+C_2R^6\Delta h_\eps(0)\leq C(R).
\]
Consequently,
\[
\left|\int_{B_{2R}}h_\eps dx\right|
\leq C(R).
\]
This together with \eqref{addcovergence} implies
\[
\begin{aligned}
\|v_\eps\|_{L^1(B_{2R})}
&=-\int_{B_{2R}}v_\eps dx\leq\left|\int_{B_{2R}}h_\eps dx\right|+\|w_\eps\|_{L^1(B_{2R})}\leq C(R).
\end{aligned}
\]
Thus the claim is proved. 

Now, by using the interior $L^p$ regularity for biharmonic operator (\cite{Gazzola}, see also Lemma 7.6 in \cite{ChLuZhu2020}), for any $1<p<\frac4\beta$,  we get 
\[
\begin{aligned}
\|v_\eps\|_{W^{4,p}(B_R)}&\leq C(R,p)\left(
\|\Delta^2v_\eps(x)\|_{L^p(B_{2R})}+\|v_\eps\|_{L^1(B_{2R})}
\right)\leq C(R,p).
\end{aligned}
\]

Let $m$ be a nonnegative integer and let $0\leq\theta<1$ satisfy $m+\theta<4-\beta.$ Choose $p$ such that
\[
1<p<\frac4\beta, m+\theta<4-\frac4p.
\]
Then by the compact Sobolev embedding and \eqref{adddelta10-3}, we have
$v_\eps\rightarrow v$ {in }$C_{\rm loc}^{m,\theta}(\mathbb R^4),$ with  $v$ satisfying
\[
\Delta^2v=0
\text{ in }\mathbb R^4,  v \leq0,
v(0)=0.
\]

On the other hand, for every $R>0$, through the weak lower semicontinuity of the
$L^2$-norm and \eqref{add9122}, we obtain
\[
\int_{B_R}|\Delta v|^2dx
\leq\liminf_{\eps\to0}\int_{B_R}|\Delta v_\eps|^2dx
\leq1.
\]
By the mean value property of the harmonic function $\Delta v$, for any fixed $x\in \mathbb{R}^4$ and $R>0$, we have 
\[|\Delta v(x)| \leq \frac{1}{|B_R|} \int_{B_R(x)} |\Delta v |dx \leq |B_R|^{-1/2} \|\Delta v\|_{L^2(\mathbb{R}^4)} \leq C R^{-2}.\]
Letting $R\rightarrow \infty$, we get  $\Delta v(x)=0$, that means that $v$ is also a harmonic function in $\mathbb{R}^4$.  Since $v\leq0$ and $v(0)=0$, the Liouville theorem for harmonic functions implies $v\equiv0.$

Finally, let $i$ be a positive integer satisfying $i<4-\beta$. Choose $\theta>0$ such that $i+\theta<4-\beta.$
Then
\[
\nabla^iv_\eps\rightarrow0
 \text{ locally uniformly in }\mathbb R^4.
\]
Since
$\nabla^iv_\eps(x)
=r_\eps^{\frac{4i}{4-\beta}}\nabla^iu_\eps\left(
x_\eps+r_\eps^{\frac4{4-\beta}}x \right),$ we conclude that  for every fixed $R>0$,
\[
\begin{aligned}
\sup_{y\in B_{Rr_\eps^{4/(4-\beta)}}(x_\eps)}
|\nabla^iu_\eps(y)|&=r_\eps^{-\frac{4i}{4-\beta}}
\sup_{x\in B_R}|\nabla^iv_\eps(x)|=
o\left(r_\eps^{-\frac{4i}{4-\beta}}\right).
\end{aligned}
\]
The proof is finished.
\end{proof}

\begin{lemma}[Gradient estimate on $B_{Rr_\eps^{\frac{4}{4-\beta}}}(x_\eps)$]\label{lemge}
For any $R > 0$, there holds
\[
c_\eps \int_{B_{Rr_\eps^{\frac{4}{4-\beta}}}} |\Delta u_\eps| dx \leq C (Rr_\eps^{\frac{4}{4-\beta}})^2.
\]
Furthermore, we have
\[
\int_{B_R} |\Delta \phi_\eps| dx
= c_\eps r_\eps^{-\frac8{4-\beta}} \int_{B_{Rr_\eps^{\frac{4}{4-\beta}}}} |\Delta u_\eps| dx
\leq C R^2.
\]
\end{lemma}

\begin{proof}
Choose \(R_0>0\) sufficiently small such that $B_{8R_0}(0)\Subset\Omega.$
Since \(x_\eps\to0\), it follows that $B_{4R_0}(x_\eps)\Subset\Omega$ for all sufficiently small \(\eps\).

Let \(u_\eps^{R_0}\) be the biharmonic replacement of
\(u_\eps\) in \(B_{4R_0}(x_\eps)\), namely,
\[
\begin{cases}
\Delta^2u_\eps^{R_0}=0
&\text{in }B_{4R_0}(x_\eps),\\
u_\eps^{R_0}=u_\eps, 
\displaystyle
\frac{\partial u_\eps^{R_0}}{\partial\nu}
=
\frac{\partial u_\eps}{\partial\nu}
&\text{on }\partial B_{4R_0}(x_\eps).
\end{cases}
\]
Using the weak biharmonicity of $u_\eps^{R_0}$ and the fact that
\[
u_\eps-u_\eps^{R_0}
\in H_0^2(B_{4R_0}(x_\eps)),
\]
we obtain
\[
\int_{B_{4R_0}(x_\eps)}
\Delta u_\eps^{R_0}
\Delta(u_\eps-u_\eps^{R_0})dx
=0.
\]
Consequently,
\[
\begin{aligned}
\int_{B_{4R_0}(x_\eps)}|\Delta u_\eps|^2dx
&=\int_{B_{4R_0}(x_\eps)}|\Delta u_\eps^{R_0}|^2dx+
\int_{B_{4R_0}(x_\eps)}|\Delta(u_\eps-u_\eps^{R_0})|^2dx.
\end{aligned}
\]
Since $\|\Delta u_\eps\|_2=1$, it follows that
\begin{equation}\label{add9123bound}
  \|\Delta u_\eps^{R_0}\|_{L^2(B_{4R_0}(x_\eps))}
+
\|\Delta(u_\eps-u_\eps^{R_0})\|_{L^2(B_{4R_0}(x_\eps))}
\leq C.  
\end{equation}

Moreover, since $u_\eps-u_\eps^{R_0}\in H_0^2(B_{4R_0}(x_\eps)),$ the Poincar\'{e} inequality gives
\begin{equation}\label{addl2}
    \|u_\eps-u_\eps^{R_0}\|_{L^2(B_{4R_0}(x_\eps))}
\leq
CR_0^2
\|\Delta(u_\eps-u_\eps^{R_0})\|_{L^2(B_{4R_0}(x_\eps))}
\leq C(R_0).
\end{equation}On the other hand, by the definition of $\lambda(\Omega)$, we have   \begin{equation}
    \label{boundforl2}\|u_\eps\|_{L^2(\Omega)}\leq\frac{1}{\sqrt{\lambda(\Omega)}}. 
\end{equation}
This, together with \eqref{addl2} implies
\[
\|u_\eps^{R_0}\|_{L^2(B_{4R_0}(x_\eps))}
\leq C(R_0).
\]
Since $u_\eps^{R_0}$ is biharmonic in
$B_{4R_0}(x_\eps)$, the interior estimates for biharmonic
functions (see \cite{Gazzola}) yield
\begin{equation}\label{eqbhcr}
\|u_\eps^{R_0}\|_{C^4(B_{3R_0}(x_\eps))}\leq C(R_0)
\|u_\eps^{R_0}\|_{L^2(B_{4R_0}(x_\eps))}\leq C(R_0).
\end{equation}

Now, we show that  $\Delta^2(u_\eps-u_\eps^{R_0})$ is bounded in $L^1(B_{2R_0}(x_\eps))$.  Note that from \eqref{eqjzfc}, $u_\eps-u_\eps^{R_0}$  satisfies the following equation
\begin{equation} \label{definforf}
\Delta^2(u_\eps-u_\eps^{R_0})
=\frac{d_\eps}{\lambda_\eps}\frac{u_\eps e^{\alpha_\eps u_\eps^2}}{|x|^\beta}+\gamma_\eps u_\eps=:f_\eps
 \text{ in }B_{4R_0}(x_\eps).
\end{equation}
We first claim that  $f_\eps$ is bounded in $ L(\log L)^{\frac{1}{2}}(B_{4R_0}(x_\eps))$, that is,
\begin{equation} \label{addclaim104}
\int_{B_{4R_0}(x_\eps)}|f_\eps|\left(\log(2+|f_\eps|)\right)^{1/2}dx\leq C.
\end{equation}
Indeed, since $d_\eps$, $\gamma_\eps$ and
$\lambda_\eps^{-1}$ are uniformly bounded, we have
\begin{equation}\label{addcontf104}
|f_\eps|\leq\frac{C}{\lambda_\eps}\frac{|u_\eps|e^{\alpha_\eps u_\eps^2}}{|x|^\beta}+C|u_\eps|.
\end{equation}
Consequently,
\begin{equation*} 
\begin{aligned}
\log(2+|f_\eps|)&\leq C+\log(1+|u_\eps|)
+\alpha_\eps u_\eps^2+\beta|\log{|x|}|\leq C\left( 1+u_\eps^2+|\log{|x|}|\right).
\end{aligned}
\end{equation*}
Therefore, we have
\[
\left(\log(2+|f_\eps|)\right)^{1/2}
\leq C \left(1+|u_\eps|+\left(|\log{|x|}|\right)^{1/2}\right)
\]
and 
\begin{equation} \label{add2104}
\begin{split}
|f_\eps| \left(\log(2+|f_\eps|)\right)^{1/2}
\leq & C \left(|f_\eps| +|f_\eps u_\eps|+|f_\eps|\left(|\log{|x|}|\right)^{1/2}\right). 
\end{split}\end{equation}


Splitting \(\Omega\) into
\(\{|u_\eps|\leq1\}\) and \(\{|u_\eps|>1\}\), from the definition of $\lambda_\eps$ and  Lemma \ref{lem3.4}, we obtain
\begin{equation}\begin{split}\label{eq091101}
\frac1{\lambda_\eps}\int_\Omega\frac{|u_\eps|e^{\alpha_\eps u_\eps^2}}{|x|^\beta}dx
&\leq\frac1{\lambda_\eps}
\int_{\{|u_\eps|\leq1\}}\frac{e^{\alpha_\eps}}{|x|^\beta}dx+\frac1{\lambda_\eps}\int_{\{|u_\eps|>1\}}\frac{u_\eps^2e^{\alpha_\eps u_\eps^2}}{|x|^\beta}dx\leq C.
\end{split}\end{equation}
Moreover, using H\"older's inequality and \eqref{boundforl2},  we have
\begin{equation}\label{eq091102}
\int_{B_{4R_0}(x_\eps)}|u_\eps|dx\leq|\Omega|^{1/2}\|u_\eps\|_{L^2(\Omega)}\leq C.
\end{equation}
Combining \eqref{addcontf104}, \eqref{eq091101} and \eqref{eq091102}, we see that
\begin{equation}\label{eq091103}
\int_{B_{4R_0}(x_\eps)}|f_\eps|dx\leq C.
\end{equation}
 Moreover, by \eqref{boundforl2}, we have
 \begin{equation}\label{eq091104}
\int_{B_{4R_0}(x_\eps)}|f_\eps u_\eps|dx=\frac{d_\eps}{\lambda_\eps}
\int_\Omega\frac{u_\eps^2e^{\alpha_\eps u_\eps^2}}{|x|^\beta}dx+\int_\Omega\gamma_\eps u_\eps^2 dx\leq C.
\end{equation} 

It remains to consider the estimate for $|f_\eps|\left(|\log{|x|}|\right)^{1/2}$:
\begin{equation*}\begin{split}
\int_{B_{4R_0}(x_\eps)}|f_\eps|\left(|\log{|x|}|\right)^{1/2}dx
&\leq \frac {C}{\lambda_\eps}\int_{B_{4R_0}(x_\eps)}\frac{|u_\eps|e^{\alpha_\eps u_\eps^2}}{|x|^\beta}\left(|\log|x||\right)^{1/2}dx\\
&\ \ +C\int_{B_{4R_0}(x_\eps)}  |u_\eps|\left(|\log|x||\right)^{1/2} dx\\
&:=I_1+I_2.
\end{split}\end{equation*}
For $I_2$,  by \eqref{boundforl2} and H\"{o}lder's inequality, we have
\begin{equation}\begin{split}\label{eq091113}
I_2\leq\int_{B_{4R_0}(x_\eps)}|u_\eps|\left(|\log|x||\right)^{1/2}dx
&\leq\|u_\eps\|_{L^2(\Omega)}\left(\int_{\Omega}|\log|x||dx\right)^{1/2}\leq C.
\end{split}\end{equation}
Next, we estimate $I_1$. We divide ${B_{4R_0}(x_\eps)} $ into  $$\Big\{B_{4R_0}(x_\eps)\cap \{|x|\geq1\}\Big\}\cup\Big\{B_{4R_0}(x_\eps)\cap \{|x|<1\}\Big\}.$$ 
On the domain ${B_{4R_0}(x_\eps)}\cap \{|x|\geq1\}$, we have $\left(|\log|x||\right)^{1/2}\leq C$, since $\Omega$ is bounded. Then, by \eqref{eq091101}, we have
\begin{equation}\begin{split}\label{eq091112}
&\frac {C}{\lambda_\eps}\int_{B_{4R_0}(x_\eps) \cap \{|x|\geq1\}}\frac{|u_\eps|e^{\alpha_\eps u_\eps^2}}{|x|^\beta}\left(|\log|x||\right)^{1/2}dx
\leq \frac {C}{\lambda_\eps}\int_{\{|x|\geq1\}}\frac{|u_\eps|e^{\alpha_\eps u_\eps^2}}{|x|^\beta}dx\leq C.\\
\end{split}\end{equation}

 We next estimate the integral over the domain ${B_{4R_0}(x_\eps)}\cap \{|x|<1\}$. Since $0<\beta<4$, we can choose $\eta>0$ sufficiently small such that $\beta+\eta\sup_\eps\alpha_\eps<4.$  We further decompose  ${B_{4R_0}(x_\eps)}\cap \{|x|<1\}$  into the two subsets $B_{4R_0}(x_\eps)\cap \{|x|<1\} \cap\left\{u_\eps^2\geq\eta|\log|x||  \right\}$ and $B_{4R_0}(x_\eps)\cap \{|x|<1\} \cap\left\{u_\eps^2<\eta|\log|x||  \right\}$. 
 
 On the set $B_{4R_0}(x_\eps)\cap \{|x|<1\} \cap\left\{u_\eps^2\geq\eta|\log|x||  \right\},$
we have $\left(|\log|x||\right)^{1/2}\leq\frac{|u_\eps|}{\sqrt\eta}.$
Therefore,
\begin{equation}\begin{split}\label{eq091105}
&\frac1{\lambda_\eps}
\int_{B_{4R_0}(x_\eps)\cap \{|x|<1\}
\cap\{u_\eps^2\geq\eta|\log|x||\}}
\frac{|u_\eps|e^{\alpha_\eps u_\eps^2}}{|x|^\beta}
\left(|\log|x||\right)^{1/2}dx\\
&\leq\frac1{\sqrt\eta\lambda_\eps}
\int_\Omega\frac{u_\eps^2e^{\alpha_\eps u_\eps^2}}{|x|^\beta}dx
=\frac1{\sqrt\eta}\leq C.
\end{split}\end{equation}
On the complementary set $B_{4R_0}(x_\eps)\cap \{|x|<1\}\cap\left\{u_\eps^2<\eta|\log|x||\right\},$
we have
$$
|u_\eps|\left(|\log|x||\right)^{1/2}<\sqrt\eta|\log|x||.
$$
Notice $|x|<1$ on this domain, it holds
\[
\begin{aligned}
e^{\alpha_\eps u_\eps^2}
&< e^{\eta\alpha_\eps|\log|x||}=|x|^{-\eta\alpha_\eps}
\leq|x|^{-\eta\sup_\eps\alpha_\eps}.
\end{aligned}
\]
Using Lemma \ref{lem3.4} and the assumption $\beta+\eta\sup_\eps\alpha_\eps<4$, we obtain
\begin{equation*}\begin{split} 
&\frac1{\lambda_\eps}\int_{B_{4R_0}(x_\eps)\cap\{|x|<1\}
\cap\{u_\eps^2<\eta|\log|x||\}}
\frac{|u_\eps|e^{\alpha_\eps u_\eps^2}}{|x|^\beta}
\left(|\log|x||\right)^{1/2}dx\\
&\leq C\int_{B_1(0)}\frac{|\log|x||}{|x|^{\beta+\eta\sup_\eps\alpha_\eps}}dx\\
&=2\pi^2C\int_0^1r^{3-\beta-\eta\sup_\eps\alpha_\eps}|\log r|dr\leq C.
\end{split}\end{equation*}
Combining this estimate with \eqref{eq091105}, we obtain
\begin{equation}\begin{split}\label{eq091114}
&\frac1{\lambda_\eps}\int_{B_{4R_0}(x_\eps)\cap\{|x|<1\}}
\frac{|u_\eps|e^{\alpha_\eps u_\eps^2}}{|x|^\beta}
\left(|\log|x||\right)^{1/2}dx\leq C.\\
\end{split}\end{equation}
It then follows from \eqref{eq091112} and \eqref{eq091114} that
$I_1\leq C$. Finally, combining this bound with \eqref{eq091113} yields
\begin{equation*}\begin{split} 
\int_{B_{4R_0}(x_\eps)}|f_\eps|\left(|\log{|x|}|\right)^{1/2}dx
&\leq C.
\end{split}\end{equation*}
This, together with  \eqref{eq091103}, \eqref{eq091104} and \eqref{add2104} implies $f_\eps$ is bounded in  $L(\log L)^{\frac{1}{2}}(B_{4R_0}(x_\eps)),$ and the claim \eqref{addclaim104} is proved. 

By \eqref{addclaim104} and the higher-order \(L(\log L)^{1/2}\) estimate (see \cite{Martinazzi09}), we can obtain
\begin{equation}\label{eq09110812}
\left\|\nabla^j(u_\eps-u_\eps^{R_0})\right\|_{L^{4/j,2}(B_{3R_0}(x_\eps))}\leq C, \ \   j=1,2,3,
\end{equation}
where \(L^{4/j,2}\) denotes the Lorentz space.  For the definition of Lorentz spaces and their basic properties, we refer the interested reader to \cite{ONeil1963}.

We now establish the $L^1(B_{2R_0}(x_\eps))$-boundedness of $\Delta^2(u_\eps-u_\eps^{R_0})$. A direct calculation gives
\begin{equation}\begin{split}\label{eq091107}
\Delta^2
\left((u_\eps-u_\eps^{R_0})^2\right)
&=2(u_\eps-u_\eps^{R_0})\Delta^2(u_\eps-u_\eps^{R_0})+C \sum_{j=1}^{3} \left| \nabla^j (u_\eps - u_\eps^{R_0}) \right| \left| \nabla^{4-j} (u_\eps - u_\eps^{R_0}) \right|.
\end{split}\end{equation}

We first estimate the contribution of the first term on the right-hand side of \eqref{eq091107}. Since $u_\eps^{R_0}$ is biharmonic on $B_{2R_0}(x_\eps)$, we have
\begin{equation*}\begin{split} 
&\int_{B_{2R_0}(x_\eps)}\left|2(u_\eps-u_\eps^{R_0})\Delta^2(u_\eps-u_\eps^{R_0})\right|dx\\
&\leq 2\int_{B_{2R_0}(x_\eps)}\left|u_\eps\Delta^2(u_\eps)\right|dx+
2\int_{B_{2R_0}(x_\eps)}\left|u_\eps^{R_0}\Delta^2(u_\eps)\right|dx
:=I_3+I_4.
\end{split}\end{equation*}
For $I_3$, by \eqref{eqjzfc} and \eqref{boundforl2}, we get
\[
\begin{aligned}
I_3=2\int_\Omega|u_\eps\Delta^2u_\eps|dx
&\leq 2\frac{d_\eps}{\lambda_\eps}\int_\Omega\frac{u_\eps^2e^{\alpha_\eps u_\eps^2}}{|x|^\beta}dx+2|\gamma_\eps|
\int_\Omega u_\eps^2dx\leq C.
\end{aligned}
\]
For $I_4$, using the uniform \(C^4\)-bound \eqref{eqbhcr} for
\(u_\eps^{R_0}\) and \eqref{eq091103}, we have $I_4\leq C$.
Therefore, the first term of \eqref{eq091107} is bounded in \(L^1(B_{2R_0})\). In virtue of the  Lorentz estimates of \eqref{eq09110812} and H\"{o}lder type inequality of O'Neil \cite{ONeil1963}, one  can easily derive  that the remaining terms of \eqref{eq091107} are also bounded in \(L^1(B_{2R_0})\) . 
Therefore, we have
\[
\int_{B_{2R_0}(x_\eps)}\left|\Delta^2\left((u_\eps-u_\eps^{R_0})^2\right)\right|dx
\leq C.
\]
Now, carrying out the same argument as in the proof of Lemma~6 in \cite{Martinazzi09} (see also the proof of Lemma~3.3 in \cite{ChLuZhu2020}), for every fixed \(R>0\), we have
\begin{equation}\label{L1bounded104}
\begin{aligned}
&\int_{B_{Rr_\eps^{\frac4{4-\beta}}}(x_\eps)}
\left|\Delta\left((u_\eps-u_\eps^{R_0})^2\right)
\right|dx\leq C\left(Rr_\eps^{\frac4{4-\beta}}\right)^2.
\end{aligned}
\end{equation}


Next, we estimate $\int_{B_{Rr_\eps^{\frac4{4-\beta}}}(x_\eps)}|\Delta(u_\eps^2)|dx$. Since $u_\eps=(u_\eps-u_\eps^{R_0})+u_\eps^{R_0},$ we have
\[
\begin{aligned}
\Delta(u_\eps^2)
&=\Delta\left(
(u_\eps-u_\eps^{R_0})^2
\right)+2u_\eps^{R_0}\Delta(u_\eps-u_\eps^{R_0})+2(u_\eps-u_\eps^{R_0})\Delta u_\eps^{R_0}\\
&\ \ +4\nabla u_\eps^{R_0}\cdot\nabla(u_\eps-u_\eps^{R_0})+\Delta\left((u_\eps^{R_0})^2\right).
\end{aligned}
\]
Recalling the \(C^4\)-bound \eqref{eqbhcr} for \(u_\eps^{R_0}\) and the uniform \(H_0^2(B_{4R_0}(x_\eps))\)-bound of $\{u_\eps-u_\eps^{R_0}\}$ (see \eqref{add9123bound} and \eqref{addl2}), by the H\"{o}lder's inequality and \eqref{L1bounded104},  we have

\begin{equation}
\begin{split}\label{eq091110}
&\int_{B_{Rr_\eps^{\frac4{4-\beta}}}(x_\eps)}|\Delta(u_\eps^2)|dx\\
&\leq\int_{B_{Rr_\eps^{\frac4{4-\beta}}}(x_\eps)}
\left|\Delta\left((u_\eps-u_\eps^{R_0})^2\right)\right|dx
+2\int_{B_{Rr_\eps^{\frac4{4-\beta}}}(x_\eps)}
|u_\eps^{R_0}\Delta(u_\eps-u_\eps^{R_0})|dx\\
&\quad
+2\int_{B_{Rr_\eps^{\frac4{4-\beta}}}(x_\eps)}
|(u_\eps-u_\eps^{R_0})\Delta u_\eps^{R_0}|dx
+4\int_{B_{Rr_\eps^{\frac4{4-\beta}}}(x_\eps)}
|\nabla u_\eps^{R_0}\cdot\nabla(u_\eps-u_\eps^{R_0})|dx\\
&\quad
+\int_{B_{Rr_\eps^{\frac4{4-\beta}}}(x_\eps)}
|\Delta((u_\eps^{R_0})^2)|dx
\leq C\left(Rr_\eps^{\frac4{4-\beta}}\right)^2.
\end{split}
\end{equation}

By Lemma~\ref{lemvk}, $u_\eps(x)\geq\frac{c_\eps}{2}$ in \(B_{Rr_\eps^{\frac4{4-\beta}}}(x_\eps)\) for  sufficiently small \(\eps\),  Together with $\Delta(u_\eps^2)=2 u_\eps\Delta u_\eps+2|\nabla u_\eps|^2 $, it follows that 

\begin{equation}\label{add91243}
    \begin{aligned}
c_\eps|\Delta u_\eps|&\leq2u_\eps|\Delta u_\eps|\leq|\Delta(u_\eps^2)|+2|\nabla u_\eps|^2, \  \text{for } x\in B_{Rr_\eps^{\frac4{4-\beta}}}(x_\eps).
\end{aligned}
\end{equation}
 
Using the Sobolev embedding  $H_0^2(\Omega)\hookrightarrow W^{1,4}(\Omega)$, we have
\begin{equation*}\begin{split} 
\int_{B_{Rr_\eps^{\frac4{4-\beta}}}(x_\eps)}
|\nabla u_\eps|^2dx
&\leq\left|B_{Rr_\eps^{\frac4{4-\beta}}}(x_\eps)\right|^{1/2}
\|\nabla u_\eps\|_{L^4(\Omega)}^2\leq C\left(Rr_\eps^{\frac4{4-\beta}}\right)^2.
\end{split}\end{equation*}
Combining this with \eqref{add91243} and \eqref{eq091110} implies
\begin{equation}\label{addcotra3104}
c_\eps
\int_{B_{Rr_\eps^{\frac4{4-\beta}}}(x_\eps)}
|\Delta u_\eps|dx
\leq
C\left(
Rr_\eps^{\frac4{4-\beta}}
\right)^2.
\end{equation}

Finally, by the definition of \(\phi_\eps\),
\[
\Delta\phi_\eps(x)
=
c_\eps r_\eps^{\frac8{4-\beta}}
\Delta u_\eps\left(
x_\eps+r_\eps^{\frac4{4-\beta}}x
\right).
\]
Making the change of variables $y=x_\eps+r_\eps^{\frac4{4-\beta}}x$ 
and using \eqref{addcotra3104}, we obtain
\[
\begin{aligned}
\int_{B_R}|\Delta\phi_\eps|dx
&=c_\eps r_\eps^{-\frac8{4-\beta}}\int_{B_{Rr_\eps^{\frac4{4-\beta}}}(x_\eps)}|\Delta u_\eps(y)|dy\leq CR^2.
\end{aligned}
\]
This completes the proof.
\end{proof}

With this gradient estimate at hand, we now investigate the limiting behavior of $\phi_\eps$ as $\eps\to 0$.
 \begin{lemma}\label{lemxe}
For every $1<q<\frac{4}{\beta}$,  there is a subsequence such that $\phi_\eps\rightharpoonup\phi$ in $W_{\mathrm{loc}}^{4,q}(\RR^4)$ and $\phi_\eps\to\phi$ in $C^4_{loc}(\mathbb{R}^4\backslash \{0\})$.  The limiting function $\phi$ is radially symmetric and radially decreasing about the origin, and satisfies
\begin{equation}\label{eq:phi}
\Delta^2\phi(x)=\frac{e^{64\pi^2(1-\frac\beta4)\phi(x)}}{|x|^\beta}  \text{in } \RR^4, \ 
\phi\leq\phi(0)=0,
\end{equation}
in the sense of distributions. Moreover,  $\phi$ satisfies the mass identity
\begin{equation}\label{eq:phiMassShift}
\int_{\RR^4}\frac{e^{64\pi^2(1-\frac\beta4)\phi(x)}}{|x|^\beta}dx=1.
\end{equation}
\end{lemma}
\begin{proof}
We first claim  that $\frac{|x_\eps|}{r_\eps^{\frac{4}{4-\beta}}}$ is bounded.
Indeed, if $\frac{|x_\varepsilon|}{r_\varepsilon^{\frac{4}{4-\beta}}}\rightarrow +\infty$  as $\varepsilon\rightarrow0$, we define the rescaled functions as following

\begin{equation*}\begin{cases} 
\displaystyle
\widetilde\psi_\varepsilon(x):=\frac{u_\varepsilon(x_\varepsilon+r_\varepsilon |x_\varepsilon|^{\frac{\beta}{4}} x)}{c_\varepsilon},\\
\displaystyle
\widetilde\phi_\varepsilon(x):=c_\varepsilon\bigl(u_\varepsilon(x_\varepsilon+r_\varepsilon |x_\varepsilon|^{\frac{\beta}{4}} x)-c_\varepsilon\bigr).
\end{cases}\end{equation*}
Observe that $|\widetilde\psi_\varepsilon|\le1$, $\widetilde\psi_\varepsilon(0)=1$, and $\widetilde\phi_\varepsilon\le0$, $\widetilde\phi_\varepsilon(0)=0$. From the Euler--Lagrange equation \eqref{eqjzfc}, we get
\begin{equation}\begin{split}\label{eqphiw1}
\Delta^2\widetilde\phi_\varepsilon(x)
&=c_\varepsilon r_\varepsilon^4|x_\varepsilon|^\beta
\Delta^2u_\varepsilon\left(x_\varepsilon+r_\varepsilon |x_\varepsilon|^{\beta/4}x\right)\\
&=\frac{|x_\varepsilon|^\beta}{\left|x_\varepsilon+r_\varepsilon |x_\varepsilon|^{\frac{\beta}{4}}x\right|^\beta}
\widetilde\psi_\varepsilon(x)
\exp\left[\alpha_\varepsilon(1+\widetilde\psi_\varepsilon(x))\widetilde\phi_\varepsilon(x)\right] +\gamma_\varepsilon c_\varepsilon^2 r_\varepsilon^4 |x_\varepsilon|^\beta \widetilde\psi_\varepsilon(x).
\end{split}\end{equation}

From \eqref{eqr0speed}, one can easily verify that the second term of  \eqref{eqphiw1} vanishes when $\varepsilon\rightarrow0$. Since  $\frac{|x_\varepsilon|}{r_\varepsilon^{\frac{4}{4-\beta}}}\rightarrow +\infty$, for every fixed $R>0$,
$\frac{|x_\varepsilon|^\beta}{\left|x_\varepsilon+r_\varepsilon |x_\varepsilon|^{\beta/4}x\right|^\beta}\rightarrow 1$  uniformly on $ B_R(0)$. Thus, the weight appearing in \eqref{eqphiw1} is asymptotically constant on every fixed ball. Repeating the arguments used in the proofs of 
Lemmas~\ref{lem4.1}, \ref{lemvk} and~\ref{lemge}, with \(r_\eps^{4/(4-\beta)}\) replaced by \(r_\eps|x_\eps|^{\beta/4}\), we can obtain 
\[
u_\eps\left(
x_\eps+r_\eps|x_\eps|^{\frac{\beta}{4}}x
\right)-c_\eps
\rightarrow0
\quad\text{in }C_{\rm loc}^3(\mathbb R^4),
\]
\begin{equation}\label{addl1boudn10-4}
\int_{B_R}|\Delta\widetilde\phi_\eps|dx\leq CR^2,
\end{equation}
and $\widetilde\psi_\varepsilon\to1$, $\widetilde\phi_\varepsilon\to\widetilde\phi$
locally in  $C^4_{\rm loc}(\mathbb R^4)$, where the limit solution $\widetilde\phi$ satisfies the Liouville equation
\begin{equation}\begin{split} \label{addlimingt10-4}
\Delta^2\widetilde\phi=e^{64\pi^2(1-\frac{\beta}{4})\widetilde\phi} \ \text{ in } \mathbb R^4\  \text{ and }
\widetilde\phi\le\widetilde\phi(0)=0.
\end{split}\end{equation}

Here and below, the integrands are extended by zero outside the corresponding rescaled domain. From the definition of $\lambda_\varepsilon$ and $r_\varepsilon$, we get
\begin{equation}\begin{split}\label{eqphiwjx13}
\lambda_\epsilon &= \int_\Omega \frac{u^2_\epsilon(y) e^{\alpha_\epsilon u^2_\epsilon(y)}}{|y|^\beta}   dy\\
&=\int_{\mathbb{R}^4}
\frac{c_\epsilon^2 \widetilde\psi^2_\epsilon(x) e^{\alpha_\epsilon c_\epsilon^2 \widetilde\psi^2_\epsilon(x)}}{|x_\epsilon+r_\varepsilon |x_\varepsilon|^{\beta/4} x|^\beta}
 (r_\varepsilon |x_\varepsilon|^{\beta/4})^4   dx\\
&= \frac{\lambda_\epsilon}{d_\epsilon}\int_{\mathbb{R}^4}
\frac{|x_\epsilon|^\beta}{|x_\epsilon+r_\varepsilon |x_\varepsilon|^{\beta/4} x|^\beta}
\widetilde\psi^2_\epsilon(x)e^{\alpha_\epsilon c_\epsilon^2(\widetilde\psi^2_\epsilon(x)-1)} dx.\\
\end{split}\end{equation}

Since $\widetilde\psi_\varepsilon^2-1 = (1+\widetilde\psi_\varepsilon)\widetilde\phi_\varepsilon/c_\varepsilon^2$, we have $\alpha_\varepsilon c_\varepsilon^2(\widetilde\psi_\varepsilon^2-1)
= \alpha_\varepsilon(1+\widetilde\psi_\varepsilon)\widetilde\phi_\varepsilon.$ From \eqref{eq4.4}  and \eqref{eqphiwjx13} , we have
\begin{equation*}\begin{split} 
\int_{\mathbb{R}^4}
\frac{|x_\varepsilon|^\beta}{|x_\varepsilon+r_\varepsilon |x_\varepsilon|^{\beta/4} x|^\beta}
\widetilde\psi_\varepsilon^2(x)
\exp\left(\alpha_\varepsilon(1+\widetilde\psi_\varepsilon(x))\widetilde\phi_\varepsilon(x)\right)dx
=d_\varepsilon\rightarrow 1,\end{split}\end{equation*}
as $\varepsilon\rightarrow0$. This, together with \eqref{eq4.4} and the fact that $\frac{|x_\varepsilon|^\beta}{\left|x_\varepsilon+r_\varepsilon |x_\varepsilon|^{\beta/4}x\right|^\beta}\rightarrow 1$ a.e. in $\mathbb{R}^4$ as $\varepsilon\rightarrow0$, implies
\begin{equation}\label{eqmass1}
\int_{\mathbb{R}^4} e^{64\pi^2(1-\frac{\beta}{4})\widetilde\phi}dx \le 1.
\end{equation}

Set
\[
\widetilde\Phi
=
\left(1-\frac{\beta}{4}\right)\widetilde\phi
+
\frac{1}{64\pi^2}
\log\left(1-\frac{\beta}{4}\right).
\]
Then from \eqref{addl1boudn10-4} and \eqref{addlimingt10-4}, we have
\[
\Delta^2\widetilde\Phi=e^{64\pi^2\widetilde\Phi}
\quad\text{in }\mathbb R^4, \text{ with } \int_{B_R}|\Delta\widetilde\Phi|dx\leq CR^2.
\]
Therefore,  arguing as in the proof of \cite[Lemma~3.4]{ChLuZhu2020} (see also \cite{Lin1998}), we obtain
\[
\int_{\mathbb R^4}e^{64\pi^2\widetilde\Phi}dx=1.
\]
Consequently,
\begin{equation*}\begin{split} 
\int_{\mathbb{R}^4} e^{64\pi^2(1-\frac{\beta}{4})\widetilde\phi}dx
=\frac{1}{1-\frac{\beta}{4}}>1,
\end{split}\end{equation*}
which contradicts \eqref{eqmass1}.
Therefore, the claim $\frac{|x_\varepsilon|}{r_\varepsilon^{\frac{4}{4-\beta}}}\le C$ is proved.

 Now, we assume that up to a subsequence, 
\begin{equation*} 
\frac{x_\eps}{r_\eps^{\frac{4}{4-\beta}}}
\rightarrow x^*
\quad\text{in }\mathbb R^4.
\end{equation*}

By the Euler--Lagrange equation \eqref{eqjzfc}, we get
\begin{align}
\Delta^2\phi_\eps(x)
=&\frac{\psi_\eps(x)\exp\left(
\alpha_\eps(1+\psi_\eps(x))\phi_\eps(x)
\right)}{\left|x+r_\eps^{-\frac{4}{4-\beta}}x_\eps
\right|^\beta}+\gamma_\eps c_\eps^2r_\eps^{\frac{16}{4-\beta}}\psi_\eps(x).
\label{eqphiRes}
\end{align}

Since \(|\psi_\eps|\leq1\), \(\phi_\eps\leq0\) and
\(1+\psi_\eps\geq0\), we have
$\exp\left(\alpha_\eps(1+\psi_\eps)\phi_\eps\right)\leq1.$  Combining this with \eqref{add912} and \eqref{eqr0speed},   we obtain that for every \(1<p<4/\beta\),
\begin{equation}\label{eq:phiRhsLp}
\|\Delta^2\phi_\eps\|_{L^p(B_R)}
\leq C(R,p).
\end{equation}

We next derive the local \(L^\infty\)-bound of $\phi_\eps$. Combining Lemma~\ref{lemge} and \eqref{eq:phiRhsLp} with the interior $L^p$ estimate for
\(\Delta\phi_\eps\), we obtain
\[
\|\Delta\phi_\eps\|_{W^{2,p}(B_{3R})}
\leq
C\left(
\|\Delta^2\phi_\eps\|_{L^p(B_{4R})}
+
\|\Delta\phi_\eps\|_{L^1(B_{4R})}
\right)
\leq C(R,p).
\]
Since \(p>1\),  for some \(s>2\), Sobolev embedding gives
 \begin{equation*} 
     \|\Delta\phi_\eps\|_{L^s(B_{3R})}\leq C(R).
 \end{equation*}
Let \(z_\eps\) be the solution of
\[
\begin{cases}
\Delta z_\eps=\Delta\phi_\eps,&\text{in }B_{3R},\\
z_\eps=0,&\text{on }\partial B_{3R}.
\end{cases}
\]
Combining the $W^{2,s}(B_{3R})$ estimate for this Dirichlet problem with the Sobolev embedding theorem, we obtain
\[
\|z_\eps\|_{W^{2,s}(B_{3R})}+\|z_\eps\|_{L^\infty(B_{3R})}\leq C(R).
\]

Since $\phi_\eps-z_\eps$ is harmonic in $B_{3R}$ and  $\phi_\eps\leq 0$ with  $\phi_\eps(0)=0$, and 
\[
\|z_\eps\|_{L^\infty(B_{3R})}\leq C(R),
\]
applying the Harnack inequality to a suitable positive harmonic shift of $\phi_\eps-z_\eps$, we obtain
\[
\|\phi_\eps-z_\eps\|_{L^\infty(B_{2R})}\leq C(R).
\]
Consequently,
\begin{equation}\label{eq:phiLocalLinfty}
\|\phi_\eps\|_{L^\infty(B_{2R})}\leq C(R).
\end{equation}
Thus, for any $1<p<\frac4\beta,$  \eqref{eq:phiRhsLp},
\eqref{eq:phiLocalLinfty} and the interior \(L^p\)-estimate for the bi-Laplacian yield 
\begin{equation*} 
\|\phi_\eps\|_{W^{4,p}(B_R)}\leq C(R,p).
\end{equation*}
Consequently, up to a subsequence, there exists
\(\phi\in W_{\rm loc}^{4,p}(\mathbb R^4)\)  such that  for every \(1<p<4/\beta\) 
\begin{equation*} 
\phi_\eps\rightharpoonup\phi \quad\text{weakly in }W_{\rm loc}^{4,p}(\mathbb R^4)
\end{equation*}
and
\[
\phi_\eps\rightarrow\phi\quad\text{locally uniformly in }\mathbb R^4.
\]
In particular,
\begin{equation*} 
\phi\leq\phi(0)=0.
\end{equation*}

Since translations are continuous in
\(L_{\rm loc}^p(\mathbb R^4)\) for \(p<4/\beta\), we have
\[
\left|
x+r_\eps^{-\frac{4}{4-\beta}}x_\eps
\right|^{-\beta}
\rightarrow
|x+x^*|^{-\beta}
\quad\text{in }L_{\rm loc}^p(\mathbb R^4).
\]
Using also \(\psi_\eps\to1\) locally uniformly and
$\alpha_\eps\rightarrow32\pi^2\left(1-\frac{\beta}{4}\right),$
we may pass to the limit in \eqref{eqphiRes} and obtain in the sense of distributions,
\begin{equation}\label{eq:phiPDEShift}
\Delta^2\phi(x)=
\frac{e^{64\pi^2(1-\frac\beta4)\phi(x)}}{|x+x^*|^\beta} \text{ in } \mathbb R^4,
\phi\leq\phi(0)=0.
\end{equation}

The definition of \(r_\eps\) also gives
\begin{align*}
d_\eps=\int_{\Omega_\eps}
\frac{\psi_\eps^2(x)\exp\left(\alpha_\eps(1+\psi_\eps(x))\phi_\eps(x)\right)}
{\left|x+r_\eps^{-\frac{4}{4-\beta}}x_\eps\right|^\beta}dx.
\end{align*}
Since \(d_\eps\to1\), the Fatou's lemma yields
\begin{equation}\label{eq:phiMassUpper}
\int_{\mathbb R^4}
\frac{
e^{64\pi^2(1-\frac{\beta}{4})\phi(x)}
}{
|x+x^*|^\beta
}dx
\leq1.
\end{equation}

Notice that Fatou's lemma yields only the above inequality. To establish the equality, we will apply the classification theorem established in \cite{AWZ}. Indeed, by Lemma~\ref{lemge} and lower semicontinuity, for every \(R>0\), it holds
\begin{equation}\label{eqlimitDe}
\int_{B_R}|\Delta\phi|dx\leq CR^2.
\end{equation}

From \eqref{eq:phiMassUpper}, we know that \(\phi\) is a finite-volume entire solution of \eqref{eq:phiPDEShift}. After translating \(-x^*\) to the origin, applying Theorem~2.10 in \cite{AWZ}, we see that \(\phi\) must have the form
\[
 \phi=-2(1-\frac{\beta}{4})\log|x|-\sum_{j=1}^4a_j(x_j-x_j^0)^2+c_0+O(|x|^{-\tau}),
 \text{ with } a_j\geq0,
\]
 for some $\tau>0$ and for large $|x|$. Hence, 
\[
 \Delta\phi(x)=-2\sum_{j=1}^4a_j+o(1),
 \text{ as }|x|\to\infty.
\]
Using \eqref{eqlimitDe}, we obtain that
 \(a_j=0\) for every \(j\).
Therefore
\begin{equation}\label{subqug}
    |\phi(x)|=o(|x|^2),
 \text{ as }|x|\to\infty.
\end{equation}

Define
\[
\omega(x)
=16\pi^2\left(1-\frac{\beta}{4}\right)\phi(x)
+\frac14\log\left(\frac{8\pi^2(1-\frac{\beta}{4})}{3}\right).
\]
Then
\[
\Delta^2\omega=6|x+x^*|^{-\beta}e^{4\omega}\quad\text{in }\mathbb R^4.
\]
Moreover, by \eqref{eq:phiMassUpper} and \eqref{subqug}, we have 
\[
\int_{\mathbb R^4}|x+x^*|^{-\beta}e^{4\omega(x)}dx<+\infty,
\]
and \begin{equation*} 
    |\omega(x)|=o(|x|^2),
 \text{ as }|x|\to\infty.
\end{equation*}
Applying \cite[Theorem~1.3]{AWZ} to the translated function
\(x\mapsto\omega(x-x^*)\),  we obtain
\[
6\int_{\RR^4}|x+x^*|^{-\beta}e^{4\omega(x)}dx=16\pi^2\left(1-\frac{\beta}{4}\right).
\]
That is,
\begin{equation}\label{eq:phiMassShift1}
\int_{\RR^4}
\frac{e^{64\pi^2(1-\frac{\beta}{4})\phi(x)}}{|x+x^*|^\beta}dx=1.
\end{equation}
Thus equality holds in \eqref{eq:phiMassUpper}.

Moreover, Theorem~1.3 in \cite{AWZ} implies that $\omega$ is radially symmetric about $-x^*$, and admits the following integral representation:
\[
-\Delta \omega(x)=\frac{3}{2\pi^2}\int_{\RR^4}\frac{|y+x^*|^{-\beta}e^{4\omega(y)}}{|x-y|^2}dy>0.
\]
Writing $r=|x+x^*|$, we therefore have
\[
\bigl(r^3\omega'(r)\bigr)'=r^3\Delta \omega(r)<0.
\]
Since $r^3\omega'(r)\to0$ as $r\to0$, it follows that
$\omega'(r)<0$ for every $r>0$. Hence $-x^*$ is the unique maximum
point of $\omega$, and therefore also of $\phi$. On the other hand, since $x_\eps$ is a maximum point of $u_\eps$,
we have
\[\phi(0)=0=\max_{\RR^4}\phi.\]
Consequently, the unique maximum point of $\phi$ is the origin. 
Hence \begin{equation}\label{xstar}
    x^*=\lim_\eps\frac{x_\eps}{r_\eps^{\frac4{4-\beta}}}=0.
\end{equation}
Combining \eqref{xstar}, \eqref{eq:phiPDEShift} and \eqref{eq:phiMassShift1}, we derive \eqref{eq:phi} and \eqref{eq:phiMassShift}.

Finally, let \(K\Subset\mathbb R^4\setminus\{0\}\). Since $\frac{x_\eps}{r_\eps^{\frac{4}{4-\beta}}}\rightarrow0,$
the shifted weights in \eqref{eqphiRes} are smooth and uniformly bounded on a neighborhood of \(K\). Applying the interior elliptic estimates to the equations satisfied by
\(\psi_\eps\) and \(\phi_\eps\), we obtain
\[
\phi_\eps\rightarrow\phi\quad\text{in }\ C^4(K).
\]
Therefore
\[
\phi_\eps\rightarrow\phi\quad\text{in }\  C_{\rm loc}^4(\mathbb R^4\setminus\{0\}).
\]
The proof is finished.

\end{proof}



\begin{remark}
When $\beta=0$, the limiting profile $\phi$ admits the explicit expression (see \cite{LY-2009})
\[
\phi(x)=-\frac{1}{16\pi^2}\log\left(1+\frac{\pi}{\sqrt{6}}|x|^2 \right).
\]
For $0<\beta<4$, however, no explicit expression for $\phi$ is known (see \cite{AWZ}).
\end{remark}
In the following, to analyze the asymptotic behavior of $u_\eps$ away from the blow up point $0$, we introduce biharmonic truncations and establish their relevant properties. The construction of biharmonic truncations for the Adams inequality without a singular weight was first introduced in  \cite[Lemmas~4.20 and~4.21]{DG2021}. As the proof of the relevant properties is rather long and technical, we defer it to the appendix for the sake of readability.

\begin{proposition}\label{lemDM1}
For every \(M>1\) and all sufficiently small \(\eps>0\), there exist a radius \(\rho_\eps^M>0\), with
\(B_{\rho_\eps^M}(x_\eps)\Subset\Omega\), and a constant
\(C=C(M)>0\) such that
\begin{enumerate}
    \item
    $u_\eps\geq\frac{c_\eps}{M}$ in $B_{\rho_\eps^M}(x_\eps);$
    \item
   $\left|u_\eps-\frac{c_\eps}{M}\right|\leq\frac{C}{c_\eps}$ on $\partial B_{\rho_\eps^M}(x_\eps)$;
    \item for every \(1\leq l\leq3\),
    $|\nabla^l u_\eps|\leq\frac{C}{c_\eps(\rho_\eps^M)^l}$ on $\partial B_{\rho_\eps^M}(x_\eps)$;
    \item $\rho_\eps^M\to0,$  $\frac{\rho_\eps^M}{r_\eps^{\frac{4}{4-\beta}}}\to+\infty$ as $\eps\to0.$
\end{enumerate}
\end{proposition}

Let \(u_\varepsilon^{\rho_k^M} \in C^4\left(\overline{B_{\rho_\varepsilon^M}(x_\varepsilon)}\right)\) be the unique solution of
\[
\begin{cases}
    \Delta^2 u_\varepsilon^{\rho_\varepsilon^M} = 0 & \text{in } B_{\rho_\varepsilon^M}(x_\varepsilon), \\
    \partial_\nu^i u_\varepsilon^{\rho_\varepsilon^M} = \partial_\nu^i u_\varepsilon & \text{on } \partial B_{\rho_\varepsilon^M}(x_\varepsilon), i = 0, 1.
\end{cases}
\]

Define
\[
u_\eps^M=
\begin{cases}
u_\eps^{\rho_\eps^M}
&\text{in }B_{\rho_\eps^M}(x_\eps),\\
u_\eps
&\text{in }\Omega\setminus B_{\rho_\eps^M}(x_\eps).
\end{cases}
\]
Then \(u_\eps^M\in H_0^2(\Omega)\).

\begin{lemma}\label{lemDM2}
For every \(M>1\), we have
\[
u_\eps^M=\frac{c_\eps}{M}+O\left(c_\eps^{-1}\right)
\]
uniformly on \(\overline{B_{\rho_\eps^M}(x_\eps)}\).
\end{lemma}

\begin{proof}
After rescaling \(B_{\rho_\eps^M}(x_\eps)\) onto the unit ball, the standard estimate for biharmonic functions, Proposition~\ref{lemDM1}  and Proposition A.2 in \cite{DG2021} together yield 
\[
\begin{aligned}
\left\| u_\eps^{\rho_\eps^M}-\frac{c_\eps}{M} \right\|_{L^\infty(B_{\rho_\eps^M}(x_\eps))}
\leq & C\bigg(\left\|u_\eps-\frac{c_\eps}{M}\right\|_{L^\infty(\partial B_{\rho_\eps^M}(x_\eps))}
+\rho_\eps^M \|\nabla u_\eps\|_{L^\infty(\partial B_{\rho_\eps^M}(x_\eps))}\bigg)
\leq\frac{C}{c_\eps}.
\end{aligned}
\]
This proves the result.
\end{proof}

\begin{lemma}\label{lemjd}
For every \(M>1\), there holds
\[
\limsup_{\eps\to0}\int_\Omega |\Delta u_\eps^M|^2dx\leq\frac1M.
\]
\end{lemma}

\begin{proof}
Testing the Euler--Lagrange equation \eqref{eqjzfc} with \(u_\eps-u_\eps^M\), we obtain
\[
\begin{aligned}
&\int_{B_{\rho_\eps^M}(x_\eps)}\Delta u_\eps\Delta(u_\eps-u_\eps^M)dx\\
&=\frac{d_\eps}{\lambda_\eps}\int_{B_{\rho_\eps^M}(x_\eps)}\frac{u_\eps e^{\alpha_\eps u_\eps^2}
(u_\eps-u_\eps^M)}{|x|^\beta}dx+\gamma_\eps\int_{B_{\rho_\eps^M}(x_\eps)}u_\eps(u_\eps-u_\eps^M)dx .
\end{aligned}
\]
By Proposition~\ref{lemDM1} and~Lemma \ref{lemDM2}, in $B_{\rho_\eps^M}(x_\eps)$,  we have
\[
u_\eps\geq\frac{c_\eps}{M},
\text{ and }
u_\eps-u_\eps^M\geq-\frac{C}{c_\eps}.
\]
Hence, we get
   \begin{equation*}\begin{aligned}\int_{B_{\rho_\eps^M}(x_\eps)}\frac{u_\eps e^{\alpha_\eps u_\eps^2}(u_\eps-u_\eps^M)}{|x|^\beta}dx &\geq-\frac{C}{c_\eps}\int_{B_{\rho_\eps^M}(x_\eps)}\frac{u_\eps e^{\alpha_\eps u_\eps^2}}{|x|^\beta}dx\\&\geq-\frac{C}{c^2_\eps}\int_{B_{\rho_\eps^M}(x_\eps)}\frac{u_\eps^2e^{\alpha_\eps u_\eps^2}}{|x|^\beta}dx  \geq -\frac{C\lambda_\eps}{c^2_\eps}.
\end{aligned}\end{equation*}

Combining this with the convergence $\Vert   u_\eps \Vert_2\rightarrow0$ and  $$\frac{1}{\left|y+r_\eps^{-\frac{4}{4-\beta}}x_\eps\right|^\beta}\rightarrow\frac{1}{\left|y\right|^\beta} \text{ a.e. in }  \Omega,$$  we have

\begin{equation*}\begin{split}
\int_{B_{\rho_\eps^M}(x_\eps)}\Delta u_\eps\Delta(u_\eps-u_\eps^M)dx
&\geq\frac{d_\eps}{\lambda_\eps}\int_{B_{Rr_\eps^{\frac{4}{4-\beta}}}(x_\eps)}\frac{u_\eps e^{\alpha_\eps u_\eps^2}(u_\eps-u_\eps^M)}{|x|^\beta}dx+o_\eps(1)\\
&\geq \int_{B_R}
\frac{\psi_\eps(y)\left(\psi_\eps(y)-\frac1M+O(c_\eps^{-2})\right)
e^{\alpha_\eps(1+\psi_\eps(y))\phi_\eps(y)}}{\left|y+r_\eps^{-\frac{4}{4-\beta}}x_\eps\right|^\beta}dy+o_\eps(1)\\
&\geq\left(1-\frac1M\right)\int_{B_R}\frac{e^{64\pi^2(1-\frac{\beta}{4})\phi(y)}}{|y|^\beta}dy+o_\eps(1).
\end{split}\end{equation*}
Letting \(R\to+\infty\) and using \eqref{eq:phiMassShift}, we obtain
\begin{equation}\label{addupe10-6}
    \liminf_{\eps\to0}
\int_{B_{\rho_\eps^M}(x_\eps)}
\Delta u_\eps\Delta(u_\eps-u_\eps^M)dx
\geq1-\frac1M.
\end{equation}

Since $u_\eps-u_\eps^M\in H_0^2\bigl(B_{\rho_\eps^M}(x_\eps)\bigr)$ and \(u_\eps^M\) is biharmonic in
\(B_{\rho_\eps^M}(x_\eps)\), then
\[
\int_{B_{\rho_\eps^M}(x_\eps)}\Delta u_\eps^M\Delta(u_\eps-u_\eps^M)dx=0.
\]
This together with the fact that 
\(\|\Delta u_\eps\|_2^2=1\), we have
\[
\begin{aligned}
\int_\Omega|\Delta u_\eps^M|^2dx
&=1-\int_{B_{\rho_\eps^M}(x_\eps)}
\left(|\Delta u_\eps|^2-|\Delta u_\eps^M|^2\right)dx\\
&=1-\int_{B_{\rho_\eps^M}(x_\eps)}
\Delta u_\eps\Delta(u_\eps-u_\eps^M)dx.
\end{aligned}
\]
It then follows from this and \eqref{addupe10-6} that
\[
\limsup_{\eps\to0}\int_\Omega|\Delta u_\eps^M|^2dx\leq\frac1M.
\] Therefore, the proof is finished.
\end{proof}

\begin{lemma}\label{lemexp}
There holds
\[
\begin{aligned}
\lim_{\eps\to0}\int_\Omega\frac{e^{\alpha_\eps u_\eps^2}}{|x|^\beta}dx
&=\int_\Omega\frac{1}{|x|^\beta}dx+\lim_{R\to+\infty}\lim_{\eps\to0}\int_{B_{Rr_\eps^{\frac{4}{4-\beta}}}(x_\eps)}\frac{e^{\alpha_\eps u_\eps^2}}{|x|^\beta}dx\\
&=\int_\Omega \frac{1}{|x|^\beta}dx+\lim_{\eps\to0}\frac{\lambda_\eps}{c_\eps^2}.
\end{aligned}
\]
Consequently,
\[
\frac{\lambda_\eps}{c_\eps}\to+\infty,
\ \ and \ \
\sup_\eps\frac{c_\eps^2}{\lambda_\eps}<+\infty.
\]
\end{lemma}

\begin{proof}
By \eqref{eq4.4} and the definition of  $u_\eps^M$, we see that  \(u_\eps^M\to0\) almost everywhere in \(\Omega\) and
\(\alpha_\eps\) bounded, hence
\[
e^{\alpha_\eps(u_\eps^M)^2}-1\to0 \quad\text{a.e. in }\Omega.
\]
On the other hand, for any \(M>1\), choosing \(1<p<M\), it follows from \eqref{eq4.4} and Lemma~\ref{lemjd} that
\[
\limsup_{\eps\to0} p\alpha_\eps\|\Delta u_\eps^M\|_2^2
\leq \frac{p}{M} 32\pi^2\left(1-\frac{\beta}{4}\right)
< 32\pi^2\left(1-\frac{\beta}{4}\right).
\]
Then the singular Adams inequality \eqref{eqsadm} implies that 
$e^{\alpha_\eps (u_\eps^M)^2}$ is bounded in $L^p(\Omega,|x|^{-\beta}dx)$ with some $p>1$.  Hence, using Vitali's convergence theorem, we have
\[
\lim_{\eps\rightarrow 0}\int_\Omega\frac{e^{\alpha_\eps(u_\eps^M)^2}-1}{|x|^\beta}dx=0, 
\]
and then
\[
\lim_{\eps\rightarrow 0}\int_{\Omega\setminus B_{\rho_\eps^M}(x_\eps)}\frac{e^{\alpha_\eps u_\eps^2}-1}{|x|^\beta}dx=0.
\]
Together with the fact that $u_\eps>\frac{c_\eps}{M}$ on $B_{\rho_\eps^M}(x_\eps)$, this yields 
\begin{equation}\begin{split}\label{eqxiaoyu}
\int_\Omega\frac{e^{\alpha_\eps u_\eps^2}-1}{|x|^\beta}dx
&=\int_{B_{\rho_\eps^M}(x_\eps)}\frac{e^{\alpha_\eps u_\eps^2}}{|x|^\beta}dx+o_\eps(1)\\
&\leq\frac{M^2}{c_\eps^2}\int_\Omega\frac{u_\eps^2e^{\alpha_\eps u_\eps^2}}{|x|^\beta}dx+o_\eps(1)
=M^2\frac{\lambda_\eps}{c_\eps^2}+o_\eps(1).
\end{split}\end{equation}

On the other hand, performing the change of variables  $x=x_\eps+r_\eps^{\frac{4}{4-\beta}}y$  and using  Lemma~\ref{lemxe}, \(d_\eps\to1\) and
\eqref{eq:phiMassShift},  we obtain
\begin{equation}\begin{split}\label{eqdayu}
\int_\Omega\frac{e^{\alpha_\eps u_\eps^2}-1}{|x|^\beta}dx
&\geq\lim_{R\to+\infty}\lim_{\eps\to0}\int_{B_{Rr_\eps^{\frac{4}{4-\beta}}}(x_\eps)}
\frac{e^{\alpha_\eps u_\eps^2}-1}{|x|^\beta}dx\\
&=\lim_{R\to+\infty}\lim_{\eps\to0}\frac{\lambda_\eps}{d_\eps c_\eps^2}\int_{B_R}
\frac{e^{\alpha_\eps(1+\psi_\eps(y))\phi_\eps(y)}}{
\left|y+r_\eps^{-\frac{4}{4-\beta}}x_\eps\right|^\beta}dy\\
&\geq\lim_{R\to+\infty}\left(\lim_{\eps\to0}\frac{\lambda_\eps}{c_\eps^2}\right)
\int_{B_R}\frac{e^{64\pi^2(1-\frac{\beta}{4})\phi(y)}}{|y|^\beta}dy\\
&=\lim_{\eps\to0}\frac{\lambda_\eps}{c_\eps^2}.
\end{split}\end{equation}
Combining \eqref{eqxiaoyu} and \eqref{eqdayu}, we get
\[
\lim_{\eps\to0}\frac{\lambda_\eps}{c_\eps^2}
\leq\lim_{\eps\to0}\int_\Omega\frac{e^{\alpha_\eps u_\eps^2}-1}{|x|^\beta}dx
\leq M^2\lim_{\eps\to0}\frac{\lambda_\eps}{c_\eps^2}.
\]
Letting \(M\rightarrow1\), we have
\[
\begin{aligned}
\lim_{\eps\to0}\int_\Omega\frac{e^{\alpha_\eps u_\eps^2}-1}{|x|^\beta}dx
&
=\lim_{\eps\to0}\frac{\lambda_\eps}{c_\eps^2}.
\end{aligned}
\]
Consequently, we obtain
\[
\begin{aligned}
\lim_{\eps\to0}\int_\Omega\frac{e^{\alpha_\eps u_\eps^2}}{|x|^\beta}dx
&=\int_\Omega \frac{1}{|x|^\beta}dx+\lim_{\eps\to0}\frac{\lambda_\eps}{c_\eps^2}.
\end{aligned}
\]

Finally, since $S_{\Omega,q,\beta}>\int_\Omega\frac{dx}{|x|^\beta}$,  we conclude that
\[
\sup_\eps\frac{c_\eps^2}{\lambda_\eps}<+\infty
\ {\text {and }}\
\frac{\lambda_\eps}{c_\eps}+\infty,
\]
as $\eps\to 0$. The proof is finished.
\end{proof}


To proceed with the analysis of the asymptotic behavior of $u_\eps$  away from the blow-up point $0$, we introduce the following quantities:
\begin{equation*} 
\begin{cases}
\displaystyle
b_\eps = \lambda_\eps \big/ \int_\Omega |u_\eps| \frac{\exp({\alpha_\eps u_\eps^2})}{|x|^{\beta}}dx,\\
\displaystyle
\tau = \lim_{\eps\to 0}\frac{c_\eps}{b_\eps},\\
\displaystyle
\sigma = \lim_{\eps\to 0}\frac{\int_\Omega \frac{u_\eps \exp({\alpha_\eps u_\eps^2})}{|x|^{\beta}}dx}{\int_\Omega \frac{|u_\eps| \exp({\alpha_\eps u_\eps^2})}{|x|^{\beta}}dx}.
\end{cases}
\end{equation*}

Obviously, we have $\tau\ge 1$ and $|\sigma|\le 1$. In what follows, we will show that   $\sigma= 1$ and  $\tau=1$.

\begin{lemma} 
There holds $\sigma=1.$
\end{lemma}

\begin{proof}

For any \(M>1\), by Lemma \ref{lemjd} and the singular Adams inequality \eqref{eqsadm}, there exists some \(p>1\) such that
\[
\limsup_{\eps}
\int_\Omega
\frac{e^{p\alpha_\eps(u_\eps^M)^2}}{|x|^\beta}dx
<+\infty.
\]
Choosing \(1<q<p\), for every \(t\in\mathbb R\), it holds
\[
\left(
\bigl(|t|+t^2\bigr)e^{\alpha_\eps t^2}\right)^q\leq C e^{p\alpha_\eps t^2}.
\]
Consequently, we have
\begin{equation}\label{addboundp108}
  \sup_\eps
\int_\Omega\frac{\left(\bigl(|u_\eps^M|+(u_\eps^M)^2\bigr)e^{\alpha_\eps(u_\eps^M)^2}\right)^q}{|x|^\beta}dx<+\infty.  
\end{equation}

Since \(\rho_\eps^M\to0\) and \(x_\eps\to0\), hence for every \(x\in\Omega\setminus\{0\}\) and all sufficiently small \(\eps\),
\[
x\notin B_{\rho_\eps^M}(x_\eps),
u_\eps^M(x)=u_\eps(x)\to0.
\]
It follows that
\begin{equation}\label{addae108}
    \bigl(|u_\eps^M|+(u_\eps^M)^2\bigr)
e^{\alpha_\eps(u_\eps^M)^2}
\to0\text{ a.e. in }\Omega.
\end{equation}

Combining \eqref{addboundp108} and \eqref{addae108}, by Vitali's convergence theorem, we have
\[
\lim_{\eps\rightarrow 0}\int_\Omega\frac{\bigl(|u_\eps^M|+(u_\eps^M)^2\bigr)e^{\alpha_\eps(u_\eps^M)^2}}{|x|^\beta}dx=0.
\]
Notice that  \(u_\eps^M=u_\eps\) in
\(\Omega\setminus B_{\rho_\eps^M}(x_\eps)\), hence we have
\begin{equation}\begin{split}\label{eqtrunc}
\int_{\Omega\setminus B_{\rho_\eps^M}(x_\eps)}\frac{\bigl(|u_\eps|+u_\eps^2\bigr)e^{\alpha_\eps u_\eps^2}}{|x|^\beta}dx
&=\int_{\Omega\setminus B_{\rho_\eps^M}(x_\eps)}\frac{\bigl(|u_\eps^M|+(u_\eps^M)^2\bigr)e^{\alpha_\eps(u_\eps^M)^2}}{|x|^\beta}dx\\
&\leq\int_\Omega\frac{\bigl(|u_\eps^M|+(u_\eps^M)^2\bigr)e^{\alpha_\eps(u_\eps^M)^2}}{|x|^\beta}dx
=o_\eps(1).
\end{split}\end{equation}

Since \(u_\eps>0\) in \(B_{\rho_\eps^M}(x_\eps)\), it follows from \eqref{eqtrunc} that
\begin{equation}\label{adduni108}
\begin{aligned}
\int_\Omega\frac{u_\eps e^{\alpha_\eps u_\eps^2}}{|x|^\beta}dx
&=\int_{B_{\rho_\eps^M}(x_\eps)}\frac{|u_\eps|e^{\alpha_\eps u_\eps^2}}{|x|^\beta}dx+o_\eps(1),\\
\int_\Omega\frac{|u_\eps|e^{\alpha_\eps u_\eps^2}}{|x|^\beta}dx
&=\int_{B_{\rho_\eps^M}(x_\eps)}\frac{|u_\eps|e^{\alpha_\eps u_\eps^2}}{|x|^\beta}dx+o_\eps(1).
\end{aligned}
\end{equation}
Moreover, by Proposition~\ref{lemDM1} and~Lemma \ref{lemexp}, 
\begin{equation}\label{addloci108}
    \begin{aligned}
\int_{B_{\rho_\eps^M}(x_\eps)}\frac{|u_\eps|e^{\alpha_\eps u_\eps^2}}{|x|^\beta}dx&\geq\frac{c_\eps}{M}\int_{B_{\rho_\eps^M}(x_\eps)}\frac{e^{\alpha_\eps u_\eps^2}}{|x|^\beta}dx\\&\geq\frac{c_\eps}{M}\int_{B_{Rr_\eps^{\frac{4}{4-\beta}}}(x_\eps)}\frac{e^{\alpha_\eps u_\eps^2}}{|x|^\beta}dx\geq\frac{\lambda_\eps}{Mc_\eps}\rightarrow\infty.
\end{aligned}
\end{equation}
Combining \eqref{adduni108} and \eqref{addloci108}, we have
\[
\sigma = \lim_{\eps\to 0}\frac{\int_\Omega \frac{u_\eps \exp({\alpha_\eps u_\eps^2})}{|x|^{\beta}}dx}{\int_\Omega \frac{|u_\eps| \exp({\alpha_\eps u_\eps^2})}{|x|^{\beta}}dx}=1.
\]
\end{proof}

\begin{lemma}\label{lemtau}
There holds $\tau=1.$
\end{lemma}

\begin{proof}
For any  \(M>1\), it follows from \eqref{eqtrunc} that
\[
\lambda_\eps=\int_{B_{\rho_\eps^M}(x_\eps)}\frac{u_\eps^2e^{\alpha_\eps u_\eps^2}}{|x|^\beta}dx+o_\eps(1)
\]
and
\[
\int_\Omega\frac{|u_\eps|e^{\alpha_\eps u_\eps^2}}{|x|^\beta}dx=\int_{B_{\rho_\eps^M}(x_\eps)}\frac{u_\eps e^{\alpha_\eps u_\eps^2}}{|x|^\beta}dx+o_\eps(1).
\]
Combining the above relations with the fact that 
$
c_\eps/M\leq u_\eps\leq c_\eps\ \text{in }B_{\rho_\eps^M}(x_\eps),
$ we deduce that
\[
\lambda_\eps\geq\frac{c_\eps^2}{M^2}
\int_{B_{\rho_\eps^M}(x_\eps)}\frac{e^{\alpha_\eps u_\eps^2}}{|x|^\beta}dx+o_\eps(1)
\]
and
\[
\int_\Omega\frac{|u_\eps|e^{\alpha_\eps u_\eps^2}}{|x|^\beta}dx\leq
c_\eps\int_{B_{\rho_\eps^M}(x_\eps)}\frac{e^{\alpha_\eps u_\eps^2}}{|x|^\beta}dx+o_\eps(1).
\]
By Lemma \ref{lemexp} and the above estimates, we obtain
\[
b_\eps=\frac{\lambda_\eps}
{\int_\Omega\frac{|u_\eps|e^{\alpha_\eps u_\eps^2}}{|x|^\beta}dx}
\geq\frac{c_\eps}{M^2}\bigl(1+o_\eps(1)\bigr).
\]
On the other hand, since \(u_\eps^2\leq c_\eps|u_\eps|\) in \(\Omega\), and hence $b_\eps\leq c_\eps.$ Therefore, we have 
 \[1\leq\frac{c_\eps}{b_\eps}\leq M^2+o_\eps(1).\]
Letting first \(\eps\to0\) and then \(M\rightarrow1\), we conclude that $\tau=1.$
\end{proof}


We now investigate the asymptotic behavior of $u_\eps$ away from the blow-up point $0$. To this end, we first recall that the Green function $G(x,y)$ associated with the biharmonic operator
$\Delta^2$ with Dirichlet boundary conditions is defined by
\begin{equation*}
\begin{cases}\Delta^2G(x,y)=\delta_x(y),&  \text{ in }\Omega,\\
G(x,y)=\frac{\partial G(x,y)}{\partial\nu}=0,
 & \text{ on }\partial\Omega.
\end{cases}
\end{equation*}
Every sufficiently regular solution $u\in H_0^2(\Omega)$  of $\Delta^2u=f$ can be represented as
\[u(x)=\int_\Omega G(x,y)f(y)dy.\]
We shall use the following useful estimates for $G(x,y)$ (see,  \cite{DS2004}): there exists a constant $C>0$ such that  for all $x,y\in\Omega$ with $x\neq y$,
\begin{equation}\label{eq4.2}
\begin{aligned}
|G(x,y)|
&\leq C\log\left(2+\frac{1}{|x-y|}\right),\
|\nabla^i G(x,y)|
&\leq C|x-y|^{-i},
  i\geq 1.
\end{aligned}
\end{equation}

Using these estimates of the Green function $G(x,y)$, we can obtain the following result.
\begin{lemma}\label{lem:weighted43}
For any $1<r<2$, the sequence $\{c_\eps u_\eps\}$ is
bounded in $W_0^{2,r}(\Omega)$. In particular, there exists a constant $C$ depending only on $\Omega$, $\lambda(\Omega)$ and $\alpha_0$ such that $\|c_\eps u_\eps\|_{W_0^{2,r}(\Omega)}\le C$ uniformly for $a_1\in[0,\alpha_0]$ with $\alpha_0<\lambda(\Omega)$.
\end{lemma}

\begin{proof}
Let $\zeta_\eps$ solve the following biharmonic problem
\[
\begin{cases}
\displaystyle
\Delta^2\zeta_\eps=\frac{d_\eps b_\eps}{\lambda_\eps}\frac{u_\eps e^{\alpha_\eps u_\eps^2}}{|x|^\beta}
&\text{in }\Omega,\\
\zeta_\eps=\partial_\nu \zeta_\eps=0
&\text{on }\partial\Omega.
\end{cases}
\]
By the definition of $b_\eps$, we have
\begin{equation}\label{eq:vMass}
\frac{d_\eps b_\eps}{\lambda_\eps}
\int_\Omega
\frac{|u_\eps|e^{\alpha_\eps u_\eps^2}}{|x|^\beta}dx
=d_\eps\leq C.
\end{equation}
For $i=1,2$, using Green’s representation formula and the Green-function estimate \eqref{eq4.2}, we have
\[
|\nabla^i\zeta_\eps(x)|\leq C\frac{d_\eps b_\eps}{\lambda_\eps}
\int_\Omega\frac{|u_\eps(y)|e^{\alpha_\eps u_\eps^2(y)}}{|y|^\beta|x-y|^i}dy.
\]
By H\"older's inequality and \eqref{eq:vMass}, for any $1<r<2$,
\begin{align*}
|\nabla^i\zeta_\eps(x)|^r
&\leq C_r\left(\frac{d_\eps b_\eps}{\lambda_\eps}
\int_\Omega\frac{|u_\eps(y)|e^{\alpha_\eps u_\eps^2(y)}}
{|y|^\beta}dy\right)^{r-1}\frac{d_\eps b_\eps}{\lambda_\eps}\int_\Omega\frac{|u_\eps(y)|e^{\alpha_\eps u_\eps^2(y)}}{|y|^\beta|x-y|^{ir}}dy \\
&\leq C_r\frac{d_\eps b_\eps}{\lambda_\eps}\int_\Omega\frac{|u_\eps(y)|e^{\alpha_\eps u_\eps^2(y)}}{|y|^\beta|x-y|^{ir}}dy.
\end{align*}
Since $ir<4$ for $i=1,2$, using the Fubini's theorem, we obtain
\begin{align*}
\int_\Omega|\nabla^i\zeta_\eps(x)|^rdx
&\leq C_r\frac{d_\eps b_\eps}{\lambda_\eps}\int_\Omega\frac{|u_\eps(y)|e^{\alpha_\eps u_\eps^2(y)}}{|y|^\beta}
\left(\int_\Omega|x-y|^{-ir}dx\right)dy\leq C_r.
\end{align*}
From this and the Poincar\'e inequality, it follows that 
\begin{equation}\label{eq:vHr}
\|\zeta_\eps\|_{W_0^{2,r}(\Omega)}\leq C_r, \text{ for any }
 1<r<2.
\end{equation}

Set $w_\eps=b_\eps u_\eps-\zeta_\eps.$ 
Then $w_\eps$ satisfies
\[
\begin{cases}
\Delta^2w_\eps=\gamma_\eps(w_\eps+\zeta_\eps), \text{ in }\Omega\\
w_\eps=\partial_\nu w_\eps=0,
\quad\text{ on }\partial\Omega.
\end{cases}\]
Testing by $w_\eps$, using the definition of $\lambda(\Omega)$
and Poincar\'e inequality, we obtain
\begin{align*}
\|\Delta w_\eps\|_2^2
&=\gamma_\eps\|w_\eps\|_2^2+\gamma_\eps\int_\Omega \zeta_\eps w_\eps dx\leq\frac{\gamma_\eps}{\lambda(\Omega)}
\|\Delta w_\eps\|_2^2+C\|\zeta_\eps\|_2\|\Delta w_\eps\|_2.
\end{align*}
Since $\gamma_\eps\to a_1<\lambda(\Omega)$, for sufficiently small $\eps$, we obtain 
\[
\left(1-\frac{\gamma_\eps}{\lambda(\Omega)}\right)
\|\Delta w_\eps\|_2^2
\leq C\|\Delta w_\eps\|_2,
\]
which implies $\|\Delta w_\eps\|_2\leq C$. Using the Sobolev embedding, we immediately have  $\|w_\eps\|_{W_0^{2,r}(\Omega)}\leq C$ for any $1<r<2$. Combining this with \eqref{eq:vHr}, we conclude that
\[
\|b_\eps u_\eps\|_{W_0^{2,r}(\Omega)}
\leq\|w_\eps\|_{W_0^{2,r}(\Omega)}+\|\zeta_\eps\|_{W_0^{2,r}(\Omega)}
\leq C_r,
\]
where $C_r$ is independent of $\eps$ and can be chosen
uniformly for $a_1\in[0,\alpha_0]$.

Finally, from Lemma \ref{lemtau}, we have  $\lim\limits_{\eps\to 0}\frac{c_\eps}{b_\eps}=1$. Therefore, for any $1<r<2$, we have $\|c_\eps u_\eps\|_{W_0^{2,r}(\Omega)}\le C_r$ uniformly for $a_1\in[0,\alpha_0]$ with $\alpha_0<\lambda(\Omega)$. 
\end{proof}


Next,  we prove that $c_\eps u_\eps$ converges to some Green function.

\begin{lemma}\label{lem:weighted-source}
For every \(\varphi\in C_c^\infty(\Omega)\), there holds
\[
\lim_{\eps\to0}\int_\Omega\varphi(x)\frac{d_\eps c_\eps}{\lambda_\eps}\frac{u_\eps e^{\alpha_\eps u_\eps^2}}{|x|^\beta}dx=\varphi(0).
\]
\end{lemma}

\begin{proof}

Fix \(M>1\) and \(R>0\). Since $\frac{\rho_\eps^M}{r_\eps^{\frac4{4-\beta}}}\to+\infty,$ we have $B_{Rr_\eps^{\frac4{4-\beta}}}(x_\eps)\subset B_{\rho_\eps^M}(x_\eps)$ for sufficiently small \(\eps\). We split the integral as follows
\begin{equation*}\begin{split} 
\int_\Omega
\varphi(x)\frac{d_\eps c_\eps}{\lambda_\eps}
\frac{u_\eps e^{\alpha_\eps u_\eps^2}}{|x|^\beta}dx&=\left(\int_{B_{\rho_\eps^M}\setminus B_{Rr_\eps^{\frac4{4-\beta}}}}+\int_{B_{Rr_\eps^{\frac4{4-\beta}}}}+\int_{\Omega\setminus B_{\rho_\eps^M}}\right)
\varphi(x)\frac{d_\eps c_\eps}{\lambda_\eps}
\frac{u_\eps e^{\alpha_\eps u_\eps^2}}{|x|^\beta}dx\\&=:I_1+I_2+I_3.
\end{split}\end{equation*}
Since \(u_\eps\geq c_\eps/M\) in \(B_{\rho_\eps^M}(x_\eps)\), we have
\(c_\eps u_\eps\leq M u_\eps^2\). Hence, using the change of variables, we have
\begin{equation*}\begin{split}
I_1&\leq M\|\varphi\|_\infty\frac{d_\eps}{\lambda_\eps}\int_{B_{\rho_\eps^M}(x_\eps)
\setminus B_{Rr_\eps^{\frac4{4-\beta}}}(x_\eps)}\frac{u_\eps^2e^{\alpha_\eps u_\eps^2}}{|x|^\beta}dx\\
&\leq M\|\varphi\|_\infty\left(d_\eps-\frac{d_\eps}{\lambda_\eps}\int_{B_{Rr_\eps^{\frac4{4-\beta}}}(x_\eps)}\frac{u_\eps^2e^{\alpha_\eps u_\eps^2}}{|x|^\beta}dx\right)\\
&\leq M\|\varphi\|_\infty \left(d_\eps- \int_{B_R}
\frac{\psi_\eps^2(y)e^{\alpha_\eps(1+\psi_\eps(y))\phi_\eps(y)}}{\left|y+r_\eps^{-\frac4{4-\beta}}x_\eps\right|^\beta}dy \right).
\end{split}\end{equation*}
From Lemma~\ref{lemxe}, using the fact that  $$\frac{1}{\left|y+r_\eps^{-\frac{4}{4-\beta}}x_\eps\right|^\beta}\rightarrow\frac{1}{\left|y\right|^\beta} \text{ a.e. in }  \Omega,$$ and Fatou's lemma,  we obtain
\begin{equation}\begin{split}\label{eq091402}
\lim_{R\to+\infty}\limsup_{\eps\to0}I_1
\leq M\|\varphi\|_\infty\left(1-\int_{\mathbb R^4}\frac{e^{64\pi^2(1-\frac\beta4)\phi(y)}}{|y|^\beta}dy\right)=0.
\end{split}\end{equation}
For $I_2$, performing the same change of variables and using Lemma~\ref{lemxe},  we get the convergence $$\frac{1}{\left|y+r_\eps^{-\frac{4}{4-\beta}}x_\eps\right|^\beta}\rightarrow\frac{1}{\left|y\right|^\beta} \text{ in }  L^p(\Omega), \text{ for some }p>1.$$ 
By Vitali's convergence theorem, we obtain
\begin{equation}\begin{split}\label{eq091403}
\lim_{R\to+\infty}\lim_{\eps\to0}I_2&=\lim_{R\to+\infty}\lim_{\eps\to0}\int_{B_R}\varphi\left(x_\eps+r_\eps^{\frac4{4-\beta}}y\right)\frac{\psi_\eps(y)
e^{\alpha_\eps(1+\psi_\eps(y))\phi_\eps(y)}}{\left|y+r_\eps^{-\frac4{4-\beta}}x_\eps\right|^\beta}dy\\
&=\varphi(0)\int_{\mathbb R^4}\frac{e^{64\pi^2(1-\frac\beta4)\phi(y)}}{|y|^\beta}dy
=\varphi(0).
\end{split}\end{equation}
Finally, by \eqref{eqtrunc} and the fact that \(\lambda_\eps/c_\eps\to+\infty\)(see Lemma \ref{lemexp}), we have
\begin{equation}\begin{split}\label{eq091404}
I_3&\leq\|\varphi\|_\infty\frac{d_\eps c_\eps}{\lambda_\eps}\int_{\Omega\setminus B_{\rho_\eps^M}(x_\eps)}
\frac{|u_\eps|e^{\alpha_\eps u_\eps^2}}{|x|^\beta}dx\to0.
\end{split}\end{equation}

Combining the estimates \eqref{eq091402}, \eqref{eq091403} and \eqref{eq091404}, we completed the proof.
\end{proof}

\begin{lemma}\label{lem4.4}
For any $1<r<2$,
$ c_\eps u_\eps \rightharpoonup G_{a_1}(\cdot,0) \text{ in } W_0^{2,r}(\Omega), $
where $G_{a_1}$ satisfies
\begin{equation*} 
\begin{cases}
\Delta^2 G_{a_1}=\delta_0+a_1G_{a_1} & \text{in }\Omega,\\
G_{a_1}=\dfrac{\partial G_{a_1}}{\partial\nu}=0 & \text{on }\partial\Omega.
\end{cases}
\end{equation*}
Moreover, $c_\eps u_\eps\to G_{a_1}(\cdot,0) \text{ in } C_{\rm loc}^4(\overline{\Omega}\setminus\{0\})$. 
The function $G_{a_1}$ admits the expansion
\begin{equation}\label{eq:GreenExpansionSigma}
G_{a_1}(x,0) =-\frac{1}{8\pi^2}\log|x|+A_0+H(x),
\end{equation}
where $A_0$ is a constant depending on $a_1$. Here  $H(x)\in C^3(\overline{\Omega})$ satisfies $H(0)=0$.
\end{lemma}

\begin{proof}

From the Euler--Lagrange equation \eqref{eqjzfc}, we obtain
\begin{equation}\label{eq:scaled-Euler}
\Delta^2(c_\eps u_\eps)
=\frac{d_\eps c_\eps}{\lambda_\eps}\frac{u_\eps e^{\alpha_\eps u_\eps^2}}{|x|^\beta}+\gamma_\eps c_\eps u_\eps
\quad\text{in }\Omega.
\end{equation}
By Lemma~\ref{lem:weighted-source}, we have
\begin{equation}\label{adddela108}
   \frac{d_\eps c_\eps}{\lambda_\eps}\frac{u_\eps e^{\alpha_\eps u_\eps^2}}{|x|^\beta}\rightharpoonup\delta_0, 
\end{equation}
in the sense of distributions. 

On the other hand, Lemma \ref{lem:weighted43} shows that for every \(1<r<2\)  there exists \(G_{a_1}(\cdot,0)\in W_0^{2,r}(\Omega)\) such that
\begin{equation}\label{addweakc108}
c_\eps u_\eps\rightharpoonup G_{a_1}(\cdot,0)  \text{ in }W_0^{2,r}(\Omega).
\end{equation}
 
 Since \(\gamma_\eps\to a_1\), combining \eqref{adddela108} and \eqref{addweakc108}, passing to the limit in \eqref{eq:scaled-Euler}, we obtain
\[
\Delta^2G_{a_1}=\delta_0+a_1G_{a_1} \text{ in }\Omega,
\]
with $G_{a_1}=\frac{\partial G_{a_1}}{\partial\nu}=0  \text{ on }\partial\Omega.$

For any fixed $r>0$ such that $B_{2r}(0)\Subset\Omega$, using the singular Adams inequality \eqref{eqsadm} and a cut-off function argument,  we see that
 $\frac{e^{\alpha_\eps u_\eps^2}}{|x|^\beta}$  is bounded in $L^s(\Omega\setminus B_r(0))$ for some $s>1$. This, together with the boundedness of $c_\eps/\lambda_\eps\leq C$ and  the standard  regularity theory for biharmonic operator $\Delta^2$, gives $c_\eps u_\eps\rightarrow G_{a_1}(\cdot,0)$
{in }$C_{\mathrm{loc}}^4(\overline\Omega\setminus\{0\}).$

To prove \eqref{eq:GreenExpansionSigma}, we choose
$\eta\in C_0^4(\Omega)$ such that $ 0\leq\eta\leq1 $ on $\Omega$, 
$$\eta\equiv1 \text{ on }B_r(0), \text{and}~\eta\equiv0\text{ on }\Omega\setminus B_{2r}(0). $$
Define
\[g(x)=G_{a_1}(x,0)+\frac{1}{8\pi^2}\eta(x)\log|x|.\]
Then $g$ satisfies
\[
\Delta^2g=f \ \text{in }\ \Omega,  g=\frac{\partial g}{\partial\nu}=0\ \text{on }\ \partial\Omega,\]
 where
\begin{equation*} 
\begin{aligned}
f(x) &=-\frac{\sigma}{8\pi^2}\bigl(\Delta^2\eta\log|x|+2\nabla\Delta\eta\cdot\nabla\log|x|\\
&+2\Delta\eta\Delta\log|x|+2\Delta(\nabla\eta\cdot\nabla\log|x|)+2\nabla\eta\cdot\nabla\Delta\log|x|\bigr)+a_1G_{a_1}(x,0).
\end{aligned}
\end{equation*}

Since
$G_{a_1}\in W_0^{2,s}(\Omega)$ for every $1<s<2$, and hence
$G_{a_1}\in L^q(\Omega)$ for every $q>1$. Therefore,
$f\in L^q(\Omega)$ for every $q>1$. Then by the
standard regularity theory,  we have $g\in C^3(\overline\Omega).$

Let $A_0=g(0)$ and define
\[
H(x)=g(x)-g(0)+\frac{1}{8\pi^2}(1-\eta(x))\log|x|.
\]
Since $1-\eta$ vanishes near the origin,
$H\in C^3(\overline\Omega)$, $\psi(0)=0$.
It follows that
\[G_{a_1}(x,0)=-\frac{1}{8\pi^2}\log|x|+A_0+H(x). \]
This completes the proof.
\end{proof}


\section{Neck analysis and the sharp upper bound for $S_{\Omega,q,\beta}$ in the blow-up case}

In this section, we will study the sharp upper bound for $S_{\Omega,q,\beta}$ when the blow-up occurs, and complete the proof  of Theorem~\ref{thm1}. To this end, we employ a direct comparison argument involving biharmonic functions on suitably chosen annular domains, following the approach introduced by DelaTorre and Mancini \cite{DG2021} (see also \cite{ChLuZhu2026}).

‌The main result of this section is as follows:
\begin{proposition}
If blow-up occurs, then the supremum
\begin{equation}\label{eqly5.23t}
\begin{aligned}
S_{\Omega,q,\beta}\leq\int_\Omega\frac{1}{|x|^{\beta}}dx+\exp\Big(32\pi^2\big(1-\frac \beta4\big)\left(A_0-C_0-J_0\right)\Big),
\end{aligned}
\end{equation}
where $A_0$ is the constant appearing in \eqref{eq:GreenExpansionSigma}, and
\begin{equation*}\begin{cases}
\displaystyle
C_0=\frac{1}{8\pi^2}\int_{\RR^4}
\frac{\log|y|e^{64\pi^2(1-\frac{\beta}{4})\phi(y)}}{|y|^\beta}dy,\\
\displaystyle
J_0=\int_{\mathbb R^4}\frac{\phi(x)e^{64\pi^2(1-\frac \beta4)\phi(x)}}{|x|^{\beta}}dx,
\end{cases}\end{equation*}
here $\phi$ is the limiting function appearing in Lemma \ref{lemxe}.
\end{proposition}

\begin{remark}
When $\beta=0$, the limiting profile is explicitly given by
\[
\phi(x)=-\frac{1}{16\pi^2}\log\left( 1+\frac{\pi}{\sqrt6}|x|^2\right).
\]
Through direct calculation, we have
\[
C_0=-\frac{1}{16\pi^2}\log\frac{\pi}{\sqrt6}.
\]
Moreover,
\[
\begin{aligned}
J_0&=\int_{\RR^4}\phi(x)e^{64\pi^2\phi(x)}dx=
-\frac{3}{8\pi^2}\int_0^\infty\frac{t\log(1+t)}{(1+t)^4}dt=-\frac{5}{96\pi^2}.
\end{aligned}
\]
Substituting these values of $C_0$ and $J_0$ into the preceding
 upper bound \eqref{eqly5.23t}, we recover
\[
S_{\Omega,q,0}\leq|\Omega|+\frac{\pi^2}{6}\exp\left(\frac53+32\pi^2A_p\right)
\]
when the blow-up occurs, which coincides with the upper bound obtained in Section~5 of
\cite{LY-2009}.
\end{remark}

\begin{proof}[Proof of Proposition \ref{eqly5.23t}]

Fix \(R>1\) and \(0<\delta<\frac12d(0,\partial\Omega)\). Set
\[
 A_\eps(R,\delta) :=B_\delta(x_\eps)\setminus
 B_{Rr_\eps^{\frac4{4-\beta}}}(x_\eps).
\]
 Since \(x_\eps\to0\) and \(r_\eps\to0\), the set $A_\eps(R,\delta)$ is contained in \(\Omega\) for all sufficiently small \(\eps\). To the sharp upper bound for $S_{\Omega,q,\beta}$ when the blow-up occurs, we will compare the Dirichlet energies of \(u_\eps\) and the biharmonic comparison function \begin{equation*}
 W_\eps(x):=-\frac1{8\pi^2c_\eps}\log|x-x_\eps| 
\end{equation*} on $A_\eps(R,\delta)$.

Using the definition of \(r_\eps\), we obtain the exact logarithmic identity
\begin{equation}\label{eqreps}
 \log r_\eps^{\frac4{4-\beta}}
 =\frac1{4-\beta}
 \left(
  \log\frac{\lambda_\eps}{d_\eps c_\eps^2}
  -\alpha_\eps c_\eps^2
 \right).
\end{equation}
By Lemma~\ref{lemxe}, on \(\partial B_{Rr_\eps^{4/(4-\beta)}}(x_\eps)\) we have 
\begin{equation}\label{eq5.4t}
\begin{aligned}
 c_\eps(u_\eps-W_\eps)
 &=c_\eps^2-
 \frac{\alpha_\eps c_\eps^2}{8\pi^2(4-\beta)}+\phi(R)+\frac{\log R}{8\pi^2}
+\frac1{8\pi^2(4-\beta)}\log\frac{\lambda_\eps}{d_\eps c_\eps^2}+o_\eps(1),\\
 c_\eps\partial_r(u_\eps-W_\eps)
 &= R^{-1}r_\eps^{\frac4{\beta-4}}\left(R\phi'(R)+\frac1{8\pi^2}+o_\eps(1)\right), 
\end{aligned}
\end{equation}
 where $\partial_r$ denotes the outward normal derivative on the sphere \(\partial B_{Rr_\eps^{4/(4-\beta)}}(x_\eps)\). While, on the outer boundary \(\partial B_\delta(x_\eps)\), using Lemma \ref{lem4.4} and the fact  \(x_\eps\to0\), we have
\begin{equation}\label{eq5.5t}
\begin{aligned}
 c_\eps(u_\eps-W_\eps)
 &=A_0+O(\delta)+o_\eps(1),\\
 c_\eps\delta\partial_r(u_\eps-W_\eps)
&=O(\delta)+o_\eps(1).
\end{aligned}
\end{equation}

Using the Young's inequality, we have
\begin{equation}\label{eq5.6t}
\begin{aligned}
 &c_\eps^2\int_{A_\eps(R,\delta)}|\Delta u_\eps|^2dx
 -c_\eps^2\int_{A_\eps(R,\delta)}|\Delta W_\eps|^2dx\geq
 2c_\eps^2\int_{A_\eps(R,\delta)}
 \Delta(u_\eps-W_\eps)\Delta W_\eps dx.
\end{aligned}
\end{equation}
Since \(W_\eps\) is biharmonic on the annulus,  integrating by parts twice, we obtain
\begin{equation}\label{equweps}
\begin{aligned}
 c_\eps^2\int_{A_\eps(R,\delta)}
 \Delta(u_\eps-W_\eps)\Delta W_\eps dx
 =c_\eps^2\int_{\partial A_\eps(R,\delta)}
 \left(\Delta W_\eps\partial_\nu(u_\eps-W_\eps)
 -(u_\eps-W_\eps)\partial_\nu\Delta W_\eps\right)d\omega.
\end{aligned}
\end{equation}
Using \(|\partial B_r|=2\pi^2r^3\) and \eqref{eq5.5t}, the contribution of the outer boundary \(\partial B_\delta(x_\eps)\) in \eqref{equweps} is
\[
\begin{aligned}
 &c_\eps^2\int_{\partial B_\delta(x_\eps)}
 \left(\Delta W_\eps\partial_r(u_\eps-W_\eps)-(u_\eps-W_\eps)\partial_r\Delta W_\eps\right)d\omega
=-A_0+O(\delta)+o_\eps(1).
\end{aligned}
\]
Using \eqref{eq5.4t}, the corresponding contribution of the outer boundary $\partial B_{Rr_\eps^{4/(4-\beta)}}(x_\eps)$ is
\[
\begin{aligned}
 &-c_\eps^2\int_{\partial B_{Rr_\eps^{4/(4-\beta)}}(x_\eps)}\left(
 \Delta W_\eps\partial_r(u_\eps-W_\eps)-(u_\eps-W_\eps)\partial_r\Delta W_\eps\right) d\omega\\
 &=c_\eps^2-\frac{\alpha_\eps c_\eps^2}{8\pi^2(4-\beta)}+\left(\phi(R)+\frac R2\phi'(R)\right)+\frac1{8\pi^2}\log R\\
 &\ \ +\frac1{8\pi^2(4-\beta)}\log\frac{\lambda_\eps}{d_\eps c_\eps^2}+\frac1{16\pi^2}+o_\eps(1).
\end{aligned}
\]
 Summing the above two contributions yields
\begin{equation}\label{eq5.8t}
\begin{aligned}
 &c_\eps^2\int_{A_\eps(R,\delta)}\Delta(u_\eps-W_\eps)\Delta W_\eps dx\\
  &=c_\eps^2-\frac{\alpha_\eps c_\eps^2}{8\pi^2(4-\beta)}+
 \left(\phi(R)+\frac R2\phi'(R)+\frac1{16\pi^2}-A_0+\frac1{8\pi^2}\log\delta\right)\\
 &\quad+\frac1{8\pi^2}\log\frac R\delta+\frac1{8\pi^2(4-\beta)}\log\frac{\lambda_\eps}{d_\eps c_\eps^2}+o_\eps(1)+O(\delta|\log\delta|).
\end{aligned}
\end{equation}

We next calculate the left-hand side of \eqref{eq5.6t}.  A change of
variables gives
\[
 c_\eps^2\int_{B_{Rr_\eps^{\frac4{4-\beta}}}(x_\eps)}|\Delta u_\eps|^2dx=\int_{B_R}|\Delta\phi_\eps|^2dx.
\]
By Lemma \ref{lemxe} and Fatou's Lemma, this implies 
\begin{equation} \label{eq5.9t}
 \liminf_{\eps\to0}
 c_\eps^2\int_{B_{Rr_\eps^{\frac4{4-\beta}}}(x_\eps)}|\Delta u_\eps|^2dx
 \geq\int_{B_R}|\Delta\phi|^2dx,
\end{equation}
where $\phi$ is the limit function appearing in Lemma~\ref{lemxe}. 

While,  Lemma~\ref{lem4.4} and integration by parts
give
\begin{equation} \label{eqly5.14t}
\begin{aligned}
 &c_\eps^2\int_{\Omega\setminus B_\delta(x_\eps)}|\Delta u_\eps|^2dx
 =\int_{\Omega\setminus B_\delta(0)}|\Delta G_{a_1}|^2dx+o_\eps(1)\\
 &= -\int_{\partial B_\delta}\Delta G_{a_1}\partial_rG_{a_1}d\omega+\int_{\partial B_\delta}G_{a_1}\partial_r\Delta G_{a_1}d\omega+a_1\int_{\Omega\setminus B_\delta}G_{a_1}^2dx+o_\eps(1)\\
 &=-\frac1{16\pi^2}-\frac1{8\pi^2}\log\delta+A_0+a_1\|G_{a_1}\|_2^2+O(\delta|\log\delta|)+o_\eps(1).
\end{aligned}\end{equation}

On the other hand, a direct computation combined with \eqref{eqreps} yields the exact identity  
\begin{equation}\label{eq5.10t}
\begin{aligned}
 c_\eps^2\int_{A_\eps(R,\delta)}|\Delta W_\eps|^2dx
 &=\frac1{8\pi^2}\log\frac{\delta}{Rr_\eps^{\frac4{4-\beta}}}\\
 &=\frac{\alpha_\eps c_\eps^2}{8\pi^2(4-\beta)}
 -\frac1{8\pi^2(4-\beta)}\log\frac{\lambda_\eps}{d_\eps c_\eps^2}+\frac1{8\pi^2}\log\frac\delta R,
\end{aligned}
\end{equation}
where $\alpha_\eps$ admits the following asymptotic expansion:
\begin{equation}\label{eqly5.15t}
 \alpha_\eps=\left(32\pi^2\left(1-\frac\beta4\right)-\eps\right)
 \left(1+\frac{a_1\|G_{a_1}\|_2^2}{c_\eps^2}+o(c_\eps^{-2})\right),
\end{equation}
which can be easily verified by the fact that 
 \(c_\eps u_\eps\to G_{a_1}\) in \(L^2(\Omega)\).

Combining \eqref{eq5.6t}, \eqref{eq5.8t}--\eqref{eqly5.15t},  we obtain
\begin{equation}\label{eqly5.16t}
\begin{aligned}
 \frac1{8\pi^2(4-\beta)}\log\frac{\lambda_\eps}{d_\eps c_\eps^2}
  &\leq-\int_{B_R}|\Delta\phi|^2dx-2\phi(R)-R\phi'(R)-\frac1{16\pi^2}+A_0\\
 &\quad-\frac1{8\pi^2}\log R+o_\eps(1)+O(\delta|\log\delta|).\\
\end{aligned}
\end{equation}
 Letting first \(\eps\to0\) and then \(\delta\to0\) in \eqref{eqly5.16t},   for every fixed \(R>1\), we have
\begin{equation}\label{eq5.22tb}
\begin{aligned}
 \limsup_{\eps\to0}
 \log\frac{\lambda_\eps}{d_\eps c_\eps^2}
 \leq{}&32\pi^2\left(1-\frac\beta4\right)
 \Bigg(
 A_0-\int_{B_R}|\Delta\phi|^2dx
 -2\phi(R)-R\phi'(R)\\
 &\ \
 -\frac1{8\pi^2}\log R-\frac1{16\pi^2}
 \Bigg).
\end{aligned}
\end{equation}

It remains to evaluate the limit of the following expression
in \eqref{eq5.22tb}
\begin{equation*} 
\int_{B_R}|\Delta\phi|^2dx+2\phi(R)+R\phi'(R)+\frac{1}{8\pi^2}\log R, \ \text{as} \ R\to +\infty. 
\end{equation*}

Using  asymptotic results  proved in \cite{AWZ} (see Theorem~1.3, Lemma~2.7 and equations
\textup{(2.37)--(2.39)} of \cite{AWZ}),  there exists some constant $C_0$ such that for some \(k>0\),
\begin{equation}\label{eq5.22tc}
\begin{aligned}
 \phi(R)&=-\frac1{8\pi^2}\log R+C_0+O(R^{-k}),\\
 R\phi'(R)&=-\frac1{8\pi^2}+O(R^{-k}),\\
 \Delta\phi(R)&=-\frac1{4\pi^2R^2}+O(R^{-2-k}).
\end{aligned}
\end{equation}

We claim that  $C_0$ has the following integral representation: 
\begin{equation}\label{eq:C0-integral}
C_0=\frac{1}{8\pi^2}\int_{\RR^4}\frac{\log|y|
e^{64\pi^2(1-\frac{\beta}{4})\phi(y)}}{|y|^\beta}dy.
\end{equation}
Indeed,  applying Lemma~2.4 of \cite{AWZ} and using $\phi(0)=0$, we obtain
\begin{equation}\label{eqphiintegral}
\phi(x)=\frac{1}{8\pi^2}\int_{\RR^4}\log\frac{|y|}{|x-y|}\frac{e^{64\pi^2(1-\frac{\beta}{4})\phi(y)}}{|y|^\beta}dy.
\end{equation}
Since $0<\beta<4$, $\phi\in C_{\mathrm{loc}}^0(\RR^4)$ and
$\phi\leq0$, we have
\[
\int_{B_1}\frac{(1+|\log|y||)e^{64\pi^2(1-\frac{\beta}{4})\phi(y)}}{|y|^\beta}dy
\leq 2\pi^2\int_0^1(1+|\log r|)r^{3-\beta}dr<+\infty.
\]
 On the other hand, the asymptotic expansion
\eqref{eq5.22tc} gives
\[
e^{64\pi^2(1-\frac{\beta}{4})\phi(y)}=O\bigl(|y|^{-8+2\beta}\bigr)\text{ as }|y|\to\infty.
\]
Consequently, again by  $0<\beta<4$, we  obtain
\[
\begin{aligned}
&\int_{\RR^4\setminus B_1}
\frac{(1+|\log|y||)
e^{64\pi^2(1-\frac{\beta}{4})\phi(y)}}{|y|^\beta}dy\leq
C\int_1^\infty(1+\log r)r^{\beta-5}dr
<+\infty.
\end{aligned}
\]
 Therefore,
\begin{equation}\label{eq:log-moment}
\int_{\RR^4}
\frac{(1+|\log|y||)
e^{64\pi^2(1-\frac{\beta}{4})\phi(y)}}{|y|^\beta}dy
<+\infty.
\end{equation}
Using the mass identity in \eqref{eq:phiMassShift} and \eqref{eq:log-moment},
we may rewrite \eqref{eqphiintegral} as

\begin{align}\label{addlimforphi1020}
\phi(x)+\frac{1}{8\pi^2}\log|x|
=\frac{1}{8\pi^2}
\left(\int_{\RR^4}\frac{\log|y|
e^{64\pi^2(1-\frac{\beta}{4})\phi(y)}
}{|y|^\beta}dy+\int_{\RR^4}
\log\frac{|x|}{|x-y|}
\frac{e^{64\pi^2(1-\frac{\beta}{4})\phi(y)}}{|y|^\beta}dy\right).
\end{align}

Now, we show  \begin{equation}\label{addlim1020}
    \lim_{|x|\to\infty}
\int_{\RR^4}\log\frac{|x|}{|x-y|}\frac{e^{64\pi^2(1-\frac{\beta}{4})\phi(y)}}{|y|^\beta}dy
=0.
\end{equation}
By the following key observation that
\begin{equation}\label{addasym1020}
    f(y):=\frac{e^{64\pi^2(1-\beta/4)\phi(y)}}{|y|^\beta}
=O(|y|^{\beta-8})
\ \ \text{ as }\ \ |y|\to\infty,
\end{equation}
which can be easily verified from \eqref{eq5.22tc}. 

Write $r=|x|$ and decompose
\[
\int_{\mathbb R^4}\log\frac{|x|}{|x-y|}f(y)dy
=\left(\int_{|y|<r/2}+\int_{|y|\ge r/2}\right)\log\frac{|x|}{|x-y|}f(y)dy
=:I_1(x)+I_2(x).
\]
On $\{|y|<r/2\}$, we have
\[
\left|\log\frac{|x|}{|x-y|}\right|\le\log 2.
\]
Moreover, for each fixed $y$,
\[
\log\frac{|x|}{|x-y|}f(y)\rightarrow0 \text{ as }|x|\to\infty.
\]
Since $\beta<4$ and $\phi$ is locally bounded, it follows that $f\in L^1(\mathbb R^4)$. Using the dominated convergence theorem, we therefore obtain
\[
I_1(x)\rightarrow0, \text{ as } |x|\to\infty.
\]

To estimate $I_2$, we write $x=re$, where $|e|=1$, and use the
change of variables $y=rz$.  Then, for all sufficiently large $r$, it follows
from \eqref{addasym1020} and the assumption $\beta<4$ that
\[
\begin{aligned}
|I_2(x)|
&\le C\int_{|y|\ge r/2}
 \left|\log\frac{|x|}{|x-y|}\right||y|^{\beta-8}dy\\
&=Cr^{\beta-4} \int_{|z|\ge1/2} |\log|e-z|||z|^{\beta-8}dz\\
 &\leq Cr^{\beta-4}\rightarrow0, 
\end{aligned} 
\]
as $|x|\rightarrow \infty$.
Combining the two estimates, we obtain \eqref{addlim1020}.

It follows from \eqref{addlimforphi1020} and \eqref{addlim1020} that
\[
\lim_{|x|\to\infty}
\left( \phi(x)+\frac{1}{8\pi^2}\log|x| \right)
=\frac{1}{8\pi^2}\int_{\RR^4} \frac{\log|y| e^{64\pi^2(1-\frac{\beta}{4})\phi(y)}}{|y|^\beta}dy.
\]
Comparing this  with \eqref{eq5.22tc}, we conclude the integral representation \eqref{eq:C0-integral}.


On the other hand, using the asymptotic behavior of $\Delta\phi(R)$ in  \eqref{eq5.22tc}, we see that
\begin{equation*} 
\begin{aligned}
&\int_{B_R}|\Delta\phi|^2dx-\frac1{8\pi^2}\log R= \int_{B_1}|\Delta\phi|^2dx +2\pi^2\int_1^R
\left( r^3|\Delta\phi(r)|^2-\frac1{16\pi^4r} \right)dr
\end{aligned}
\end{equation*}
has a finite limit as $R\to\infty$.  Combining this with \eqref{eq5.22tc}, we get 
\begin{equation}\begin{split}\label{eqinfinitr2}
&\lim_{R\rightarrow\infty}\Big(\int_{B_R}|\Delta\phi|^2dx+2\phi(R)+R\phi'(R)+\frac{1}{8\pi^2}\log R\Big)\\
&=\lim_{R\rightarrow\infty}\Big(\int_{B_R}|\Delta\phi|^2dx-\frac{1}{8\pi^2}\log R\Big)+2C_0-\frac{1}{8\pi^2}.
\end{split}\end{equation}

Next, we will consider the limit:
\begin{equation*} 
\lim_{R\rightarrow\infty}\left(\int_{B_R}|\Delta\phi|^2dx-\frac1{8\pi^2}\log R\right).
\end{equation*}
Integrating by parts on $B_R$, we have
\begin{align*}
 \int_{B_R}\phi\Delta^2\phi dx
 &=\int_{\partial B_R}
   \phi\frac{\partial\Delta\phi}{\partial r}d\omega
   -\int_{B_R}\nabla\phi\cdot\nabla\Delta\phi dx
 \notag\\
 &=\int_{\partial B_R}
   \phi\frac{\partial\Delta\phi}{\partial r}d\omega
   -\int_{\partial B_R}
   \Delta\phi\frac{\partial\phi}{\partial r}d\omega
   +\int_{B_R}|\Delta\phi|^2dx.
\end{align*}
Hence,
\begin{equation*} 
\int_{B_R}|\Delta\phi|^2dx
=\int_{B_R} \frac{\phi(x)e^{64\pi^2(1-\frac \beta4)\phi(x)}}{|x|^{\beta}}dx
-\phi(R)+\frac1{16\pi^2}+o(1).
\end{equation*}
Consequently, using \eqref{eq5.22tc} and the assumption $\beta<4$, we derive
\begin{equation}\label{eqrefinfinitr4}
\begin{aligned}
&\lim_{R\to\infty}
\left(\int_{B_R}|\Delta\phi|^2dx -\frac1{8\pi^2}\log R\right)\\
&=\int_{\RR^4}\frac{\phi(x)e^{64\pi^2(1-\frac \beta4)\phi(x)}}{|x|^{\beta}}dx-\lim_{R\to\infty}
\left(\phi(R)+\frac1{8\pi^2}\log R\right)+\frac1{16\pi^2}\\
&=\int_{\RR^4}\frac{\phi(x)e^{64\pi^2(1-\frac \beta4)\phi(x)}}{|x|^{\beta}}dx-C_0+\frac1{16\pi^2}\\
&:=J_0-C_0+\frac1{16\pi^2}
\end{aligned}
\end{equation}

Combining  \eqref{eq5.22tb},  \eqref{eqinfinitr2} and  \eqref{eqrefinfinitr4}, we obtain
\begin{equation*} 
{
\begin{aligned}
\limsup_{\eps\to0}\frac{\lambda_\eps}{c_\eps^2}
&\leq\exp\Big\{32\pi^2\big(1-\frac \beta4\big)\left(A_0-2C_0+\frac1{16\pi^2}-
\lim_{R\to\infty}\Big(\int_{B_R}|\Delta\phi|^2dx-\frac{\log R}{8\pi^2}\Big)\right)\Big\}\\
&\leq\exp\Big(32\pi^2\big(1-\frac \beta4\big)\left(A_0-C_0-J_0\right)\Big).
\end{aligned}}
\end{equation*}
The preceding estimate and Lemma \ref{lemexp} together yield
\begin{equation*} 
\begin{aligned}
S_{\Omega,q,\beta}\leq\int_\Omega\frac{1}{|x|^{\beta}}dx+\exp\Big(32\pi^2\big(1-\frac \beta4\big)\left(A_0-C_0-J_0\right)\Big).
\end{aligned}
\end{equation*}
Therefore, the proof is finished.
\end{proof}

Finally, we can complete the proof of Theorem~\ref{thm1} as follows. 
\begin{proof}[Proof of Theorem~\ref{thm1}]

Theorem~\ref{thm1} is an immediate consequence of combining  the bounded case of $\{c_\eps\}$ analyzed at the end of Section \ref{boundedcase1020}, with the blow-up case derived in Proposition \ref{eqly5.23t}.

\end{proof}


\section{Existence of extremal functions}

In this section, we establish the existence of extremal functions attaining the supremum $S_{\Omega,q,\beta}$. To this end, we construct suitable functions $ \varphi_\epsilon \in H_0^2(\Omega) $ such that $ \|\Delta \varphi_\epsilon\|_2^2 = 1 $ and
\[
\int_{\Omega} \frac{\exp\Big({32\pi^2\big(1-\frac{\beta}{4}\big) \varphi_\epsilon^2 q(\|\varphi_\epsilon\|_2^2)}\Big)}{|x|^\beta} dx > \int_\Omega\frac{1}{|x|^{\beta}}dx+\exp\Big(32\pi^2\big(1-\frac \beta4\big)\left(A_0-C_0-J_0\right)\Big)
\]
provided that $ a_1 > 0 $ is sufficiently small. Since $ a_i \geq 0 $ $( i = 1, 2, \ldots, k )$, we only need to show
\begin{equation}\label{eq7.6t}
\int_{\Omega} \frac{\exp\Big({32\pi^2\big(1-\frac{\beta}{4}\big) \varphi_\epsilon^2 (1+a_1 \|\varphi_\epsilon\|_2^2)} \Big)}{|x|^\beta} dx > \int_\Omega\frac{1}{|x|^{\beta}}dx+\exp\Big(32\pi^2\big(1-\frac \beta4\big)\left(A_0-C_0-J_0\right)\Big),
\end{equation}
provided that $ a_1 > 0 $ small enough.  Standard regularity theory and \eqref{eqly5.23t} then give Theorem \ref{thm2}.

For any $\eps>0$, define
\begin{equation}\label{eq7.7t}
R=\left(\log\frac1\eps\right)^m,
\end{equation}
with  $m\cdot\min\{k,4-\beta\}>2$.  Here $k$ is the constant appearing in \eqref{eq5.22tc}. 

With this choice of $R$, it is easy to verify that $R\to\infty$, $R\eps\to0$, as $\eps\to 0$. Writing $r=|x|$, define
\begin{equation}\label{eq7.8t}
\varphi_\eps(x)=
\begin{cases}
\displaystyle
c+\frac{a+\phi\!\left(\frac{x}{\eps}\right)+A_0+H(x)+br^2}{c},
&0\leq r\leq R\eps,\\[3mm]
\displaystyle
\frac{G_{a_1}(x)}{c},
&r>R\eps,
\end{cases}
\end{equation}
where $a,b,c$ will be determined later.

To ensure that $\varphi_\eps(x)\in H^2_0(\Omega)$, we require that 
\[
c^2+a+\phi(R)+bR^2\eps^2=-\frac1{8\pi^2}\log(R\eps)
\]
and 
\[
\frac{\phi'(R)}{\eps}+2bR\eps=-\frac1{8\pi^2R\eps}.
\]
Consequently, we have
\begin{equation}\label{eq7.9t}
{
b=-\frac1{16\pi^2 R^2\eps^2}-\frac{\phi'(R)}{2R\eps^2},\ \ 
a=-c^2-\phi(R)-\frac1{8\pi^2}\log(R\eps)+\frac{1}{16\pi^2}+ \frac{R\phi'(R)}2.}
\end{equation}

From \eqref{eq7.8t}, on $B_{R\eps}$, we have
\[
\Delta\varphi_\eps(x)
=\frac1c\left(\eps^{-2}\Delta\phi\left(\frac{x}{\eps}\right) +\Delta H(x)+8b \right).
\]

After the change of variables $x=\eps y$, we have
\begin{equation}
\begin{split}\label{eq7.11ta}
c^2\int_{B_{R\eps}}|\Delta\varphi_\eps|^2dx
&=\int_{B_R}|\Delta\phi(y)|^2dy+2\eps^2\int_{B_R}
\Delta\phi(y)\Delta H(\eps y)dy\\
&\quad+16b\eps^2\int_{B_R}\Delta\phi(y)dy+\eps^4\int_{B_R}|\Delta H(\eps y)|^2dy\\
&\quad+16b\eps^4\int_{B_R}\Delta H(\eps y)dy +64b^2\eps^4|B_R|.
\end{split}
\end{equation}

Combining \eqref{eqlimitDe},\eqref{eq7.9t} and \eqref{eq5.22tc}, we obtain
\[
\int_{B_R}|\Delta\phi|dy=O(R^2)
\]
and
\[
8b\eps^2=-\frac4{R^2} \left(\frac1{8\pi^2}+R\phi'(R)\right) =O(R^{-2-k}).
\]
Since $\Delta H$ is bounded near the
origin, these estimates imply that the sum of the five error terms involving $\eps$ in \eqref{eq7.11ta} is 
\[
O(\eps^2R^2)+ O(R^{-k}).
\]
Therefore, we have
\begin{equation}\label{eq7.11t}
{
\int_{B_{R\eps}}|\Delta\varphi_\eps|^2dx
=\frac1{c^2}\left( \int_{B_R}|\Delta\phi|^2dy+O(\eps^2R^2)+ O(R^{-k})\right).}
\end{equation}

On the other hand, since
$\Delta^2G_{a_1}=a_1G_{a_1}$ away from the origin, integration by parts gives
\begin{align}\label{addgelin1021}
\int_{\Omega\setminus B_r}|\Delta G_{a_1}|^2dx
&=a_1\int_{\Omega\setminus B_r}G_{a_1}^2dx
-\int_{\partial B_r}\Delta G_{a_1} ~\partial_rG_{a_1}d\omega \nonumber\\
&\quad+\int_{\partial B_r}G_{a_1} ~\partial_r\Delta G_{a_1}d\omega.
\end{align}
By Lemma \ref{lem4.4}, we have
\[
\begin{aligned}
G_{a_1}&=-\frac1{8\pi^2}\log r+A_0+O(r),
&\partial_rG_{a_1}&=-\frac1{8\pi^2r}+O(1),\\
\Delta G_{a_1}&=-\frac1{4\pi^2r^2}+O(1),
&\partial_r\Delta G_{a_1}&=\frac1{2\pi^2r^3}+O(1),
\end{aligned}
\]
where $\partial_r$ denotes the outward normal derivative on the sphere \(\partial B_{r}\). Combining these asymptotic behavior of the Green's function above and the fact that $|\partial B_r|=2\pi^2r^3$, we obtain
\[
\begin{aligned}
-\int_{\partial B_r}\Delta G_{a_1} ~\partial_rG_{a_1}d\omega
&=-\frac1{16\pi^2}+O(r),\\
\int_{\partial B_r}G_{a_1} ~\partial_r\Delta G_{a_1}d\omega
&=-\frac1{8\pi^2}\log r+A_0+O(r|\log r|),\\
a_1\int_{\Omega\setminus B_r}G_{a_1}^2dx
&=a_1\|G_{a_1}\|_2^2+O(r^4|\log r|^2).
\end{aligned}
\]
Substituting $r=R\eps$ in \eqref{addgelin1021} and using the estimates above, we get
\begin{equation}\label{eq7.12t}
\begin{aligned}
\int_{\Omega\setminus B_{R\eps}}|\Delta\varphi_\eps|^2dx
=\frac1{c^2}\Bigg(
&-\frac1{16\pi^2}
-\frac1{8\pi^2}\log(R\eps)
+A_0+a_1\|G_{a_1}\|_2^2+O(R\eps \log(R\eps))
\Bigg).
\end{aligned}
\end{equation}

Combining \eqref{eq7.11t} and \eqref{eq7.12t}, we  have
\begin{equation*} 
{
\begin{aligned}
\int_\Omega|\Delta\varphi_\eps|^2dx=\frac1{c^2}\Big(
&\int_{B_R}|\Delta\phi|^2dx-\frac1{16\pi^2}-\frac1{8\pi^2}\log(R\eps)+A_0\\
&+a_1\|G_{a_1}\|_2^2+O(\eps^2R^2)+ O(R^{-k})+O(R\eps \log(R\eps)) \Big).
\end{aligned}}
\end{equation*}
We now choose the  $c$ so that $\|\Delta\varphi_\eps\|_2=1$. Then it holds 
\begin{equation}\label{eq7.13t}
\begin{aligned}
c^2 &=\int_{B_R}|\Delta\phi|^2dx-\frac1{16\pi^2}-\frac1{8\pi^2}\log(R\eps)+A_0+a_1\|G_{a_1}\|_2^2\\
&\ \ +O(\eps^2R^2)+ O(R^{-k})+O(R\eps \log(R\eps)).
\end{aligned}
\end{equation}

On the other hand, using \eqref{eq5.22tc}, for  $|x|$ large enough, we have
\begin{equation*}
e^{64\pi^2(1-\frac{\beta}{4})\phi(x)} \leq C|x|^{-8+2\beta}, \text{ and } |\phi(x)| \leq C(1+\log|x|).
\end{equation*}
Then using the L'Hospital's Rule, we obtain
\[
\begin{aligned}
\int_{\mathbb{R}^4\setminus B_R(0)} \frac{|\phi(x)|e^{64\pi^2(1-\frac{\beta}{4})\phi(x)}}{|x|^{\beta}} dx &\leq C\int_R^\infty (1+\log r)r^{-\beta}r^{-8+2\beta}r^3 dr\\& = C\int_R^\infty (1+\log r)r^{\beta-5} dr\\&= O( R^{\beta-4} \log R ).
\end{aligned}
  \]
From the preceding estimate and \eqref{eqrefinfinitr4}, \eqref{eq5.22tc}, we deduce 
\begin{equation}\label{addconver1022}
   \int_{B_R}|\Delta\phi|^2dx
=\frac{1}{8\pi^2}\log R+J_0-C_0+\frac{1}{16\pi^2}+O(R^{-k})+O(R^{\beta-4}\log R), 
\end{equation}
where  $J_0$ and $C_0$ are the constants appearing in \eqref{eqly5.23t}.

Combining \eqref{eq7.13t} and \eqref{addconver1022}, we have
\[
\begin{aligned}
c^2
&=\frac{1}{8\pi^2}\log\frac1\eps+J_0-C_0+A_0+a_1\|G_{a_1}\|_2^2+O(R^{-k})+O(R^{\beta-4}\log R).
\end{aligned}
\]
With $R$ chosen as in \eqref{eq7.7t}, the sum of the above error terms is indeed $o(c^{-2})$. Thus, we obtain
\begin{equation}\label{eq7.15t}
\begin{aligned}
c^2
&=\frac{1}{8\pi^2}\log\frac1\eps+J_0-C_0+A_0+a_1\|G_{a_1}\|_2^2+o(c^{-2}).
\end{aligned}
\end{equation}

Combining \eqref{eq7.9t} and  \eqref{eq7.15t},  we also obtain
\begin{equation}\label{eq7.16t}
{
a+A_0
=-J_0
-a_1\|G_{a_1}\|_2^2+o(c^{-2}).}
\end{equation}

Since
 $G_{a_1}(x)=O(1+|\log|x||)$ near the origin.  Therefore, we have
\[
\int_{B_{R\eps}}\varphi_\eps^2dx
=O\left(c^2(R\eps)^4\right)
+O\left(\eps^4\log R\right)
=o(c^{-4}),
\]
and 
\[
\int_{\Omega\setminus B_{R\eps}}\varphi_\eps^2dx
=\frac1{c^2}\Big(
\|G_{a_1}\|_2^2
+O\!\left((R\eps)^4|\log(R\eps)|^2\right)
\Big).
\]
Consequently,
\begin{equation}\label{eq7.17t}
{
\|\varphi_\eps\|_2^2
=\frac{\|G_{a_1}\|_2^2}{c^2}+o(c^{-4}).}
\end{equation}

For $x=\eps y\in B_{R\eps}$, \eqref{eq7.8t} and \eqref{eq7.16t} give, uniformly on
$B_R$,
\begin{equation}\label{eq7.19t}
\begin{aligned}
\varphi_\eps(\eps y)
=c+\frac1c\Bigg(
\phi(y)
&-J_0-a_1\|G_{a_1}\|_2^2+o(c^{-2})
\Bigg).
\end{aligned}
\end{equation}
Here, the term $H(\eps y)$ is absorbed into the error because
$|H(\eps y)|\leq C R\eps=o(c^{-2})$.

Combining  \eqref{eq7.17t} and \eqref{eq7.19t}, we obtain
\begin{equation*}\begin{split} 
\varphi_\eps^2(\eps y)\left(1+a_1\|\varphi_\eps\|_2^2\right)
&=c^2+2\Bigg(\phi(y)-J_0-a_1\|G_{a_1}\|_2^2\Bigg)\\
&\ \ +a_1\|G_{a_1}\|_2^2+\frac1{c^2}\Bigg\{
\Bigg(\phi(y)-J_0\Bigg)^2-a_1^2\|G_{a_1}\|_2^4\Bigg\}+o(c^{-2}).
\end{split}\end{equation*}
By \eqref{eq7.15t}, this implies
\begin{align*}\begin{split} 
&32\pi^2\left(1-\frac \beta4\right)
\varphi_\eps^2(\eps y)\left(1+a_1\|\varphi_\eps\|_2^2\right)\\
&=(4-\beta)\log\frac1\eps+32\pi^2\left(1-\frac \beta4\right)
\Bigg(A_0-C_0-J_0 +2\phi(y)\Bigg)\\
&\ \ +\frac{32\pi^2(1-\frac \beta4)}{c^2}\Bigg(\Big(\phi(y)-J_0\Big)^2
-a_1^2\|G_{a_1}\|_2^4\Bigg)+o(c^{-2}).
\end{split}\end{align*}

Therefore, we have
\begin{equation}\label{eq7.22t}
\begin{aligned}
&\int_{B_{R\eps}}
\frac{e^{32\pi^2(1-\frac \beta4)
\varphi_\eps^2(1+a_1\|\varphi_\eps\|_2^2)}}{|x|^{\beta}}dx\\
& \geq\exp\Bigg(32\pi^2\left(1-\frac \beta4\right) \big(A_0-C_0-J_0 \big)\Bigg)\\
&\ \ \Bigg(1+\frac{32\pi^2(1-\frac \beta4)}{c^2} \Bigg(\int_{\RR^4}\Big(\phi(y)-J_0
\Big)^2\frac{e^{64\pi^2(1-\frac \beta4)\phi(y)}}{|y|^\beta}dy
-a_1^2\|G_{a_1}\|_2^4\Bigg)+o(c^{-2})\Bigg),
\end{aligned}
\end{equation}
where we have used the following observation that for each $j=0,1,2,4$,
\[
\int_{\RR^4\setminus B_R}
\frac{|\phi(y)|^j e^{64\pi^2(1-\frac \beta4)\phi(y)}}{|y|^\beta}dy
=O\!\left(R^{-(4-\beta)}\log^jR\right)=o(c^{-2}),
\]
which can be derived from \eqref{eq5.22tc}.

 Since $e^s\geq1+s$, then on $\Omega\setminus B_{R\eps}$  we have
\begin{equation}\label{eq7.23t}
\begin{split}
&\int_{\Omega\setminus B_{R\eps}}\frac{e^{32\pi^2(1-\frac \beta4)
\varphi_\eps^2(1+a_1\|\varphi_\eps\|_2^2)}}{|x|^{\beta}}dx\\
&\geq\int_{\Omega\setminus B_{R\eps}}|x|^{-\beta}dx
+\frac{32\pi^2(1-\frac \beta4)}{c^2}\int_{\Omega\setminus B_{R\eps}}
\frac{G_{a_1}^2(x)}{|x|^{\beta}}dx\\
&=\int_\Omega|x|^{-\beta}dx+\frac{32\pi^2(1-\frac \beta4)}{c^2}
\int_\Omega\frac{G_{a_1}^2(x)}{|x|^{\beta}}dx+o(c^{-2}).
\end{split}\end{equation}


Finally, combining \eqref{eq7.22t} and \eqref{eq7.23t}, we obtain 
\begin{equation*} 
\begin{aligned}
&\int_\Omega\frac{e^{32\pi^2(1-\frac\beta4)\varphi_\eps^2(1+a_1\|\varphi_\eps\|_2^2)}}{|x|^{\beta}}dx\\
\geq &\int_\Omega\frac{1}{|x|^{\beta}}dx+\exp\Big(32\pi^2\left(1-\frac \beta4\right)\big(A_0-C_0-J_0\big)\Big)\\
&\ \ +\frac{32\pi^2(1-\frac \beta4)}{c^2}
\Bigg\{\int_\Omega\frac{G_{a_1}^2(x)}{|x|^{\beta}}dx+\exp\Big(32\pi^2\left(1-\frac \beta4\right)
\big(A_0-C_0-J_0\big)\Big)\\
&\ \ \times\Bigg(
\int_{\RR^4}\Bigg(\phi(y)-J_0\Bigg)^2
\frac{e^{64\pi^2(1-\frac \beta4)\phi(y)}}{|y|^\beta}dy
-a_1^2\|G_{a_1}\|_2^4\Bigg)\Bigg\}+o(c^{-2}).
\end{aligned}
\end{equation*}
When $a_1=0$, since the expression inside  braces is strictly positive, hence for sufficiently small $\eps$, we desire \eqref{eq7.6t}. While, when $a_1>0$, by the Proposition~\ref{prop:weighted-green-uniform} below, there eixsts some $a_*>0$ such that for any $a_1\leq a_*$,
\[
\int_\Omega\frac{G_{a_1}^2}{|x|^{\beta}}dx
\geq\frac12\int_\Omega\frac{G_0^2}{|x|^{\beta}}dx>0
\]
and $\|G_{a_1}\|_2^4\leq C $, with $C$  independent of
$a_1\in[0,a_*]$. 
  Thus if  $a_*$ is selected to be so small that
\[
Ca_*^2 <\frac14\int_\Omega\frac{G_0^2}{|x|^{\beta}}dx.
\]
Then, for sufficiently small $\eps$, we have
\[
\int_{\Omega} \frac{\exp\Big({32\pi^2\big(1-\frac{\beta}{4}\big) \varphi_\epsilon^2 q(\|\varphi_\epsilon\|_2^2)}\Big)}{|x|^\beta} dx > \int_\Omega\frac{1}{|x|^{\beta}}dx+\exp\Big(32\pi^2\big(1-\frac \beta4\big)\left(A_0-C_0-J_0\right)\Big). 
\]
Therefore, \eqref{eq7.6t} is derived.

 \begin{proposition}\label{prop:weighted-green-uniform}
For each $0<\alpha_0<\lambda(\Omega)$, there exists a constant $C$ depending only on
$\Omega$, $\lambda(\Omega)$ and $\alpha_0$ such that for any $0\leq a_1\leq \alpha_0$,
\[
\|G_{a_1}\|_2^2\leq C.
\]
Moreover,
\[
\int_\Omega\frac{G_{a_1}^2}{|x|^{\beta}}dx
\rightarrow \int_\Omega\frac{G_0^2}{|x|^{\beta}}dx>0,
\] as $a_1\rightarrow0.$
\end{proposition}

\begin{proof}
Taking $r=\frac{3}{2}$ in Lemmas~\ref{lem:weighted43} and~\ref{lem4.4}, weak lower semicontinuity gives
\[
\|G_{a_1}\|_{W_0^{2,3/2}(\Omega)}
\leq\liminf_{\eps\to0}\|c_\eps u_\eps\|_{W_0^{2,3/2}(\Omega)}\leq C
\]
uniformly for $0\leq a_1\leq\alpha_0$.
The Sobolev embedding $W_0^{2,3/2}(\Omega)\hookrightarrow L^2(\Omega)$ implies for any $0\leq a_1\leq\alpha_0,$
\begin{equation}\label{l2eat1028}
   \|G_{a_1}\|_2^2\leq C. 
\end{equation}

Now, let $G_0(x,y)$ denote the  Green kernel of $\Delta^2$ at $y$
on $\Omega$ with Dirichlet boundary conditions. That is, for each fixed $y\in\Omega$,  it satisfies
\[
\begin{cases}
\Delta_x^2G_0(x,y)=\delta_y & \text{in }\Omega,\\
G_0(x,y)=\partial_{\nu_x}G_0(x,y)=0
& \text{on }\partial\Omega.
\end{cases}
\]
 By Lemma~\ref{lem4.4}, we have
\begin{equation*} 
G_{a_1}(x,0)=-\frac{1}{8\pi^2}\log|x|+H_{a_1}(x),
\end{equation*}
where $H_{a_1}\in C^3(\overline{\Omega})$. Similarly, $G_0$ admits the expansion
\[G_0(x,0)=-\frac{1}{8\pi^2}\log|x|+H_{0}(x),\]
where $H_{0}\in C^3(\overline{\Omega})$ denotes the regular part of $G_0$.
Thus, for any $x\in\Omega\setminus\{0\}$,
\begin{equation*} 
G_{a_1}(x,0)-G_0(x,0)=H_{a_1}(x)-H_{0}(x).
\end{equation*}
We write $G_0(x)=G_0(x,0)$ and $G_{a_1}(x)=G_{a_1}(x,0)$, and for any $x\in\Omega\setminus\{0\}$ define
\[
\omega_{a_1}(x):=G_{a_1}(x)-G_0(x).
\]
By the expansions of the two Green functions, we have 
$$\omega_{a_1}(x)=H_{a_1}(x)-H_{0}(x), x\in\Omega\setminus\{0\}.$$
Moreover, \(\omega_{a_1}\) can be extended continuously to the origin by 
\[
\omega_{a_1}(0):=\lim_{x\to0}\omega_{a_1}(x)=H_{a_1}(0)-H_{0}(0).
\]
Hence $\omega_{a_1}\in C^3(\overline{\Omega})$.
Besides, $\omega_{a_1}$ is a distributional solution of the equation 
\[
\begin{cases}
\Delta^2\omega_{a_1}=a_1G_{a_1}
& \text{in }\Omega,\\
\omega_{a_1}=\partial_\nu\omega_{a_1}=0 & \text{on }\partial\Omega.
\end{cases}
\]
Applying the Green representation, we have
\begin{equation}\label{gerrrep1028}
    \omega_{a_1}(x)=a_1\int_\Omega G_0(x,y)G_{a_1}(y)\,dy.
\end{equation}
Using the logarithmic estimate for the Green kernel $G_0(x,y)$ at $y$, we have
\[
\sup_{x\in\Omega}\|G_0(x,\cdot)\|_2^2\leq C\int_0^{\operatorname{diam}(\Omega)}r^3\bigl(1+|\log r|^2\bigr)dr<\infty.
\]
Using the Cauchy--Schwarz inequality, \eqref{gerrrep1028} and \eqref{l2eat1028}, we obtain
\[
\|\omega_{a_1}\|_{L^\infty(\Omega)}
\leq a_1\sup_{x\in\Omega}\|G_0(x,\cdot)\|_2\|G_{a_1}\|_2\leq Ca_1.
\]
Therefore, $\|G_{a_1}-G_0\|_{L^\infty(\Omega)}\leq Ca_1$, which gives $G_{a_1}(x)\to G_0(x)$ for $x\neq0$ as $a_1\to0$.

Moreover, since
\[\begin{aligned}
G^2_{a_1}(x)&=\bigl(G_0(x)+\omega_{a_1}(x)\bigr)^2\leq 2G^2_0(x) +2\omega^2_{a_1}(x)
\leq 2G^2_0(x)+2C^2\alpha_0^2.
\end{aligned}\]
Hence 
\[0\leq\frac{G_{a_1}(x)^2}{|x|^\beta}\leq\frac{2G_0(x)^2+2C^2\alpha_0^2}{|x|^\beta},\]
where the  right hand  is independent of $a_1$ and integrable over $\Omega$. The dominated convergence theorem therefore yields
\[
\int_\Omega\frac{G_{a_1}^2}{|x|^\beta}dx\rightarrow\int_\Omega\frac{G_0^2}{|x|^\beta}dx>0
\ \ \ \text{as}\  \ a_1\to0.
\]
This completes the proof.
\end{proof}

\section{Appendix-Pointwise gradient estimates and proof of Proposition \ref{lemDM1}}
In this section, we present the proof of Proposition~\ref{lemDM1}, which
relies on the following gradient estimates.
\begin{lemma}\label{lemDMsquare}Let $1\le j\le3$.  For any $a\in\RR^4$ and $r>0$, 
there exists $C>0$ such that
\begin{equation*} 
\int_{B_r(a)\cap\Omega}|\nabla^j(u_\eps^2)|dx\le Cr^{4-j}.
\end{equation*}
\end{lemma}

\begin{proof}
We first  prove that $f_\eps$ is bounded in $ L(\log L)^{\frac{1}{2}}(\Omega)$. In the proof of Lemma \ref{lemge},  we have showed that $f_\eps$  in \eqref{definforf} is bounded in $ L(\log L)^{\frac{1}{2}}(B_{4R_0}(x_\eps))$. Thus, it remains to show that $f_\eps$ is bounded in $ L(\log L)^{\frac{1}{2}}(\Omega\setminus B_{4R_0}(x_\eps))$.

Fix $R_0>0$, since $x_\eps\to0$,  we have $|x_\eps|\le R_0$ for all sufficiently small $\eps$.
Consequently,  for any $x\in\Omega\setminus B_{4R_0}(x_\eps)$, we have
\[
|x|\ge |x-x_\eps|-|x_\eps|\ge 3R_0.
\]
Thus the weight $|x|^{-\beta}$ is uniformly bounded on $\Omega\setminus B_{4R_0}(x_\eps)$. Hence, on $\Omega\setminus B_{4R_0}(x_\eps)$, we have
\[
|f_\eps|\le C|u_\eps|e^{\alpha_\eps u_\eps^2}+C|u_\eps|\le C |u_\eps|e^{\alpha_\eps u_\eps^2}, 
\]
and \[
    \bigl(\log(2+|f_\eps|)\bigr)^{1/2}\le C(1+|u_\eps|), \ \ \text{in }\Omega\setminus B_{4R_0}(x_\eps).
\]
It follows that
\begin{equation*}\begin{split} 
&\int_{\Omega\setminus B_{4R_0}(x_\eps)}
|f_\eps|\bigl(\log(2+|f_\eps|)\bigr)^{1/2}dx
\le\frac{C}{\lambda_\eps}\int_\Omega\frac{(|u_\eps|+u_\eps^2)e^{\alpha_\eps u_\eps^2}}{|x|^\beta}dx
+C\int_\Omega(|u_\eps|+u_\eps^2)dx.
\end{split}\end{equation*}
For the first term, we have 
\begin{equation}\begin{split}\label{eq100802}
&\frac1{\lambda_\eps}\int_\Omega\frac{(|u_\eps|+u_\eps^2)e^{\alpha_\eps u_\eps^2}}{|x|^\beta}dx\\
&=\frac1{\lambda_\eps}\left(\int_{\Omega\cap \{|u_\eps|\le1\}}+\int_{\Omega\cap \{|u_\eps|>1\}}   \right)\frac{(|u_\eps|+u_\eps^2)e^{\alpha_\eps u_\eps^2}}{|x|^\beta}dx\\
&\le\frac{C}{\lambda_\eps} \int_\Omega\frac{dx}{|x|^\beta}+\frac2{\lambda_\eps}
\int_\Omega\frac{u_\eps^2e^{\alpha_\eps u_\eps^2}}{|x|^\beta}dx
\le C.
\end{split}\end{equation}
Using the uniform $L^2(\Omega)$ bound of $u_\eps$,  we immediately have
$\int_\Omega(|u_\eps|+u_\eps^2)dx\le C.$ 
From this and \eqref{eq100802}, we obtain 
\[
\int_{\Omega\setminus B_{4R_0}(x_\eps)}|f_\eps|\bigl(\log(2+|f_\eps|)\bigr)^{1/2}dx\le C.
\]
Hence, $f_\eps$ is bounded in $ L(\log L)^{\frac{1}{2}}(\Omega)$, that is
\begin{equation*} 
\int_\Omega|f_\eps|\bigl(\log(2+|f_\eps|)\bigr)^{1/2}\,dx\le C.
\end{equation*}

Now, we consider the function $u_\eps^2$. Since $u_\eps$ is a maximizer of the supremum $S_{\Omega,q,\beta,\eps}$ of subcritical singular Adams functionals, hence, for each fixed $\eps$, we have $u_\eps\in H_0^2(\Omega)\cap L^\infty(\Omega)$. Using the Sobolev embedding, we obtain
 $\nabla u_\eps\in L^4(\Omega)$ and 
$u_\eps^2\in L^2(\Omega)$,  which imply  $\nabla(u_\eps^2)=2u_\eps\nabla u_\eps\in L^2(\Omega).$
Moreover, we have
\[
\|\nabla^2(u_\eps^2)\|_{L^2(\Omega)}
\le
2\|u_\eps\|_{L^\infty(\Omega)}
  \|\nabla^2u_\eps\|_{L^2(\Omega)}
+
2\|\nabla u_\eps\|_{L^4(\Omega)}^2
<\infty.
\]
Thus $u_\eps^2\in H_0^2(\Omega)$.

Next, we claim that $\|\Delta^2(u_\eps^2)\|_1\le C$. Indeed, applying  Proposition A.10 of \cite{DG2021}, for $ j=1,2,3$, we have
\begin{equation}\label{eq:lorentz}
\|\nabla^j u_\eps\|_{L^{4/j,2}(\Omega)}\le C.
\end{equation}
By a direct computation, we have
\[
\Delta^2(u_\eps^2)=2u_\eps\Delta^2u_\eps
+2(\Delta u_\eps)^2+8\nabla u_\eps\cdot\nabla\Delta u_\eps
+4|\nabla^2u_\eps|^2.
\]
Using \eqref{eqjzfc} and the definition of $\lambda_{\eps}$, we can obtain    $\int_\Omega|u_{\eps}\Delta^2u_\eps|\le d_\eps+|\gamma_\eps|\|u_\eps\|_2^2\le C$. Since $u_\eps$ is bounded in $H^2_0(\Omega)$, the terms $2(\Delta u_\eps)^2$ and $4|\nabla^2u_\eps|^2$ are both bounded in $L^1$. Applying the H\"{o}lder inequality for Lorentz spaces and \eqref{eq:lorentz}, we get
\[
\int_\Omega|\nabla u_\eps|\,|\nabla\Delta u_\eps|\,dx
\le C\|\nabla u_\eps\|_{L^{4,2}}
       \|\nabla\Delta u_\eps\|_{L^{4/3,2}}\le C.
\]
Then the claim is proved.

Let $G(z,x)$ be the Green function of the associated with the biharmonic
operator $\Delta^2$ with Dirichlet boundary conditions at $z\in \Omega$. 
From the Proposition 2.2 of \cite{DG2021}, for $1\le j\le3$ and $ x\ne z$, we have
\[
|\nabla_x^j G(z,x)|
\le \frac{C}{|x-z|^j},
\]
where $C$ is independent of $x$ and $z$.  Moreover, for any $a,z\in\mathbb R^4$, $r>0$ and $1\le j\le3$, we have
\begin{align*}
\int_{B_r(a)}\frac{dx}{|x-z|^j}
&=\int_{B_r(a)\cap B_r(z)}\frac{dx}{|x-z|^j}
+\int_{B_r(a)\setminus B_r(z)}\frac{dx}{|x-z|^j}\\
&\le\int_{B_r(z)}\frac{dx}{|x-z|^j}+r^{-j}|B_r(a)|
=\left(\frac{|\mathbb S^3|}{4-j}+|B_1|\right)r^{4-j}.
\end{align*}
In particular, we obtain 
\begin{equation}\label{eq100805}
\sup_{z\in\mathbb R^4}\int_{B_r(a)}\frac{dx}{|x-z|^j}\le Cr^{4-j}.
\end{equation}

Now, using  the Green representation, we have
\[
|\nabla^j(u_\eps^2)(x)|\le C\int_\Omega\frac{|\Delta^2(u_\eps^2)(z)|}{|x-z|^j}dz.
\]
 Integrating over $B_r(a)\cap\Omega$, applying Tonelli's theorem, \eqref{eq100805} and the fact that $\|\Delta^2(u_\eps^2)\|_1\le C$, we obtain
\begin{align*}
\int_{B_r(a)\cap\Omega}|\nabla^j(u_\eps^2)(x)|dx
&\le C\int_{B_r(a)\cap\Omega}\int_\Omega
\frac{|\Delta^2(u_\eps^2)(z)|}{|x-z|^j}dzdx\\
&=C\int_\Omega|\Delta^2(u_\eps^2)(z)|\left(\int_{B_r(a)\cap\Omega}\frac{dx}{|x-z|^j}\right)dz\\
&\le Cr^{4-j}\|\Delta^2(u_\eps^2)\|_{L^1(\Omega)}\le Cr^{4-j},
\end{align*}
where the constant $C$ is independent of $\eps$, $a$ and $r$. The proof of this lemma is finished.
\end{proof}

\begin{lemma}\label{lemDMdensity}
There exists $C>0$ such that for any $x\in\Omega$,
\begin{equation*} 
|x-x_\eps|^{4-\beta}\frac{d_\eps}{\lambda_\eps}
 u_\eps^2e^{\alpha_\eps u_\eps^2}\le C.
\end{equation*}
\end{lemma}

\begin{proof}
Denote
\[
L_\eps:=\max_{\overline\Omega}
\frac{d_\eps}{\lambda_\eps}|x-x_\eps|^{4-\beta}
 u_\eps(x)^2e^{\alpha_\eps u_\eps(x)^2}.
\]
Assume by contradiction that $L_\eps\rightarrow \infty$ as $\eps\rightarrow 0$. For each $\eps$, observe that the function being maximized is continuous, including at $0$, and vanishes at $x_\eps$ and on $\partial\Omega$. Hence, we can choose a maximizing point $\widetilde x_\eps\in\Omega\setminus\{x_\eps\}$ 
and set
$\widetilde c_\eps=u_\eps(\widetilde x_\eps)$.  

 It follows from \eqref{eq4.4}  and \eqref{lem3.4}  that $|\widetilde c_\eps|\to\infty$.
Since $|u_\eps|\le c_\eps$, from the definition of $r_\eps$, we obtain
\[
L_\eps\le r_\eps^{-4}|\widetilde x_\eps-x_\eps|^{4-\beta}
=\left(\frac{|\widetilde x_\eps-x_\eps|}
{r_\eps^{4/(4-\beta)}}\right)^{4-\beta}.
\]
By \eqref{xstar}, this implies
\begin{equation}\label{eq:separatecore}
\frac{|\widetilde x_\eps-x_\eps|}{r_\eps^{4/(4-\beta)}}\to\infty,
\frac{|x_\eps|}{|\widetilde x_\eps-x_\eps|}\to0 \text{ and }\frac{|\widetilde x_\eps|}{|\widetilde x_\eps-x_\eps|}\to1, \ as\  \eps \to 0.
\end{equation}
In particular, one can obtain $\widetilde x_\eps\ne0$ for small $\eps$.

Define $\delta_\eps>0$ by
\begin{equation}\label{eq100806}
\delta_\eps^4\frac{d_\eps}{\lambda_\eps}
|\widetilde x_\eps|^{-\beta}\widetilde c_\eps^2
 e^{\alpha_\eps\widetilde c_\eps^2}=1.
\end{equation}
Then
\begin{equation}\label{adddelta1015}
    \frac{\delta_\eps}{|\widetilde x_\eps-x_\eps|}
=L_\eps^{-1/4}\left(\frac{|\widetilde x_\eps|}{|\widetilde x_\eps-x_\eps|}\right)^{\beta/4}\to0.
\end{equation}
Hence $\delta_\eps\to0$ and $\delta_\eps/|\widetilde x_\eps|\to0$.

Now, we define a new rescaled function by
$\widetilde u_\eps(y):=u_\eps(\widetilde x_\eps+\delta_\eps y)$  on the corresponding rescaled domain $\hat{\Omega}_\eps: = \{x\in \RR^4: \widetilde x_\eps+\delta_\eps x\in \Omega\}$. Since $\widetilde x_\eps$ is the maximizing point for $L_\eps$, we have
\begin{equation}\label{eq:maximality}
\frac{\widetilde u^2_\eps(y)}{\widetilde c_\eps^2}e^{\alpha_\eps(\widetilde u^2_\eps(y)-\widetilde c_\eps^2)}
\le\left(\frac{|\widetilde x_\eps-x_\eps|}
{|\widetilde x_\eps+\delta_\eps y-x_\eps|}\right)^{4-\beta}.
\end{equation}
On each fixed ball $B_R(0)$ intersected with the domain $\hat{\Omega}_\eps$, it follows from  \eqref{eq:separatecore} and the fact that $\delta_\eps/|\widetilde x_\eps|\to0$,   the right-hand side of \eqref{eq:maximality} converges to $1$. Moreover, 
\begin{equation}\label{eq100901}
\left(\frac{|\widetilde x_\eps|}
{|\widetilde x_\eps+\delta_\eps y|}\right)^\beta\to1 
\end{equation}uniformly.
Fix $R>0$. By \eqref{eq:maximality} and  \eqref{adddelta1015}, on  $B_R(0)\cap\hat{\Omega}_\eps$ we have
\begin{equation}\label{addbound1015}
    \frac{\widetilde u^2_\eps(y) }{\widetilde c_\eps^2}e^{\alpha_\eps(\widetilde u^2_\eps(y)-\widetilde c_\eps^2)}
\le C_R,
\end{equation}
where $C_R>1$ is independent of $\eps$.

If $|\widetilde u_\eps|\le|\widetilde c_\eps|$, then both the factors $\widetilde u^2_\eps(y) / \widetilde c_\eps^2$ and $e^{\alpha_\eps(\widetilde u^2_\eps(y)-\widetilde c_\eps^2)}$ are at most $1$. Otherwise, both factors are greater than $1$, so by \eqref{addbound1015} each factor is bounded by $C_R$. Thus, we obtain
\begin{equation}\label{eq:upperdensity}
\left|\frac{\widetilde u_\eps}{\widetilde c_\eps}\right|
\le C_R,  
\widetilde u_\eps^2-\widetilde c_\eps^2\le C_R
\end{equation}
for some $C_R$ depending only on $R$.

Next, we consider the normalized functions $\widetilde u_\eps/\widetilde c_\eps$. From \eqref{eqjzfc},  on  $B_R(0)\cap\hat{\Omega}_\eps$  we have
\begin{equation}\label{eq:normaldden}
\begin{aligned}
\Delta^2\left(\frac{\widetilde u_\eps}{\widetilde c_\eps}\right)
&=\frac1{\widetilde c_\eps^2}\frac{\widetilde u_\eps}{\widetilde c_\eps}
 e^{\alpha_\eps(\widetilde u_\eps^2-\widetilde c_\eps^2)}
\left(\frac{|\widetilde x_\eps|}
{|\widetilde x_\eps+\delta_\eps y|}\right)^\beta+\gamma_\eps\delta_\eps^4\frac{\widetilde u_\eps}{\widetilde c_\eps}.
\end{aligned}
\end{equation}
From \eqref{eq:upperdensity} and \eqref{eq100901}, we see that the right-hand side of \eqref{eq:normaldden} tends uniformly to zero on each fixed $R>0$. 
Furthermore, using the $L^4(\Omega)$ boundedness of $|\nabla u_\eps|$, we have
\begin{equation}\label{eq:smallg}
\left\|\nabla\left(\frac{\widetilde u_\eps}{\widetilde c_\eps}\right)\right\|_{L^4(\hat{\Omega}_\eps)}
\le C/|\widetilde c_\eps|\to0.
\end{equation}

We now claim that $\frac{\operatorname{dist}(\widetilde x_\eps,\partial\Omega)}{\delta_\eps}\to\infty.$ 
Otherwise, we select a subsequence,
\[
\frac{\operatorname{dist}(\widetilde x_\eps,\partial\Omega)}
{\delta_\eps}\to a\in[0,\infty).
\]
Since $\delta_\eps\to0$ and $\partial\Omega$ is smooth, the rescaled domains $\hat{\Omega}_\eps$ converge locally, after a rotation, to a half-space whose boundary has distance $a$ from the origin.
By \eqref{eq:normaldden}, \eqref{eq:upperdensity}, and the Schauder  estimates of the $\Delta^2$ (see Theorem 2.19 in \cite{Gazzola}), there exists some $\hat{\psi}\in C_{loc}^1(\hat{\Omega}_\eps)$ such that     $\widetilde u_\eps/\widetilde c_\eps\rightarrow \hat{\psi} \text{ in } C_{loc}^1(\hat{\Omega}_\eps)$, up to the boundary.
It follows from \eqref{eq:smallg} that $\nabla\hat{\psi}=0$, so $\hat{\psi}$ is constant. The zero boundary condition then forces $\hat{\psi}\equiv 0$. But this contradicts the fact that
$\widetilde u_\eps(0)/\widetilde c_\eps=1$. This proves the claim.

For any $0<\eta<1$,  using the interior $W^{4,p}$ estimates, \eqref{eq:smallg} and the fact that $\widetilde u_\eps(0)/\widetilde c_\eps=1$, we obtain 
\begin{equation}\label{eq:normalize}
\frac{\widetilde u_\eps}{\widetilde c_\eps}\to1
 \text{ in }C^{3,\eta}_{\mathrm{loc}}(\RR^4).
\end{equation}
Using the uniform  $L^4(\Omega)$ bound on $u_\eps$, we obtain from the above that $\delta_\eps^4\widetilde c_\eps^4\le C$.  
Hence 
\begin{equation}\label{eq:lineart}
\delta_\eps^4\widetilde c_\eps^2\to0, \ as \ \eps\to 0.
\end{equation}

Set $$\zeta_\eps:=\widetilde c_\eps
(\widetilde u_\eps-\widetilde c_\eps)=\frac{\widetilde u_\eps^2-\widetilde c_\eps^2}
{\widetilde u_\eps/\widetilde c_\eps+1}.$$
For each fixed $R$, \eqref{eq:normalize} implies
$\widetilde u_\eps/\widetilde c_\eps\ge1/2$ on $B_{4R}$ for small enough $\eps$.  Then,  we have
 $\zeta_\eps\le C_R$  and $\zeta_\eps(0)=0$. 
Moreover, from \eqref{eqjzfc},  $\zeta_\eps$ satisfies the following equation
\begin{equation}\label{eq:zetapde}
\begin{aligned}
\Delta^2\zeta_\eps
&=\frac{\widetilde u_\eps}{\widetilde c_\eps}
 e^{\alpha_\eps(\widetilde u_\eps^2-\widetilde c_\eps^2)}
\left(\frac{|\widetilde x_\eps|}
{|\widetilde x_\eps+\delta_\eps y|}\right)^\beta
+\gamma_\eps\delta_\eps^4\widetilde c_\eps\widetilde u_\eps.
\end{aligned}
\end{equation}
 Using \eqref{eq100901}, \eqref{eq:upperdensity}, \eqref{eq:normalize} and \eqref{eq:lineart}, we obtain \begin{equation}\label{adddetasq1015}
    \| \Delta^2\zeta_\eps\|_{L^\infty(B_{4R})}\le C_R.
 \end{equation}

On the other hand, by Lemma~\ref{lemDMsquare}, the uniform $L^4$ bound of $|\nabla u_\eps|$, and a scaling argument, we obtain
\[
\int_{B_{4R}}|\Delta(\widetilde u_\eps^2)|\,dy\le CR^2, 
\int_{B_{4R}}|\nabla\widetilde u_\eps|^2\,dy\le CR^2.
\]
This, together with the identity
\[
\Delta\zeta_\eps=
\frac{\widetilde c_\eps}{2\widetilde u_\eps}
\Delta(\widetilde u_\eps^2)
-\frac{\widetilde c_\eps}{\widetilde u_\eps}
 |\nabla\widetilde u_\eps|^2
\]
 gives $\|\Delta\zeta_\eps\|_{L^1(B_{4R})}\le CR^2$. Combining this  with \eqref{adddetasq1015} and  applying the interior elliptic estimate for the Poisson equation
$\Delta(\Delta\zeta_\eps)=\Delta^2\zeta_\eps$,   we obtain
\begin{align*}
\|\Delta\zeta_\eps\|_{L^\infty(B_{3R})}
&\le C\left(R^{-4}\|\Delta\zeta_\eps\|_{L^1(B_{4R})}
+R^2\|\Delta^2\zeta_\eps\|_{L^\infty(B_{4R})}\right)\le C_R.
\end{align*}

Let $w_\eps$ solve $\Delta w_\eps=\Delta\zeta_\eps$ in $B_{2R}$ with $w_\eps=0$ on $\partial B_{2R}$. Then
$\|w_\eps\|_\infty\le C_R$ and $\zeta_\eps-w_\eps$ is harmonic and bounded from above. Apply Harnack's inequality to a sufficiently large constant minus $\zeta_\eps-w_\eps$.
Its value at $0$ is bounded, because $\zeta_\eps(0)=0$.
Hence $\zeta_\eps$ is also bounded below on $B_R$. Combining this with \eqref{eq:zetapde} and using the interior $W^{4,p}$ estimates for $\Delta^2$, we obtain for every $0<\eta<1$,  there exists some $\zeta_0\in C^{3,\eta}_{\mathrm{loc}} (\RR^4)$ such that 
\[
\zeta_\eps\to\zeta_0
\quad\text{in }C^{3,\eta}_{\mathrm{loc}}(\RR^4),
\text{ with } \zeta_0(0)=0.
\]
 In particular, using \eqref{eq:normalize}, we have
 \[
\widetilde u_\eps^2-\widetilde c_\eps^2
=\left(1+\frac{\widetilde u_\eps}{\widetilde c_\eps}\right)
\zeta_\eps \to2\zeta_0
\quad\text{locally uniformly}.
\]
Equations \eqref{eq100806} and \eqref{eq100901} then give
\begin{align*}
\frac1{\lambda_\eps}\int_{B_{\delta_\eps}(\widetilde x_\eps)}
\frac{u_\eps^2e^{\alpha_\eps u_\eps^2}}{|x|^\beta}dx
&=\frac1{d_\eps}\int_{B_1}
\frac{\widetilde u_\eps^2}{\widetilde c_\eps^2}
 e^{\alpha_\eps(\widetilde u_\eps^2-\widetilde c_\eps^2)}
\left(\frac{|\widetilde x_\eps|}
{|\widetilde x_\eps+\delta_\eps y|}\right)^\beta dy\\
&\rightarrow\int_{B_1}e^{64\pi^2(1-\beta/4)\zeta_0(y)}dy>0.
\end{align*}

 From  \eqref{eq:separatecore} and \eqref{adddelta1015}, we see that for every fixed $R$, $B_{\delta_\eps}(\widetilde x_\eps)$ and
$B_{Rr_\eps^{4/(4-\beta)}}(x_\eps)$ are disjoint for small $\eps$. Hence, using the definition of $\lambda_\eps$ and \eqref{eq:phiMassShift}, we obtain 
\[
\begin{aligned}
1&\ge \lim_{R\to\infty}\lim_{\eps\to0}\frac1{\lambda_\eps}
\left(\int_{B_{Rr_\eps^{4/(4-\beta)}}(x_\eps)}\frac{u_\eps^2e^{\alpha_\eps u_\eps^2}}{|x|^\beta}dx
+\int_{B_{\delta_\eps}(\widetilde x_\eps)}\frac{u_\eps^2e^{\alpha_\eps u_\eps^2}}{|x|^\beta}dx\right)\\
&=1+\int_{B_1}e^{64\pi^2(1-\beta/4)\zeta_0(y)}dy>1,
\end{aligned}
\]
which is a contradiction. Therefore, we complete the proof.
\end{proof}

\begin{lemma}[Pointwise gradient estimates]\label{lemDMpoint}
For each $l\in\{1,2,3\}$, there exists a constant $C>0$ such that
\begin{equation}\label{eq:pointtwo}
\min\{|x|,|x-x_\eps|\}^{l}
|u_\eps(x)|\,|\nabla^l u_\eps(x)|\le C \ 
  \text{ for all } \ x\in\Omega\setminus\{0,x_\eps\}.
\end{equation}
In particular, for every $x\in\Omega\setminus\{0,x_\eps\}$ satisfying $|x-x_\eps|\ge2|x_\eps|$, it holds
\begin{equation}\label{eq:pointcenter}
|x-x_\eps|^{l}|u_\eps(x)|\,|\nabla^l u_\eps(x)|\le C.
\end{equation}
\end{lemma}

\begin{proof}
For each fixed $\eps$,  since $u_\eps$ is bounded,  we have
\[
|\Delta^2u_\eps(x)|\le C_\eps(1+|x|^{-\beta}).
\]
 We first claim that, for each $l\in{1,2,3}$ and all $x$ sufficiently close to the origin, one has
\begin{equation}\label{addgrad1016}
   |\nabla^l u_\eps(x)|\le C_\eps\bigl(1+|\log|x||+|x|^{4-\beta-l}\bigr).
\end{equation}
Indeed, for each $l\in\{1,2,3\}$, the Green representation formula gives

\begin{equation}\label{eq:green-derivative-bound}
|\nabla^l u_\eps(x)|
\le C_\eps\left(
1+\int_\Omega |x-y|^{-l}|y|^{-\beta}\,dy\right).
\end{equation}
Choose $\rho>0$ such that $\Omega\subset B_{\rho}(0)$. For $0<r:=|x|$, we decompose
\[
\begin{aligned}
\int_\Omega |x-y|^{-l}|y|^{-\beta}dy&=\left(\int_{B_{r/2}(0)}+\int_{B_{2r}(0)\setminus B_{r/2}(0)}+\int_{B_\rho(0)\setminus B_{2r}(0)}\right) |x-y|^{-l}|y|^{-\beta}dy\\
&=I_1+I_2+I_3+I_4.
\end{aligned}\]
On $B_{r/2}(0)$, we have $|x-y|\ge r/2$, and hence
\[
I_1
\le Cr^{-l}\int_0^{r/2}t^{3-\beta}\,dt
\le Cr^{4-\beta-l}.
\]
On $B_{2r}(0)\setminus B_{r/2}(0)$, we have
$|y|^{-\beta}\le Cr^{-\beta}$. Since
$B_{2r}(0)\subset B_{3r}(x)$, it follows that
\[
I_2
\le Cr^{-\beta}\int_{B_{3r}(x)}|x-y|^{-l}\,dy
\le Cr^{4-\beta-l}.
\]
On $B_\rho(0)\setminus B_{2r}(0)$, the inequality
$|x-y|\ge |y|/2$ gives
\[
I_3
\le C\int_{2r}^{\rho}t^{3-\beta-l}\,dt
\le
\begin{cases}
C, & \beta+l<4,\\
C(1+|\log r|), & \beta+l=4,\\
Cr^{4-\beta-l}, & \beta+l>4.
\end{cases}
\]
 Combining these estimates with \eqref{eq:green-derivative-bound}, the claim is proved.

Since $\beta<4$, using \eqref{addgrad1016}, we see that the expression in \eqref{eq:pointtwo} can be extended continuously by zero at $0$ and $x_\eps$. Hence, this extension  is continuous on the compact set $\overline{\Omega}$ and therefore attains its maximum.
 

For any fixed $r>0$, since $x_\eps\to0$, for sufficiently small $\eps$,   on $\Omega\setminus B_{r/2}(0)$  we have
\[
 |x-x_\eps|\ge r/4.
\]
On the set $(\Omega\setminus B_{r/2}(0))\cap\{|u_\eps|\ge1\}$, using Lemma~\ref{lemDMdensity}, we have
\[
\frac{d_\eps}{\lambda_\eps}\frac{|u_\eps|e^{\alpha_\eps u_\eps^2}}{|x|^\beta}\le
\frac{C}{|x|^\beta|x-x_\eps|^{4-\beta}|u_\eps|}\le C_r.
\]
While on  the set  $(\Omega\setminus B_{r/2}(0))\cap\{|u_\eps|<1\}$, we can obtain the same bound follows directly from the
uniform bounds for $d_\eps/\lambda_\eps$ and $\alpha_\eps$.
Thus, by the Sobolev embedding, the right-hand side of \eqref{eqjzfc} is uniformly
bounded in $L^p(\Omega\setminus B_{r/2}(0))$ for every finite $p>1$.
Combining the \(W^{4,p}\) estimates for $\Delta^2$ and  the Sobolev embedding, we obtain
\[
\sup_{\Omega\setminus B_r(0)}
\left(|u_\eps|+\sum_{l=1}^3|\nabla^lu_\eps|\right)\le C_r.
\]
Hence the expression in \eqref{eq:pointtwo} is uniformly bounded away from the origin.
Consequently, if the maximum in \eqref{eq:pointtwo} tends to infinity  as $\eps\to 0$, then  the corresponding maximizing points $\widetilde x_\eps$ must satisfy $\widetilde x_\eps\to0$.

For each $l=1,2,3,$ denote 
\[
L_{l,\eps}:=\max_{x\in\overline\Omega}\min\{|x|,|x-x_\eps|\}^{l}|u_\eps(x)|\,|\nabla^l u_\eps(x)|.
\]
The proof proceeds in two steps.

\smallskip
\noindent\emph{Step 1: $l=1$.} We show $L_{1,\eps}<C.$ Suppose $L_{1,\eps}\to\infty,$
choose a maximizing point $\widetilde x_\eps$  for $L_{1,\eps}$. Set
\[
\delta_\eps=\min\{|\widetilde x_\eps|,
|\widetilde x_\eps-x_\eps|\},
\widetilde c_\eps=u_\eps(\widetilde x_\eps),
\]
and define $\widetilde u_\eps(y):=u_\eps(\widetilde x_\eps+\delta_\eps y).$ 

From the discussion above and \eqref{xstar}, we see that $\widetilde x_\eps\to0$ and $\delta_\eps\to0$. Since $0\in\Omega$, we have
$B_{\delta_\eps/2}(\widetilde x_\eps)\Subset\Omega$ for sufficiently small $\eps$. Furthermore,
for $y\in B_{1/2}$, we see that both $|\widetilde x_\eps+\delta_\eps y|$ and $|\widetilde x_\eps+\delta_\eps y-x_\eps|$
are at least $\delta_\eps/2$. Then, using the fact that the maximum is attained at $\widetilde x_\eps$, we obtain
\begin{equation}\label{eq:firstmax}
|\widetilde u_\eps|\,|\nabla\widetilde u_\eps|\le2L_{1,\eps}=2|\widetilde c_\eps| |\nabla\widetilde u_\eps(0)| \text{ on }B_{1/2}.
\end{equation}
Moreover,  by Lemma~\ref{lemDMdensity}, on $B_{1/2}$ we  obtain
\begin{equation}\label{eq100905}
\begin{aligned}
\frac{d_\eps\delta_\eps^4}{\lambda_\eps}
\frac{\widetilde u_\eps^2e^{\alpha_\eps\widetilde u_\eps^2}}{|\widetilde x_\eps+\delta_\eps y|^\beta}
&\le\frac{C\delta_\eps^4}
{|\widetilde x_\eps+\delta_\eps y|^\beta|\widetilde x_\eps+\delta_\eps y-x_\eps|^{4-\beta}}
\le C.
\end{aligned}
\end{equation}
Hence, for every fixed $T>1$,
\[
\sup_{B_{1/2}}\left|\frac{d_\eps\delta_\eps^4}{\lambda_\eps}
\frac{\widetilde u_\eps e^{\alpha_\eps\widetilde u_\eps^2}}
{|\widetilde x_\eps+\delta_\eps y|^\beta}\right|
\le C\delta_\eps^{4-\beta}T e^{CT^2}+\frac CT.
\]
 Here, the first term provides the bound on ${|\widetilde u_\eps|\le T}$, while the second provides the bound on its complement from \eqref{eq100905}. Letting first $\eps\to0$ and then $T\to\infty$, we obtain
\begin{equation}\label{eq:smallrhs}
\Delta^2\widetilde u_\eps=o(1)+\gamma_\eps\delta_\eps^4\widetilde u_\eps
\quad\text{uniformly on }B_{1/2}.
\end{equation}

Set $D_\eps=\max_{0\le j\le3}|\nabla^j\widetilde u_\eps(0)|$.
Then $L_{1,\eps}\le D_\eps^2$, hence $D_\eps\to\infty$  as $\eps\to 0$.
Integration along the segment from $0$ to $y$ and using \eqref{eq:firstmax}, we have
\[
|\widetilde u^2_\eps(y)-\widetilde c_\eps^2|
\le\int_0^1 2|\widetilde u_\eps(ty)|
 |\nabla\widetilde u_\eps(ty)|\,|y|\,dt
\le4L_{1,\eps}|y|,
\]
and then $|\widetilde u_\eps|\le\sqrt3D_\eps$ on $B_{1/2}$. Hence, the functions $\widetilde u_\eps/D_\eps$ are uniformly bounded on $B_{1/2}$. Therefore, by \eqref{eq:smallrhs}, we have $\Delta^2(\widetilde u_\eps/D_\eps)\rightarrow 0$ uniformly bounded on $B_{1/2}$. Besides, from the uniform $L^4$ bound of $|\nabla u_\eps|$, and a scaling argument, we obtain
\[
\|\nabla(\widetilde u_\eps/D_\eps)\|_{L^4(B_{1/2})}
\le C/D_\eps\to0,
\]as $\eps\to 0$. Combining this and the interior $W^{4,p}$ estimates with $p>4$ for $\Delta^2$, we obtain $\widetilde u_\eps/D_\eps$ convergence to some constant $A$ in $C^3(\overline{B_{1/4}})$.
The normalization
\[
\max_{0\le j\le3}
|\nabla^j(\widetilde u_\eps/D_\eps)(0)|=1
\]
forces $|A|=1$ and  $\lim_{\eps\to 0}D_\eps/ |\widetilde c_\eps|=1$. Hence, it follows that, as $\eps\to 0$,
\begin{equation}\label{eq:ratio}
|\widetilde c_\eps|\to\infty,
\frac{\widetilde u_\eps}{\widetilde c_\eps}\to1
\text{ in }C^3(\overline{B_{1/4}}).
\end{equation}
  In particular,
$\frac12|\widetilde c_\eps|\le|\widetilde u_\eps|\le2|\widetilde c_\eps|$
{on }$B_{1/4}$ for sufficiently small $\eps$. This, together with
the uniform $L^4(\Omega)$ bound for $u_\eps$ implies
\[
C\ge
\delta_\eps^4\int_{B_{1/4}}|\widetilde u_\eps|^4\,dy
\ge C^{-1}\delta_\eps^4|\widetilde c_\eps|^4.
\]
Consequently, we get $\delta_\eps^4\widetilde c_\eps^4\le C$.  Hence
\begin{equation}\label{eq100906}
\delta_\eps^4\widetilde c_\eps^2\to0,\  \text{as} \ \eps\to 0.
\end{equation}
  Define
\[
\widetilde w_\eps=\frac{\widetilde u_\eps-\widetilde c_\eps}{|\nabla\widetilde u_\eps(0)|}.
\]
By \eqref{eq:firstmax} and \eqref{eq:ratio}, we have
$\|\nabla\widetilde w_\eps\|_{L^\infty(B_{1/4})}\le C$.
Since $\widetilde w_\eps(0)=0$, we also have $\|\widetilde w_\eps\|_{L^\infty(B_{1/4})}\le C$.
Furthermore, using \eqref{eq100905} and \eqref{eq:ratio}, we get
\begin{equation}\label{addlap1016}
  |\Delta^2\widetilde u_\eps|\le\frac{C}{|\widetilde c_\eps|}+C\delta_\eps^4|\widetilde c_\eps|
 \text{ on }B_{1/4}.  
\end{equation}
Dividing \eqref{addlap1016}  by $|\nabla\widetilde u_\eps(0)|$ and using
\eqref{eq100906}, we obtain
\[
\|\Delta^2\widetilde w_\eps\|_{L^\infty(B_{1/4})}
\le\frac{C(1+\delta_\eps^4\widetilde c_\eps^2)}{L_{1,\eps}}\to0.
\]
On the other hand, by Lemma~\ref{lemDMsquare},
\[
\int_{B_{1/8}}|\nabla(\widetilde u_\eps^2)|dy
=\delta_\eps^{-3}\int_{B_{\delta_\eps/8}(\widetilde x_\eps)}|\nabla(u_\eps^2)|dx
\le C.
\]
Hence, from \eqref{eq:ratio}, we have
\[\begin{aligned}
\int_{B_{1/8}}|\nabla\widetilde w_\eps|dy
&=\int_{B_{1/8}}\frac{|\nabla(\widetilde u_\eps^2)|}{2|\widetilde u_\eps|\,
|\nabla\widetilde u_\eps(0)|}dy\le\frac{C}{L_{1,\eps}}\to0.
\end{aligned}\]
 Combining this and the interior $W^{4,p}$ estimates for $\Delta^2$, we obtain $\widetilde w_\eps$ convergence to $0$ in $C^3(\overline{B_{1/8}})$,  contradicting $|\nabla\widetilde w_\eps(0)|=1$. This proves \eqref{eq:pointtwo} for $l=1$.

\bigskip
\noindent\emph{Step 2: $l=2,3$.}
Fix $l\in\{2,3\}$ and suppose
$L_{l,\eps}\to\infty.$   With a slight abuse of notation, as in Step~1, we choose a maximizing
point and define $\widetilde x_\eps$, $\delta_\eps$, $\widetilde c_\eps$,
$\widetilde u_\eps$, and $D_\eps$ such that
\[
\widetilde x_\eps\to0, \delta_\eps\to0,
\]
and
\[
L_{l,\eps}=|\widetilde c_\eps|\,|\nabla^l\widetilde u_\eps(0)|
\le D_\eps^2.
\]
In particular, $D_\eps\to\infty$.
From the Step 1 and \eqref{eq:firstmax}, we see that
\begin{equation}\label{addbound1016}
   |\widetilde u_\eps|\,|\nabla\widetilde u_\eps|
\le C
 \text{ on }B_{1/2}. 
\end{equation}
Integrating along a segment, and by \eqref{addbound1016} we obtain
\begin{equation}\label{eq100807}
|\widetilde u^2_\eps(y)-\widetilde c_\eps^2|
\le C, \text{ for } y\in B_{1/2}.
\end{equation}
Also, one can easily check that the estimates \eqref{eq100905} and \eqref{eq:smallrhs} remain valid in this case.
Since \eqref{eq100807} implies $|\widetilde u_\eps|\le CD_\eps$ on $B_{1/2}$,
the same normalization argument as in Step~1 implies that  \eqref{eq:ratio} and \eqref{eq100906} still hold.

Set
\[
\widetilde w_\eps
=\frac{\widetilde u_\eps-\widetilde c_\eps}{|\nabla^l\widetilde u_\eps(0)|}.
\]
Combining \eqref{addbound1016} with \eqref{eq:ratio}, we have
\begin{equation}\label{addboundgr1016}
    \|\nabla\widetilde w_\eps\|_{L^\infty(B_{1/4})}+\|\widetilde w_\eps\|_{L^\infty(B_{1/4})}\le\frac{C}{L_{l,\eps}},
\end{equation}
where we have used the condition $\widetilde w_\eps(0)=0$.
As in Step~1, we also have
\[
\|\Delta^2\widetilde w_\eps\|_{L^\infty(B_{1/4})}
\le\frac{C(1+\delta_\eps^4\widetilde c_\eps^2)}{L_{l,\eps}}\to0.
\]Combining the interior $W^{4,p}$ estimates with $p>4$ for $\Delta^2$ and \eqref{addboundgr1016}, we obtain
$\widetilde w_\eps\to0  \text{ in }C^3(\overline{B_{1/8}}),$
contradicting $|\nabla^l\widetilde w_\eps(0)|=1$.
This proves \eqref{eq:pointtwo} for $l=2,3$.

Finally, if $|x-x_\eps|\ge2|x_\eps|$, then
\[
|x|\ge|x-x_\eps|-|x_\eps|\ge\tfrac12|x-x_\eps|,
 \min\{|x|,|x-x_\eps|\}\ge\tfrac12|x-x_\eps|.
\]
Thus the estimate \eqref{eq:pointcenter} follows. 

\end{proof}

Using the pointwise gradient estimates established above, we can now complete the proof of Proposition~\ref{lemDM1}.
\begin{proof}[Proof of  Proposition~\ref{lemDM1}]
 For each \(\theta\in\mathbb S^3\), let
\[
t_\eps^*(\theta)
=
\sup\bigl\{
t>0:x_\eps+s\theta\in\Omega
\text{ for every }0\leq s<t
\bigr\}.
\]
Since
\[
u_\eps(x_\eps)=c_\eps>\frac{c_\eps}{M},
\qquad
u_\eps\bigl(x_\eps+t_\eps^*(\theta)\theta\bigr)=0,
\]
we may define
\[
t_\eps^M(\theta)
=
\inf\left\{
t\in(0,t_\eps^*(\theta)):
u_\eps(x_\eps+t\theta)=\frac{c_\eps}{M}
\right\}.
\]
We claim that the map
\(\theta\mapsto t_\eps^M(\theta)\) is lower semicontinuous.
Indeed, fix any $a$ with $0<a<t_\eps^M(\theta)$. Since $u_\eps-c_\eps/M$ is continuous and strictly positive on the compact segment $\{x_\eps+t\theta:0\le t\le a\}$, there exists $\delta>0$ such that
\[
u_\eps(x_\eps+t\theta')-\frac{c_\eps}{M}>0
\text{ for all }t\in[0,a]
\]
whenever $|\theta'-\theta|<\delta$. Consequently, $t_\eps^M(\theta')>a$ for all such directions $\theta'$. Since $a<t_\eps^M(\theta)$ was arbitrary, we conclude that
\[
\liminf_{\theta'\to\theta}t_\eps^M(\theta')
\geq t_\eps^M(\theta),
\]
and the claim is proved. Therefore, there exists \(\theta_\eps^M\in\mathbb S^3\) such that
\[
\rho_\eps^M:=t_\eps^M(\theta_\eps^M)
=\min_{\theta\in\mathbb S^3}t_\eps^M(\theta)>0.
\]
Set $y_\eps^M=x_\eps+\rho_\eps^M\theta_\eps^M.$
Then $u_\eps(y_\eps^M)=\frac{c_\eps}{M}.$  Choose a nearest boundary point $z_\eps\in\partial\Omega$
to $x_\eps$. Using the continuity of $u_\eps$ and the fact that $u_\eps(z_\eps)=0$, we have
\[
\rho_\eps^M\le t_\eps^M\left(\frac{z_\eps-x_\eps}{|z_\eps-x_\eps|}\right)
<|z_\eps-x_\eps|=\operatorname{dist}(x_\eps,\partial\Omega).
\]
Thus $B_{\rho_\eps^M}(x_\eps)\Subset\Omega$, and 
\begin{equation}\label{point1}
    u_\eps>\frac{c_\eps}{M} \text{ in }B_{\rho_\eps^M}(x_\eps),
u_\eps\ge\frac{c_\eps}{M}
 \text{ on }\partial B_{\rho_\eps^M}(x_\eps).
\end{equation}

From Lemma~\ref{lemvk}, for every fixed $R>0$ we have
\[
u_\eps=c_\eps+o_\eps(1) \text{ uniformly in }B_{Rr_\eps^{4/(4-\beta)}}(x_\eps).
\]
Since $M>1$ and $c_\eps\to\infty$, we see that $u_\eps>c_\eps/M$ in $B_{Rr_\eps^{4/(4-\beta)}}(x_\eps)$ for sufficiently small $\eps$.
Therefore $\rho_\eps^M\ge Rr_\eps^{4/(4-\beta)}$.
By the arbitrary of $R>0$, we obtain \begin{equation}\label{point41}\frac{\rho_\eps^M}{r_\eps^{4/(4-\beta)}}\to\infty, \ \text{as} \  \eps\to 0. \end{equation} 
Since
\[
\frac{c_\eps^2}{M^2}|B_1|(\rho_\eps^M)^4\le\int_\Omega u_\eps^2 dx\le C,
\]
we immediately get
\begin{equation}\label{point42}\rho_\eps^M\to0.\end{equation}  
Furthermore, by using \eqref{point41} and \eqref{xstar}, we obtain $\frac{|x_\eps|}{\rho_\eps^M}$ as $\eps\to 0$. Hence, for sufficiently small $\eps$,
\[
|x-x_\eps|=\rho_\eps^M\ge2|x_\eps|, \text{ on } \partial B_{\rho_\eps^M}(x_\eps).
\]

By \eqref{eq:pointcenter} in Lemma \ref{lemDMpoint} and the fact that  $u_\eps\ge c_\eps/M$ on this sphere, this implies, for each $ l=1,2,3,$

\begin{equation}\label{point3}
    |\nabla^l u_\eps|
\le \frac{CM}{c_\eps(\rho_\eps^M)^l}, 
\text{ on }\partial B_{\rho_\eps^M}(x_\eps).
\end{equation}

Finally, for any $x\in\partial B_{\rho_\eps^M}(x_\eps)$, join $y_\eps^M$ to $x$ by a shortest geodesic arc $\Gamma$ on the sphere. Since $u_\eps(y_\eps^M)=c_\eps/M$ and $\operatorname{length}(\Gamma)\le\pi\rho_\eps^M$, we obtain
\begin{equation}\label{point2}\begin{aligned}
\left|u_\eps(x)-\frac{c_\eps}{M}\right|
&=|u_\eps(x)-u_\eps(y_\eps^M)|\leq\int_\Gamma|\nabla u_\eps|ds\leq
\pi\rho_\eps^M\sup_{\partial B_{\rho_\eps^M}(x_\eps)}|\nabla u_\eps|
\leq\frac{C(M)}{c_\eps}.
\end{aligned}\end{equation}

Combining \eqref{point1}  \eqref{point2} \eqref{point3} \eqref{point41} and \eqref{point42}, we complete the proof.
\end{proof}



\vspace{12pt}


\noindent College of Science, Beijing Forestry University,  Beijing, 100083, P. R. China\\
\textit{Email address:}   xmwang@bjfu.edu.cn

\vspace{8pt}
\noindent School of Mathematics and Statistics, Nanjing University of Science and Technology, Nanjing, 210094, P. R. China.\\
\textit{Email address:}  zhumaochun2006@126.com

\end{document}